\input amstex
\documentstyle{amsppt}
\magnification=1200
\define\Z{{\Bbb Z}}
\define\Q{{\Bbb Q}}
\define\R{{\Bbb R}}
\define\G {{\Bbb G}}
\define\bS{\bold{S}}
\define\PP{{\Bbb P}}
\define\quat{{\Bbb H}}
\define\C{{\Bbb C}}
\define\sli{\frak{sl}}
\define\so{\frak{so}}
\define\sy{\frak{sp}}
\define\sku{\frak{sku}}
\define\su{\frak{su}}
\define\slt{\frak{sl} (2)}
\define\tot{\text{tot}}
\define\la{\langle}
\define\ra{\rangle}
\define\g{{\frak g}}
\define\a{{\frak a}}
\define\uf{{\frak u}}
\define\s{{\frak s}}
\define\gl{\frak{gl}}

\define\p{{\frak p}}
\define\q{{\frak q}}
\define\h{{\frak h}}

\define\Mb{M_{\bullet}}
\define\hor{\text{hor}}
\define\ver{\text{ver}}
\define\odd{\text{odd}}
\define\ev{\text{ev}}

\define\Hom{\operatorname{Hom}}
\define\depth{\operatorname{depth}}
\define\Sym{\operatorname{Sym}}
\define\sym{\operatorname{sym}}
\define\GL{\operatorname{GL}}
\define\End{\operatorname{End}}
\define\Ad{\operatorname{Ad}}
\define\Prim{\operatorname{Prim}}
\define\ad{\operatorname{ad}}
\define\Tr{\operatorname{Tr}}

\define\aut{\operatorname{\frak{aut}}}
\define\Gr{\operatorname{Gr}}
\define\Ker{\operatorname{Ker}}
\define\NS{\operatorname{NS}}
\define\MT{\operatorname{MT}}
\define\Nm{\operatorname{Nm}}
\define\Hm{\operatorname{H}}
\define\IH{\operatorname{IH}}
\define\Spec{\operatorname{Spec}}
\define\diag{\operatorname{diag}}
\define\Hdg{\operatorname{Hdg}}

\define\refer#1{(#1)}
\newcount\headnumber
\newcount\labelnumber
\define\section{\global\advance\headnumber
by1\global\labelnumber=0{{\the\headnumber}.\ }}
\define\label{(\global\advance\labelnumber by1 \the\headnumber
.\the\labelnumber )\enspace}
\NoBlackBoxes

\document

\topmatter
\title
A Lie algebra attached to a projective variety
\endtitle
\rightheadtext{Lie algebra attached to a projective variety}
\leftheadtext{Eduard Looijenga and Valery Lunts}

\author
Eduard Looijenga and Valery A.\ Lunts$^*$
\endauthor

\address
Faculteit Wiskunde en Informatica, Rijksuniversiteit Utrecht,
PO Box 80.010, 3508 TA Utrecht, The Netherlands
\endaddress
\email
looijeng\@math.ruu.nl
\endemail

\address
Department of Mathematics, Indiana University, Bloomington, IN 47405, USA
\endaddress
\email
vlunts\@ucs.indiana.edu
\endemail

\toc\nofrills
{Table of Contents}
{\eightpoint
\head \;\; Introduction\hfill 1\endhead
\head 1. Lefschetz modules\hfill 4\endhead
\head 2. Jordan--Lefschetz modules\hfill 12\endhead
\head 3. Geometric examples of Jordan type I: complex tori\hfill 18\endhead
\head 4. Geometric examples of Jordan type II: hyperk\"ahlerian manifolds\hfill
24\endhead
\head 5. Filtered Lefschetz modules\hfill 30\endhead
\head 6. Frobenius--Lefschetz modules\hfill 33\endhead
\head 7. Appendix: a property of the orthogonal and symplectic Lie
algebra's\hfill\; 41\endhead}
\endtoc

\thanks
$^*$Supported by the National Science Foundation.
\endthanks


\abstract
Each choice of a K\"ahler class on a compact complex manifold
defines an action of the Lie algebra $\slt$ on its total complex cohomology. If
a nonempty set of such  K\"ahler classes is given, then we prove that the
corresponding $\slt$-copies generate a semisimple Lie algebra. We investigate
the formal properties of the resulting representation and we work things out
explicitly in the case of complex tori, hyperk\"ahler manifolds and flag
varieties. We pay special attention to the cases where this leads to a Jordan
algebra structure or a graded Frobenius algebra.
\endabstract

\endtopmatter

\head
Introduction
\endhead

\noindent
Let $X$ be a projective manifold of complex dimension $n$. If $\kappa\in
\Hm ^2(X)$ is an ample class, then cupping with it defines an operator
$e_{\kappa}$ in the total complex cohomology (denoted here by $\Hm (X)$) of
degree $2$ and the hard Lefschetz theorem asserts that for $s=0,\dots ,n$,
$e_{\kappa}^s$ maps $\Hm ^{n-s}(X)$ isomorphically onto $\Hm ^{n+s}(X)$. As is
well-known, this is equivalent to the exististence of a (unique)
operator $f_{\kappa}$ in $\Hm (X)$ of degree $-2$ such that the commutator
$[e_{\kappa},f_{\kappa}]$ is the operator $h$ which on $\Hm ^k(X)$ is
multiplication by $k-n$. The elements $e_{\kappa},h,f_{\kappa}$ make up a Lie
subalgebra $\g _{\kappa}$ of $\gl (\Hm (X))$ isomorphic to $\sli (2)$ and the
decomposition of $\Hm (X)$ as a $\g _{\kappa}$-module into isotypical summands
is just the primitive decomposition: the primitive cohomology in degree $n-s$
generates the isotypical summand associated to the irreducible representation
of
dimension $s+1$. As these operators respect the Hodge decomposition (in the
sense that $e_{\kappa}$ resp.\ $f_{\kappa}$ has bidegree $(1,1)$ resp.\
$(-1,-1)$), the Hodge structure on $\Hm (X)$ is entirely determined
by the Hodge structure on the primitive cohomology.
However, the primitive decomposition usually depends in a
nontrivial way on the choice of $\kappa$. This we regard as a fortunate fact,
as it often leads to finding an even smaller Hodge substructure of
$\Hm (X)$ that determines the one on $\Hm (X)$. To be explicit, let us define
 the {\it
N\'eron--Severi Lie algebra}  $\g _{NS}(X)$ as the Lie subalgebra of  $\gl (\Hm
(X))$ generated by the $\g _{\kappa}$'s with $\kappa$ an ample class. This Lie
algebra is defined over $\Q$ and is evenly graded by the adjoint action of its
semisimple element $h$ (with its degree $2k$ summand acting as transformations
of bidegree $(k,k)$). We prove in this paper that it is also semisimple.  So if
we regard $\Hm (X)$ as a representation of this Lie algebra, then the subspace
of $\Hm (X)$ annihilated  by the negative degree part of $\g _{NS}(X)$ is a
Hodge substructure that determines the one on $\Hm (X)$. Notice that this Hodge
substructure is itself still invariant under the degree zero part of $\g
_{NS}(X)$ (which is a reductive Lie subalgebra). Despite its naturality, this
idea appears to be new (although a note by \cite{Verbitsky 1990}, of which we
were not aware of when we started this research, is suggestive in this repect.)

Whereas the $e_{\kappa}$'s commute, the corresponding $f_{\kappa}$'s don't
in general. This makes it difficult to compute the N\'eron--Severi Lie algebra
in practice. It is often helpful when we know of a morphism from $X$ to another
projective manifold $Y$ whose base and fibers are well-understood: for example,
the fact that the associated Leray spectral sequence degenerates yields (among
other things) the existence of a copy of $\g _{NS}(Y)$ in $\g _{NS}(X)$. This
is
an ingredient of our proof that the N\'eron--Severi Lie algebra of flag variety
of a simple complex Lie group  is ``as big as possible'' (reflected by the fact
that its Hodge structure is as simple as possible): it is the Lie algebra of
infinitesimal automorphisms of a naturally defined  bilinear form  (which is
either symmetric or skew) on its cohomology.

But if the $f_{\kappa}$'s happen to commute, then we are in a very interesting
situation: the  N\'eron--Severi Lie algebra has degrees $-2$, $0$ and $2$ only
and the (complexified) N\'eron--Severi group  acquires the structure of a
Jordan
algebra without preferred unit element. For abelian varieties this is a
classical fact, although, as far as we know, it had not been seen from this
point of view. The N\'eron--Severi Lie algebra appears here as a natural
companion of the Mumford--Tate group: the latter helps us to find the Hodge
ring
as a ring of invariants, whereas the decomposition of the Hodge ring into
$\Q$-irreducible representations of the N\'eron--Severi Lie  algebra helps us
to
say more about its structure. (For example, the subring generated by the
divisor
classes is one such irreducible summand.)

Here are some variants of this construction: instead of working with complex
projective manifolds, we could do this for compact complex manifolds that admit
a K\"ahler metric and replace the ample classes by K\"ahler classes. Or we
could
even take all cohomology classes of degree $2$ that have the Lefschetz
property;
clearly, the complex structure has now become irrelevant. The resulting Lie
algebra's are again semisimple and we call them the {\it K\"ahler Lie algebra}
and the {\it total Lie algebra} of the manifold respectively.  Examples of
interest here are complex tori and hyperk\"ahler manifolds; in both cases we
get
Jordan algebra structures. In a different direction, we can take for $X$ a
projective variety and take instead its complex intersection homology, even
with
values in a variation of polarized Hodge structure.

These examples lead us to formalize the situation by means of what we have
called a {\it Lefschetz module}. This is essentially a graded vector space
equipped with a set of commuting degree two operators that have
the Lefschetz property, such that the Lie algebra generated by the
corresponding
$\slt$-triples is semisimple. So this vector space becomes a
representation of a semisimple Lie algebra, and it was one of our goals to
classify the representations that so arise. Although we found some rather
restrictive properties, we did not succeed in this.

\smallskip
We now briefly decribe the contents of the separate sections.

In section $1$
we introduce the notion that is central to this paper, that of a Lefschetz
module, and discuss its basic properties. If a Lefschetz module has a
compatible Hodge structure, as is the case for the cohomology of a projective
manifold, then there is also defined its Mumford--Tate group and we compare the
two notions. We next define and discuss the closely related notion of a
Lefschetz pair. This is followed by a partial classification of such pairs in
case the associated Lie algebra is of classical type.

In section $2$ we concentrate on the case when the $f_{\kappa}$'s commute. We
show that the resulting structure is essentially that of a Jordan algebra and
that is why a complete classification is available. We are also led to a
remarkable class of Frobenius algebra's associated to each Jordan algebra, some
of which we describe explicitly.

The next two sections are devoted to examples of K\"ahler manifolds that
give rise to Lefschetz modules of Jordan type.
First we compute the total Lie algebra and the K\"ahler Lie algebra of a
complex
torus. Then we turn our attention to the N\'eron--Severi Lie algebra of an
abelian variety $A$ and express it in terms of the endomorphism algebra of
$A$.
We find that that this N\'eron--Severi Lie algebra intersects $\End
(A)\otimes\C$ in a Lie ideal of $\End (A)\otimes\C$ and we describe the
complementary ideal.

Our treatment of hyperk\"ahler manifolds (in section 4) follows essentially
\cite{Verbitsky 1995}, a preprint that in turn is partly based on a
preliminary version of the present paper. As an application we show how the
Hodge structure on the cohomology algebra of a compact hyperk\"ahlerian
manifold
is expressed in terms of the Hodge structure on its degree two part. We also
give an alternative description of the Beauville-Bogomolov quadratic form on
the
N\'eron--Severi group.

Thus the abelian varieties and the hyperk\"ahler manifolds produce the
classical Jordan algebra's. The exceptional Jordan algebra can be realized
topologically and we ask whether it is realizable by a Calabi-Yau threefold.

Section $5$ is about filtered Lefschetz modules. The example to keep in mind
here
is the Leray filtration on $\Hm (X)$ defined by a surjective morphism $f:X\to
Y$
of projective manifolds. We apply this to the case where $f$ is a projective
space bundle. In combination with a theorem proved in the appendix we are then
able to determine the N\'eron--Severi Lie algebra of a flag variety. It would
be
interesting to do the same for the intersection homology of Schubert
varieties.

In section $6$ we investigate  another interesting class of Lefschetz modules,
which we have called {\it Frobenius--Lefschetz} modules. These  arise as the
Lefschetz submodule of the cohomology of a projective manifold generated by its
unit element. The Jordan--Lefschetz modules are among them  and we suspect that
the remaining simple Frobenius--Lefschetz modules are  ``tautological
representations'' of orthogonal or symplectic Lie algebra's. The main result
\refer{6.8} of this section supports this belief: it says that any other simple
Frobenius--Lefschetz module must be a representation of an exceptional Lie
algebra.

\smallskip
We began this work in the Fall of 1990, when both of us were at the
University of Michigan in Ann Arbor. We would like to thank its Mathematics
Department for providing so stimulating working conditions. One of us
(Looijenga) thanks in particular Igor Dolgachev for many
inspiring discussions (then and later) regarding the subject matter of this
paper as well as closely related questions. Our work continued during the
Spring of
1991, while Looijenga was at the University of Utah. He gratefully remembers
the friendly atmosphere he encountered there. After an interruption we
resumed work on this paper in the academic year 1994-95. We thank Tonny
Springer for some helpful references and Bertram Kostant, Misha
Verbitsky and Yuri Zahrin for useful conversations.

\head
\section Lefschetz modules
\endhead

\noindent\label
We fix a field $K$ of characteristic zero. Let $\Mb$ be a $\Z$-graded
$K$-vector space of finite dimension and denote by $h:M\to M$ the
transformation
that is multiplication
by $k$ in degree $k$. So a linear transformation $u:M\to M$ has degree
$k$ if and only if $[h,u]=ku$.

We say that a linear transformation $e:M\to M$ of degree $2$ has the {\it
Lefschetz property} if for all integers $k\ge 0$, $e^k$ maps
$M_{-k}$ isomorphically onto $M_k$. According to the Jacobson--Morozov lemma
this is equivalent to the existence of  $K$-linear transformation $f$ in $M$ of
degree $-2$  such that $[e,f]=h$. This $f$ is then unique and $(e,h,f)$ is a
$\sli (2)$-triple: the assignment
$$
\pmatrix
0 & 1\\
0 & 0
\endpmatrix
\mapsto e,\quad
\pmatrix
1 & 0\\
0 & -1
\endpmatrix
\mapsto h,\quad
\pmatrix
0 & 0\\
1 & 0
\endpmatrix
\mapsto f
$$
defines a representation of $\sli (2)$. If $h$ and $e$ happen to be contained
in a semisimple Lie subalgebra $\g\subset \gl (M)$, then so is $f$.

Now let $\a$ be a finite dimensional $K$-vector space. We regard $\a$ as
a graded abelian Lie algebra which is homogeneous of degree two. We say that a
graded Lie homomorphism $e :\a\to \gl (M)$ has the Lefschetz property
if for some $a\in\a$, $e_a$ has that property. Notice that the set of $a\in\a$
with the Lefschetz property is always Zariski open in $\a$. For $a$ in this
open set, we have defined the operator $f_a$ such that
$(e_a,h,f_a)$ is $\sli (2)$-triple. This defines a rational map $f:\a\to\gl
(M)$ in
the sense of algebraic geometry. We let $\g (\a ,M)$ denote the Lie
subalgebra of $\gl (M)$ generated by the transformations $e_a,f_a$. If $\a$ is
merely an abelian group that acts on $M$ by operators of degree $2$, then the
linear extension $\a\otimes K\to\gl (M)$ is a Lie homomorphism and we then
often
write $\g (\a ,M)$ for $\g (\a\otimes K ,M)$. The following example shows that
this Lie algebra need not act reductively in $M$.

\smallskip
{\it Example.} Consider the graded $\slt$-representation $M=\slt \oplus K^2$,
where $\slt =Ke+Kh+Kf$ has the adjoint representation (with its usual
grading) and $K^2$ is the trivial representation in degree zero. Define an
operator $e'$ of degree $2$ in $M$ by $e'(xf+yh+zf,u,v)=(ve,z,0)$.
Then $ee'=e'e=0$, so that $e$ and $e'$  span an abelian Lie algebra $\a$. Since
$e$ has the Lefschetz property, the Lie algebra  $\g=\g (\a ,M)$ is defined.
Now $\a$ kills $(0,1,0)$, but the line spanned by this vector has no
$\g$-invariant complement. This example was chosen as to make $\g$
infinitesimally preserve a nondegenerate quadratic form on $M$ (namely
$(xf+yh+zf,u,v)\mapsto -2xz+y^2+2uv$). A smaller example without that property
is the submodule $\slt\oplus K\oplus 0$.

\smallskip
Notice that $\g (\a ,M)$ is evenly graded and that the grading is induced
from the action of $\ad _h$. We say that $(\a ,M)$ is a {\it Lefschetz module}
if $\g (\a ,M)$ is semisimple. In case $M\not=0$, we call
greatest integer $n$ with $M_n\not= 0$ (or equivalently, $M_{-n}\not= 0$) the
{\it depth} of $M$.

The collection
of Lefschetz modules is closed under direct sums, tensor products and taking
duals. Also, a Lefschetz module $M$ has always the decomposition $M
=M_{\ev}\oplus M_{\odd}$, where  $M_{\ev}$ (resp.\ $M_{\odd}$) is the direct
sum
of the $M_k$'s with $k$ even (resp.\ odd). Since any representation of a
semisimple Lie algebra is reductive, the category of Lefschetz modules of $\a$
is semisimple. A Lefschetz $\a$-module $M$ is irreducible as a
Lefschetz module if and only if it is irreducible as a  $\g (\a ,M)$-module.

There is also an exterior direct sum and tensor product: if $(\a ',M')$
and $(\a '',M'')$ are  Lefschetz modules, then we have defined Lefschetz
modules
$$
\align
(\a '\times\a '',M'\boxplus M''),&\quad
e_{(a',a'')}(m',m'')=(e_{a'}m',e_{a''}m'')\\
(\a '\times\a '',M'\boxtimes M''),&\quad
e_{(a',a'')}(m'\otimes
m'')=e_{a'}m'\otimes m'' + m'\otimes e_{a''}m''.
\endalign
$$
The associated Lie algebra is in the first case equal to  $\g (\a ',M')\times\g
(\a '',M'')$. This is also true in the second case if both factors are
nonzero.

The preceding discussion showed that when studying Lefschetz modules we may
restrict ourselves  to irreducible ones. The following lemma allows the further
reduction  of having the associated Lie algebra {\it simple}.

\proclaim{\label Lemma}
Let $M$ be an irreducible Lefschetz $\a$-module and let $\g (\a ,M)=\g
'\times\g ''$ be a decomposition of Lie algebra's. Then this decomposition is
graded and there exist irreducible
Lefschetz $\a$-modules $M'$ and $M''$ such that $M\cong M'\otimes M''$
as Lefschetz $\a$-modules with $\g '$ resp.\ $\g ''$ corresponding to $\g (\a
,M')$ resp.\ $\g (\a ,M'')$.
\endproclaim
\demo{Proof}
Since the grading of $\g$ is the eigen space decomposition of $\ad
_h$ it is immediate that upon writing $h=(h',h'')\in \g'\times\g''$, $\g
^{(i)}$
gets a grading from $\ad _{h^{(i)}}$ making the decomposition a graded one.

Since $M$ is an irreducible module of the semisimple Lie algebra $\g (\a ,M)$,
it must have the form $M'\otimes M''$ with $M^{(i)}$ a $\g ^{(i)}$-module. This
is compatible with the gradings. If the rational
map $f:\a\to\g _{-2}=\g '_{-2}\oplus\g ''_{-2}$ is written $(f',f'')$, then
$[h,f]=-2f$ implies $[h' ,f']=-2f'$ and $[h'' ,f'']=-2f''$. So $M^{(i)}$ a
Lefschetz module of $\a$ with the stated property.
\enddemo

\label Given a Lefschetz module $M$ of $\a$, then an invariant
bilinear form on $M$ is a bilinear map $\phi :M\times M\to K$ that defines a
morphism of Lefschetz modules $M\otimes M\to K$ (where $\a$ acts trivially on
$K$): so $\phi$ is zero on $M_k\times M_l$ unless $k+l=0$ and $\a$ preserves
the form $\phi$ infinitesimally: $\phi (e_a m,m')+\phi (m,e_a m')=0$ for all
$m,m'\in M$ and $a\in\a$. If $a$ is a Lefschetz element, then the
Jacobson--Morozov lemma implies that $f_a$ also preserves $\phi$
infinitesimally. So $\g (\a ,M)$ is
then a subalgebra of $\aut (M,\phi)$.  If $\phi$ is nondegenerate and symmetric
(resp.\ skew-symmetric), then we call $(M, \phi)$ an {\it orthogonal} (resp.\
{\it symplectic}) representation. Since a nonzero invariant bilinear form on an
irreducible representation is either orthogonal or symplectic, any
Lefschetz module  with nondegenerate bilinear form is the perpendicular direct
sum of Lefschetz modules that are irreducible orthogonal, irreducible
symplectic, or the direct sum of an irreducible Lefschetz module with its
dual.

\smallskip
\label Many Lefschetz modules have the additional structure of an
algebra. Let $A=\oplus _{i=0}^{2n} A_i$ be a graded-commutative algebra with
$A_0=K$. We say that $A$ is a {\it Lefschetz algebra of depth $n$} if $A[n]$ is
a
Lefschetz module of depth $n$ over $A_2$. Such a Lefschetz module can be
endowed
with an invariant $(-)^n$-symmetric bilinear form: let $\int :A\to K$ be a
linear
form that is an isomorphism in degree $2n$ and zero in all other degrees and
define $\phi (a, b):=(-1)^q\int (ab)$ if $a$ is homogeneous of degree $n+2q$ or
$n+2q+1$. If this form is nondegenerate (which is for instance the case when
$A[n]$ is irreducible as a Lefschetz module), then the form
$(a,b)\mapsto\int (ab)$ is also nondegenerate and so $A$ becomes a Frobenius
algebra (in the graded sense).

\medskip\label
Let $M$ be a graded real vector space. A {\it Hodge structure of total weight
$d$} on $M$ consists of a  bigrading on its
complexification: $M\otimes\C =\oplus _{p,q\in\Z} M^{p,q}$ such that
(i) $M_k\otimes\C=\oplus _{p+q=k+d} M^{p,q}$ for all $k$ and
(ii) complex conjugation interchanges $M^{p,q}$ and $M^{q,p}$.
These data are conveniently described in terms of an action of the
{\it Deligne torus} on $M$. We recall \cite{Deligne 1979} that this is
two-dimensional torus $\bS$ defined over $\Q$ that is obtained from $\GL (1)$
by
restricting scalars from $\C$ to $\R$. It comes with two characters $z$, $\bar
z$
that are each others complex conjugate and generate the character group. Their
product is for obvious reasons called the {\it norm character} and is denoted
$\Nm$. There is also a natural homomorphism $w:\GL (1)\to \bS$ which on the
real
points is given by the inclusion $\R ^{\times}\subset\C ^{\times}=\bS (\R)$. We
follow Deligne's convention by letting $\bS$ act on $M$ so that $M^{p,q}$
becomes the eigen space of $z^{-p}\bar z^{-q}$ (this action is defined over
$\R$). A
positive number $t>0$, viewed as an element of $\C ^{\times}=\bS (\R)$, acts on
$M_k$  as multiplication by $t^{-k-d}$. So $w(t)$ acts on $M$ as
$t^{-d-h}(=t^{-d}\exp(-\log (t)h))$.

The action of $\sqrt{-1}\in\C ^{\times}=\bS (\R)$ on $M_{\C}$ is on $M^{p,q}$
multiplication by $(\sqrt{-1})^{q-p}$; it is a real
operator, called the {\it Weil operator}; we denote it by $J$.

Suppose we are also given a nondegenerate $(-)^d$-symmetric form $\phi :M\times
M\to\R$ that is zero on $M^{p,q}\times M^{p',q'}$ unless $(p+p',q+q')=(d,d)$
(this
is equivalent to $\phi (gm,gm')=\Nm (g)^{-d}\phi (m,m')$ for all $g\in\bS (\C
)$).
Let  $e:M\to M$ be a real operator of bidegree $(1,1)$ which preserves $\phi$
infinitesimally. Clearly, $e$ commutes with $J$
and it is easily checked that for $k\ge 0$,
the map $H^k_e :M_{-k}\otimes\C\times M_{-k}\otimes\C\to\C$ defined by
$$
H^k_e(m,m'):=\phi (e^k m,J\overline{m'}),
$$
is a Hermitian form. We say that $e$ is a {\it polarization} of $(M,\phi)$
if for all $k\ge 0$, $H^k_e$ is definite on $\Ker (e^{k+1}|M_{-k}\otimes \C)$.

\proclaim{\label Proposition}
Let $(M,\phi)$ be as above. Let $\a$ be a real abelian Lie algebra with a pure
weight two Hodge structure that acts morphically on $(M,\phi )$  (i.e., the
action is by mutually commuting transformations of degree $2$ that preserve
$\phi$
infinitesimally and with $\a\otimes M\to M$ a morphism of Hodge structures).
Assume that for some $a\in\a ^{1,1}$, $e_a$  polarizes $(M ,\phi)$.
Then $M$ is a Lefschetz module of $\a$ and $\g (\a, M)$ is a semisimple Lie
algebra defined over $\R$ that preserves $\phi$ infinitesimally.
\endproclaim
\demo{Proof}
The nondegenerateness of the Hermitian form $H^k_a$ on $\Ker
(e_a^{k+1}| M_{-k}\otimes\C )$ implies that $e_a^k$ is injective on this
subspace.
It is easily checked that this, together with the nondegenerateness of $\phi$,
implies that $e_a$ satisfies the Lefschetz property. So $\g (\a ,M)$ is
defined.
If we regard $\phi$ as an element of $M^*\otimes M^*$ of degree zero, then the
fact that $\phi$ is killed by $e_a$ implies that it is killed by $f_a$. So $\g
(\a, M)$ preserves $\phi$ infinitesimally.

We next show that $\g (\a ,M)$ is reductive; since $\g (\a ,M)$ is generated by
commutators, it then follows that $\g (\a ,M)$ semisimple. To this end we
observe
that the image of $\a$ in $\gl (M)$ is normalized by the Weil operator $J$. So
the same is true for $\g (\a ,M)$. As $J$ is semisimple,
it is therefore enough to show that any
subspace $N\subset M\otimes\C$ that is invariant under both $\g (\a ,M)$ and
$J$
is nondegenerate with respect to $\phi$: then its $\phi$-perpendicular space
will be an invariant complement. Consider the primitive decomposition of
$N$  with respect to $e_a$: $N=\oplus _{k\ge 0}\C[e_k]P_{-k}(N)$, where
$P_{-k}(N):= \Ker (e_a^{k+1}|N_{-k})$. This decomposition is
$\phi$-perpendicular
and so we need to show that $\phi$ is nondegenerate on each summand
$\C[e_k]P_{-k}(N)$. For this we observe that $P_{-k}(N)$ is $J$-invariant.
Since
$H^k_a$ is definite on $P_{-k}(N)$, it follows from the definition of $H^k_a$
that $\phi$ is nondegenerate on $P_{-k}(N)+ e^kP_{-k}(N)$. The fact that
$e_k$ leaves $\phi$ infinitesimally invariant then implies that $\phi$ is
nondegenerate on $\C[e_k]P_{-k}(N)$.
\enddemo

\medskip\label
We briefly explain the relation between $\g (\a ,M)$
and the Mumford--Tate group. For this we have to assume that $M$, its grading,
$\a$, and the action of $\a$ on $M$ are all defined over $\Q$. Then $\g (\a
,M)$
is as a Lie subalgebra of $\gl (M)$ also defined over $\Q$. We further assume
that $\a$ acts by transformations of bidegree $(1,1)$. Then $\g (\a
,M)^{2k}$ acts by transformations of bidegree $(k,k)$, in other words, for all
$g\in\bS$ and $x\in\g (\a ,M)$ we have $gxg^{-1}=\Nm (g)^{-h}$.

Consider the image of $\bS$ in  $\GL (M)\times
\GL (1)$, where the second map is given by the norm. One defines the  {\it
Mumford--Tate group} of $M$, $\operatorname{MT}(M)$, as the smallest
$\Q$-subgroup of $\GL (M)\times \GL (1)$ containing this image. It is clear
that
this is actually a subgroup of $(\times _k\GL (M_k))\times\GL (1)$. The
projection of $\MT (M)$
onto the last factor is still called the {\it norm character} and denoted
likewise. The identity
$$ gxg^{-1}=\Nm (g)^{-h}$$
is now valid for all
$g\in\MT (M)$ and $x\in\g (\a ,M)$. This shows in particular that the adjoint
action of $\MT (M)$ on $\gl (M)$ leaves  $\g (\a ,M)$ invariant.

\medskip\label
This suggests to combine the Mumford--Tate group and  the group associated to
the
Lie algebra $\g (\a ,M)$ into a single group: if $G (\a ,M)$ denotes the closed
subgroup of $\GL (M)$ with Lie algebra $\g (\a ,M)$, then $\MT (M)G (\a ,M)$ is
a
reductive algebraic group defined over $\Q$. Its Lie algebra has a natural
Hodge
structure; it is obtained by composing the homomorphism $\bS \to \MT (M)G (\a
,M)$ with the adjoint action. In this set up the the r\^ole of the Deligne
torus
is played by the semidirect product $\bS\ltimes SL(2)$, where $s\in \bS$ acts
on  $SL(2)$ as conjugation by the diagonal matrix $\diag (z(s)^{-1},\bar
z(s))$.
So the $\Q$-homomorphism $w: GL(1)\to\bS\ltimes SL(2)$, which on the real
points
is given by   $t\in\R ^{\times}\mapsto (t,\diag (t,t^{-1}))$, maps onto a
central
subgroup. A a polarized Hodge structure of weight $d$ on  $(M ,\phi)$ can now
be
thought of as a certain representation of this group on $M$ with $w(t)$ acting
as
multiplication by $t^{-d}$.

The corresponding action of $\bS$ on $\slt$ is given by $s\mapsto (-z(s)-\bar
z(s))\ad _h$, so that for the resulting Hodge structure on $\slt$, the
bidegrees
of $e,h,f$ are $(1,1)$, $(0,0)$, $(-1,-1)$ respectively. In this spirit one can
also enhance the notion of a set of Shimura data, as defined by
\cite{Deligne}.

\medskip\label
It is high time to give the examples that motivated the preceding definitions.
Let $X$ be a compact K\"ahlerian manifold of dimension $n$. We take for $M$ its
shifted total complex cohomology $\Hm (X)[n]$ and we let $\phi$ be defined by
$$
\phi (\alpha ,\beta ):=(-1)^q\int _{X}\alpha\cup\beta
$$
if $\alpha$ is homogeneous of degree $n+2q$ or $n+2q+1$. (We shall
always suppose that $\Hm (X)$ is equipped with form and $\aut \Hm (X)$ will
stand
for the Lie algebra of endomorphisms of $\Hm (X)$ that preserve this form
infinitesimally.) The fundamental theorems of Hodge theory tell us that $M$
comes with a
Hodge structure of total weight $n$ and that $\phi$ together with cupping with
a
K\"ahler class defines a polarization of $M$. So by proposition \refer{1.6}
$\Hm (X)[n]$ is a Lefschetz module over $\Hm ^2(X)$. The
corresponding semisimple Lie subalgebra of $\aut \Hm (X)$ will be called the
{\it
total Lie algebra} of $M$ and be denoted $\g _{\tot}(X)$; it is defined over
$\Q$.
It is equivalent to say that the cohomology algebra of $X$ is a Lefschetz
algebra.
Clearly, $\g _{\tot}(X)$ is independent of the complex structure. For example,
if
$X$ is a product of an even number of circles, then  $\g _{\tot}(X)$ is
defined.

If we want to take the complex structure into account, then it is more natural
to regard $\Hm (X)[n]$ as module over $\Hm ^{1,1}(X)$. This is by \refer{1.6}
also
a Lefschetz module structure. We shall refer to the associated
Lie algebra as the {\it K\"ahler Lie algebra} of $X$ and denote it by $\g
_K(X)$.
For a complex projective manifold we can restrict further and take for $\a$ the
N\'eron--Severi group $\NS (X)$. We call the corresponding semisimple Lie
algebra
the {\it  N\'eron--Severi Lie algebra} of $X$ (denoted $\g _{NS}(X)$). It is
defined over $\Q$. Notice that the N\'eron--Severi Lie algebra and the
Mumford--Tate group behave in opposite ways under specialization: the former
gets
bigger, whereas the latter gets smaller.

If one of these Lie algebra's $\g _*(X)$ is defined over a subfield
$K\subset\C$, then we often write $\g _*(X;K)$ for the corresponding Lie
algebra
of $K$-points.

\smallskip\label
Here is another example. Let $V$ be a complex vector space and $W\subset
GL(V)$ a finite complex reflection group acting effectively (that is, $V^W=\{
0\}$). This group acts naturally in the symmetric algebra of $\Sym (V)$.
According to a theorem of Chevalley, the subalgebra of invariants, $\Sym (V)^W$
is a polynomial algebra on $\dim (V)$ homogeneous generators. Let $I$ be the
ideal generated by the invariants of positive degree. Then the quotient $\Sym
(V)/I$ is a graded complete intersection algebra. As a $W$-representation it is
isomorphic to the regular representation. In case $W$ is a Weyl group, then
after doubling the degrees, $\Sym (V)/I$ has the interpretation of the
cohomology algebra of a flag variety. From this we see that for a suitable
regrading, $\Sym (V)/I$ is a Lefschetz representation for the obvious action of
$V$. (We do not know whether this is true for an arbitrary reflection group.)
We shall determine the Lie algebra $\g (V,\Sym (V)/I)$  in \refer{5.8}.

\smallskip\label The N\'eron--Severi Lie algebra can also be defined when $X$
is
an irreducible projective variety: take for $M$ the total intersection
cohomology $\IH (X)$ with the same shift in the grading. There is in general
no such thing as a cup product on this graded vector space, but the cohomology
ring of $X$ acts on $M$ and according to \cite{Saito} polarizations have
the Hodge--Lefschetz property. This even extends to the case where we take
intersection cohomology with values in a local system defined on a Zariski
open-dense subset that underlies a polarized variation of Hodge structure.
There is an invariant form $\phi$ defined as in the previous example.  The
natural setting here is that of polarizable Hodge modules.

\medskip
We return to general properties of Lefschetz modules.

\proclaim{\label Proposition}
The Lie algebra $\g (\a ,M)$ is a Lefschetz module of $\a$. If
$a\in\a$ is such that $f_a$ is defined, then the Lie subalgebra $\g (\a
,M)_{\ge 0}$ (resp.\ $\g (\a ,M)_{\le 0}$) is generated by $\g (\a ,M)_0$ and
$e_a$ (resp.\ $f_a$). We have
$$
U\g (\a ,M) =U(\g (\a ,M)_{>0}).U(\g (\a ,M)_{0}).U(\g (\a ,M)_{<0}).
$$
\endproclaim
\demo{Proof} The first statement is clear and the second follows from this. The
third is clear also.
\enddemo

We define the {\it primitive subspace} of $M$ as the set of vectors killed by
$\g (\a ,M)_{<0}$. It is a graded $\g (\a ,M)_0$-subrepresentation of $M$ that
we denote by $\Prim (M)$. Since $\g (\a ,M)_{<0}$ is nilpotent, $\Prim (M)\not=
0$. The previous proposition yields:

\proclaim{\label Corollary}
In the situation  of the previous proposition, the primitive subspace $\Prim
(M)$ is the maximal $\g (\a ,M)_{0}$-invariant subspace contained in $\Ker
(f_a)$. Hence $M$ is irreducible as a $\g (\a ,M)$-representation if and
only if $\Prim (M)$ is irreducible as a $\g (\a ,M)_{0}$-representation (and if
$M\not=0$, then $\Prim (M)$ is the summand of $M$ of lowest degree).
\endproclaim

\smallskip The preceding  suggests to shift, in Tannakian spirit, the emphasis
from modules to Lie algebra's. For suppose that conversely, we are given a
semisimple Lie algebra $\g$, a simple element $h\in\g$ (in the sense of
appearing
as the middle element of an $\slt$-triple) and an abelian subalgebra $\a$ of
$\g$
such that \roster
\item"{(i)}" the adjoint representation of $\g$ makes $\g$ a
Lefschetz module over $\a$, i.e., there is a rational  map $f :\a \to \g_{-2}$
so
that for $e$ in the domain of $f$, we have an $\sli (2)$-triple $(e,h,f_e)$ and
\item"{(ii)}" $\g$ is as a Lie algebra generated by $\a$ and the image of $f$.
\endroster
If $M$ is a finite dimensional representation of $\g$, then $h$ determines a
grading
of $M$ and every $e\in\a$ in the domain of $f$ has the Lefschetz property in
$M$
with respect to this grading. So $M$ is then a Lefschetz module of $\a$. Since
$\g$
is generated by $\a$ and the image of $f$, it follows that $\g (\a ,M)$ is just
the
image of $\g$ in $\gl (M)$. This reduces the classification of Lefschetz
modules
to classifying triples $(\g, h,\a)$ as above. We shall call such a triple a
{\it
Lefschetz triple} and its first two items, $(\g ,h)$, a {\it
Lefschetz pair}. If we are given a Lefschetz pair $(\g ,h)$, then we say that
an
associated  Lefschetz triple $(\g, h,\a)$ is {\it saturated} if $\a$ is maximal
for
this property. The stabilizer $G_h$ of $h$ in the adjoint group $G$ permutes
these, but we do not know whether this action is transitive or even whether  it
has only finitely many orbits.

\smallskip
\label Let $(\g ,h)$ be a  Lefschetz pair. Choose a Cartan subalgebra $\h$ of
$\g$
that contains $h$. It is clear that then $\h\subset\g _{0}$.  Let $R\subset\h
^*$
denote the set of roots of $\h$ in $\g $ and let $R_k$ be the set of
$\alpha\in R$ such that $\alpha (h)=k$, or equivalently, $\g ^{\alpha}\subset\g
_k$
(remember that only even values of $k$ occur). The subset $R_0$ is a closed
root
subsystem of $R$; it is the set of roots of $\h$ in $\g _0$.  Choose a root
basis $B\subset R$ such that $\alpha (h)\ge 0$ for all $\alpha\in B$ (we then
say
that $\h$ and $B$ are {\it adapted}). Since $h$ is the semisimple element of a
$\slt$-triple, the elements of $B$ take on $h$  values in $\{ 0,1,2\}$
\cite{Bourbaki}, Ch.~VIII, \S 11, Prop.~5. Since these are also even, we get a
decomposition $B=B_0\sqcup B_2$. Clearly $B_0$ will be a root basis of $R_0$.
According to \cite{Bourbaki}, Ch.~VIII, \S 11, Prop.~8 all $\slt$-triples with
$h$ as semisimple element are conjugate under the stabilizer $G_h$ of $h$ in
the
adjoint group $G$ of $\g $. So if $(e,h,f)$ is such a triple and
$M$ is any representation of $\g$, then the isomorphism class of  $M$ as a
representation of this $\slt$-copy only depends on $(\g ,h)$. We call it the
{\it
$\slt$-type} of $M$. The following property narrows down the possible subsets
$B_2\subset B$. Let $V(k)$ denote the standard irreducible representation of
$\slt$ of dimension $k+1$ ($k=1,2,\dots$).

\proclaim{\label Proposition}
Let $(\g ,h)$ be a Lefschetz pair and let $M$ be an irreducible representation
of
$\g$ of depth $n$. Then the dimensions of the irreducible
$\slt$-representations
that occur in the $\slt$-type of $M$ make up an arithmetic progression with
increment $2$. In other words, there exists an integer $r$ with $0\le r\le
\lfloor {1\over 2}n\rfloor$ such that $\dim M_{-n}<\dim M_{-n+2}<\cdots <\dim
M_{-n+2r}=\dim M_{-n+2r+2}=\cdots =\dim M_{n-2r}<\dim M_{n-2r+2}<\cdots <\dim
M_n$. Moreover, $r>0$ unless (i) the image of $\g$ in $\gl (M)$ is reduced to
$\slt$ with $M\cong V(n)$ or (ii) the $\slt$-type of $M$ consists of a number
of
copies of $V(1)$.

\endproclaim
\demo{Proof} Denote this set of dimensions by $I$.
The irreducibility of $M$ implies that  the elements of $I$ all have  the
same parity. Suppose $I$ is not an arithmetic progression with increment $2$.
Then there exists an integer $k$ such that $M$ contains $V(k)$ (the standard
irreducible $\slt$-representation of dimension $k+1$) and $V(k+2l)$ for some
$l\ge 2$, but not $V(k+2)$. Let $(e,h,f)$ be an $\slt$-triple in $\g$
containing
$h$ and  decompose $M$ as $M=M'\oplus M''$ with $M'$ resp.\ $M''$ the sum of
the
irreducible subrepresentations of $\slt$ of $\dim\le k+1$ resp.\ $\ge k+5$. Any
linear transformation in $M$ of degree two that commutes with $e$ must preserve
this decomposition, because any $K[e]$-linear homomorphism
$$
V(n)\cong K[e]/(e^{n+1})[n]\to K[e]/(e^{m+1})[m]\cong V(m)
$$
of degree two is zero if $|n-m|>2$. If $(\g ,h,\a )$ is a Lefschetz triple with
$e\in\a$, then this applies in particular to any $e'\in\a$. If $e'$ has the
Lefschetz property, then $f_{e'}$ will also preserve this decomposition (since
$f_{e'}$ is unique) and hence $\g$ will. This contradicts our assumption that
$M$ is irreducible.

The second statement is proved in a similar way:  suppose $M\cong V(n)\otimes
P$
for some nonzero vector space $P$ with $n\ge 2$. If $e'$ is any element of
$\a$,
then the fact that $e'$ and $e$ commute, implies that $e'$ acts as $e\otimes
\sigma$ for some $\sigma\in\gl (P)$.  In this way, $\a$ maps onto subspace
$\bar\a$ of $\gl (P)$ that contains the identity of $P$. The elements of
$e\otimes\bar\a$ commute in $V(n)\otimes P$ and hence the elements of $\bar\a$
commute in $P$ (here we use that $n\ge 2$).

If $\dim P =1$, then we see that the image of $\a$ in $\gl (M)$ consists  of
multiples of $e$. This implies that the image of $\g$ in $\gl (M)$ is a copy of
$\slt$ and that $M\cong V(n)$. We now show that $\dim P\ge 2$ is impossible.
For
this we may assume that $K$ is algebraically closed. Then the commutative Lie
algebra $\bar\a$ leaves invariant a line $L\subset P$. So every $e'\in\a$
leaves
invariant $V(n)\otimes L$. If $e'$ has the Lefschetz property in $M$, then it
also has that property in $V(n)\otimes L$, and so the associated operator $f'$
leaves $V(n)\otimes L$ invariant. It follows that  $\g$ leaves $V(n)\otimes L$
invariant. This again contradicts the irreducibility of $M$. \enddemo

The question which subset  $B_2\subset B$ defines the semisimple element of an
$\slt$-triple is not difficult to answer for the classical Lie algebra's (see
\cite{Springer-St} and the discussion below) and for the exceptional Lie
algebra's this was tabulated by \cite{Dynkin 1952a}. (The weighted Dynkin
diagrams in these tables that matter here are only those that have all weights
$0$ or $2$.) The preceding proposition leads us to discard more possibilies,
but
we have not seriously studied the interesting question whether what is thus
left
actually occurs.

\smallskip
\label  Let us carry out this procedure in case $\g$ is a classical Lie
algebra.

So we assume that $\g$ is simple and classical and we let $V$ be a standard
representation of $\g$ (of dimension $l+1, 2l+1,2l,2l$ in  case $\g$ is of type
$A_l,B_l,C_l,D_l$ respectively). The element $h$ induces a grading on $V$. The
degrees that occur all have the same parity; we refer to this as the {\it
parity} of $V$. Let us recall that the $\slt$-invariant bilinear forms on
$V(k)$
are generated by a nonzero $(-)^k$-symmetric form. So an finite dimensional
$\slt$-representation of even parity always admits a nondegenerate invariant
symmetric form, whereas it admits a nondegenerate invariant
skew-symmetric form if and only if all multiplicities are even. In the case of
odd parity it is just the other way around.

The $\slt$-triple  $(e,h,f)$ in $\g$ determines a  primitive decomposition of
$V$. According to \refer{1.15} the set of positive integers $i\ge 0$ for which
$V(i)$ appears in the $\slt$-module $V$ is of the form  $\{ n,n-2,\dots
,n-2r\}$
(with $n-2r\ge 0$). So if we put $k:=\lfloor n/2\rfloor$,  then the dimensions
$d_t:=\dim V_{-n+2t}$ ($i=0,1,\dots ,k$) satisfy
$$
1\le d_0<d_1<\dots <d_r=d_{r+1}=\cdots =d_k.\tag{*}
$$
By the remark above, the $d_i$'s must all be even in the orthogonal cases with
odd
parity and in the symplectic cases with even parity.

We choose an $\h$-invariant basis of $V$ indexed as  in
\cite{Bourbaki}: $(e_1,\dots ,e_{l+1})$ in case $A_l$, $(e_1,\dots
,e_l,e_0,e_{-l},\dots ,e_{-1})$ in case $B_l$ and in the cases $C_l$ and $D_l$,
$(e_1,\dots ,e_l,e _{-l}\dots,e_{-1},)$. The same shall apply to our
labeling  of the simple roots $(\alpha _1,\alpha _2,\dots \alpha _l)$ (as
recalled below).

\smallskip
{\it Case $A_l$.} We let $\alpha _i(\diag (\lambda _1,\dots
,\lambda _{l+1}))=\lambda _i-\lambda _{i+1}$.
We find that the elements of $B_2$  are the simple roots with index
$d_0,d_0+d_1,\dots ,d_0+\cdots +d_k, d_0+\cdots +d_{k-1}+2d_k,\dots
,d_0+2(d_1+\cdots +d_{k-1}+d_k)$ (then $l=2(d_0+\cdots +d_k)$) or
$d_0,d_0+d_1,\dots ,d_0+\cdots +d_k, d_0+\cdots +2d_{k-1}+d_k,\dots
,d_0+2(d_1+\cdots +d_{k-1})+d_k$ (then $l=2(d_0+\cdots +d_{k-1})+d_k$).
So $B_2$ is symmetric with respect to the natural
involution of $B$. Notice that no two elements of $B_2$ will be adjacent in
the Dynkin diagram: if that would be the case, then $d_i=1$ for some $i>0$ and
\refer{*} shows that this is impossible.

\smallskip
{\it Case $B_l$.} Since $\dim V=2l+1$, its parity must be even. This
means that $n$ is even and $d_k$ is odd; we have
$l=d_0\cdots +d_{k-1}+{1\over 2}(d_k-1)$. We let
$$
\alpha _i(\diag (\lambda _1,\dots\lambda _l,0,-\lambda _{-l},\dots,-\lambda
_{-1}))= \cases
\lambda _i-\lambda _{i+1}&\text{if $i=1,\dots ,l-1$;}\\
\lambda _l&\text{if $i=l$.}
\endcases
$$
We find that the elements of $B_2$ are the simple roots with index
$d_0,d_0+d_1,\dots ,d_0+\cdots +d_{k-1}$. As in the previous case we see that
no two elements of $B_2$ will be adjacent in the Dynkin diagram.

\smallskip
{\it Case $C_l$.}
We let
$$
\alpha _i(\diag (\lambda _1,\dots ,\lambda _l,-\lambda _{-l},\dots,-\lambda
_{-1}))= \cases
\lambda _i-\lambda _{i+1}&\text{if $i=1,\dots ,l-1$;}\\
2\lambda _l&\text{if $i=l$.}
\endcases
$$

In the case of odd parity, $n$ is odd and $l=d_0+\cdots +d_k$. The elements
of $B_2$ are the simple roots with index $d_0,d_0+d_1,\dots ,d_0+\cdots
+d_k=l$.

In the case of  even parity $n$ and $d_0,\dots ,d_k$ are even. We have
$l=d_0+\cdots +d_{k-1}+{1\over 2}d_k$ and the elements of $B_2$ are
simple roots with index  $d_0,d_0+d_1,\dots ,d_0+\cdots +d_{k-1}$.

In both cases no two elements of $B_2$ will be adjacent in the Dynkin diagram.

\smallskip
{\it Case $D_l$, $\l\ge 4$.}
We let
$$
\alpha _i(\diag (\lambda _1,\dots ,\lambda _l,-\lambda _{-l},\dots,-\lambda
_{-1}))=  \cases
\lambda _i-\lambda _{i+1}&\text{if $i=1,\dots ,l-1$}\\
\lambda _{l-1}+\lambda _l&\text{if $i=l$;}\\
\endcases
$$
In the case of odd parity, $n$ is odd and all $d_i$'s are even. We have
$l= d_0+\cdots +d_{k-1}+d_k$ and the elements of
$B_2$ are simple roots with index  $d_0,d_0+d_1,\dots ,d_0+\cdots +d_k$. We
claim that no two such elements are adjacent in the Dynkin diagram. For this
can
only happen when $d_k=2$. In view of \refer{*}, this implies that $d_i=2$ for
all
$i$, in other words, $V\cong V(n)\oplus V(n)$. But since $n>1$, this is
excluded
by \refer{1.15}.

In the case of even parity, $n$ and $d_k$ are even. We have $l=
d_0+\cdots +d_{k-1}+{1\over 2}d_k$. Suppose first that $d_k\ge 4$.  Then the
elements of $B_2$ are the simple roots with index  $d_0,d_0+d_1,\dots
,d_0+\cdots
+d_{k-1}$. In that case $\alpha _l,\alpha _{l-1}\in B_0$ and $\alpha _1\in B_2$
if and only if $d_0=1$. Clearly no two elements of $B_2$ are adjacent.

If $d_k=2$, then by \refer{*}, $d_0=1$ and $d_2=\cdots =d_k=2$, in other words,
$V\cong V(2k)\oplus V(2k-2)$. Then $l=2k$ is even and the elements of $B_2$ are
in position $1,3,5,\dots,2k-3,2k-1,2k$.  No two elements of $B_2$ are
adjacent.

\smallskip
For future reference we record:

\proclaim{\label Corollary}
When $\g$ is simple and classical, no two elements of $B_2$ are connected.
\endproclaim

\head
\section Jordan--Lefschetz modules
\endhead

\noindent
In this section we discuss a particularly nice class of Lefschetz triples. We
introduce them via the following proposition.

\proclaim{\label Proposition}
Let $M$ be a Lefschetz $\a$-module. Then the following two properties are
equivalent
\roster
\item"{(i)}"  The graded Lie algebra  $\g (\a ,M)$ has degrees $-2$, $0$ and
$2$
only.
\item"{(ii)}" The operators $f_a$ with $a\in\a$ in the domain of $f$,
mutually commute. \endroster
\endproclaim

\demo{Proof}  We only prove the nontrivial implication $(ii)\Rightarrow (i)$.
Write $\g$ for $\g (\a ,M)$.

If $a\in\a$ is such that $f_a$ is defined, then regard $\g$ as an $\slt$-module
via the $\slt$-triple $(e_a ,h,f_a )$ and let $V(a)\subset\g$ be sum of the
irreducible summands of dimension $1$ and $3$. Notice that $V(a)$ contains the
image of $\a$.  So the intersection $V$ of the $V(a)$'s also contains the image
of $\a$. We have $V=V_{-2}\oplus V_0\oplus V_2$ where $V_{2k}$ can be
characterized as the subspace of $\g _{2k}$ that is
annihilated by all the operators $\ad ^{-k+2}(e_a)$, or alternatively, by all
the operators $\ad ^{k+2}(f_a)$. Since the operators $e_a$ resp.\ $f_a$
mutually commute it follows that $V$ is invariant
under these operators. Hence $V=\g$ and so $\g$ has degrees $-2$, $0$
and $2$ only. \enddemo

As we are interested in the graded Lie algebra's that so arise, we make the
following definition. Say that a Lefschetz pair $(\g ,h)$ is a {\it
Jordan--Lefschetz pair} if $(\g ,h,\g _2)$ is a Lefschetz triple. It is
clear that such a Lefschetz triple is saturated.

\proclaim{\label Proposition}
If $(\g ,h)$ is a Jordan--Lefschetz pair, then $\g =\g
_{-2}\oplus\g _0\oplus\g _2$ and $U\g =U\g _2.U\g_0.U\g _{-2}$.
\endproclaim
\demo{Proof}
If $e\in\g _2$ is a Lefschetz
element, then $[e,\g _2]=0$. The first assertion now follows from the primitive
decomposition under $\ad _e$. The second is a consequence of this.
\enddemo

\proclaim{\label Corollary}
Let $(\g ,h)$ be a Jordan--Lefschetz pair and let $M$ be a finite dimensional
representation of $\g$ and regard $M$ as a Lefschetz module of $\g _2$. Then
$M$
is generated as a $U\g _2$-module by  $\Prim (M)$.
\endproclaim

Our terminology is easily explained: according to \cite{Springer}, \S 2.21,
a Jordan--Lefschetz pair defines a ``Jordan algebra without unit
element''. These have been classified. Let
us give a quick proof of this classification. It is based on two lemma's.
In what follows, $(\g ,h)$ is a Jordan--Lefschetz pair with $\g$ simple and
and $\h\subset\g$ and $B$ are adapted.

\proclaim{\label Lemma} The subset $B_2$ is a singleton.
\endproclaim
\demo{Proof}
If not, then there is chain $(\beta ,\alpha _1,\dots ,\alpha _N,\beta ')$ with
$\beta ,\beta '\in B_2$ be distinct and $\alpha _i \in
B_0$. But then is $\beta +\alpha _1+\cdots \alpha _N +\beta '$ a root and in
$R_4$, which is supposed to be empty.
\enddemo

It follows that $h$ is twice a fundamental coweight.

The following lemma concerns a general property of root systems.

\proclaim{\label Lemma}
Let $\beta\in B$ and let $R(\beta )$ be the set of positive roots that have
$\beta$-coefficient one. Then $R(\beta )+R(\beta )$ does not contain a
root if and only if $R(\beta )$ contains the highest root.
\endproclaim \demo{Proof}
It is clear that if $R(\beta )$ contains the highest root, then
$R(\beta )+R(\beta )$ cannot contain a root. If on the other hand
the highest root has $\beta$-coefficient $\ge 2$, then it follows from
\cite{Bourbaki}, Ch.~VI,\S 1, Prop.~19 that there exists a sequence of roots
$\alpha _1,\dots ,\alpha _r$ such that $\alpha _1\in B$, $\alpha _{i+1}-\alpha
_i\in B$ for $i=1,\dots ,r-1$ and $\alpha _r$ is the highest root. If $\alpha
_i$ is the last root in this sequence for which the coefficient of $\beta$ is
one, then $\alpha _i$ and $\beta$ are elements of $R(\beta)$ whose sum is a
root.
\enddemo

\proclaim{\label Corollary} Let $\beta$ be the unique element of $B_2$.
Then the pair $(B,B-\{\beta\})$ is of the following type:
\roster
\item"{}" $(A_{2m-1},A_{m-1}+A_{m-1})$ ($m\ge 1$),
\item"{}" $(B_m,B_{m-1})$ ($m\ge 2$),
\item"{}" $(C_m,A_{m-1})$ ($m\ge 2$),
\item"{}" $(D_m,D_{m-1})$ ($m\ge 5$),
\item"{}" $(D_{2m},A_{2m-1})$ ($m\ge 2$) or
\item"{}" $(E_7,E_6)$.
\endroster
Conversely, every item of this list determines an isomorphism class of
Jordan--Lefschetz pairs.
\endproclaim
\demo{Proof} The pairs $(B,\beta )$ with the property of the previous lemma are
those in this list plus the following: in case $A_l$ we can take $l$
and $\beta$ arbitrary, in case $D_{2l+1}$ any end of the Dynkin diagram,
and $B$ of type $E_6$ with $\beta$ the
end of the branch of length $3$. These addional possibilities disappear if we
want $h$ to be a simple element (that is, $h=[e,f]$ for certain $e\in\g _2$ and
$f\in \g _{-2}$): for the classical cases $A_l$ and $D_{2l+1}$ this follows
from the
discussion following \refer{1.16}) and for $E_6$ this follows from table 18 of
\cite{Dynkin}. On the other hand, case $3''$ of table 19 of this reference
shows
that in case $E_7$ this element $h$ is simple. The discussion below will show
that
all the other listed cases occur as well.  \enddemo

\label We continue with the Jordan--Lefschetz pair $(\g ,h)$ and let $\h$, $B$
and $\beta\in B$ as above. Let $\varpi\in\h ^*$ be the fundamental weight
that is zero on $B^{\vee}-\{\beta ^{\vee}\}$ and $1$ on $\beta ^{\vee}$. We
say that an irreducible representation $M$ of $\g$ is a {\it Jordan--Lefschetz
module of level $k$} (where $k$ is a positive integer) if its lowest weight is
$-k\varpi$; for $k=1$, we shall also call it a {\it fundamental}
Jordan--Lefschetz
module.
Notice that $M$ then has depth $k\varpi (h)$. A more intrinsic
description of these modules is the following: given a Jordan--Lefschetz pair
$(\g ,h)$, then an irreducible representation $M$ of $\g$ is a
Jordan--Lefschetz
module if and only if it has a nonzero vector stabilized by  $\g _{-2}+\g _0$.

\proclaim{\label Lemma}
Let $M$ be an irreducible representation $M$ of the Jordan--Lefschetz pair
$(\g ,h)$ and let $n$ be its depth as a Lefschetz module.
Then $M$ is a Jordan--Lefschetz module if and only if $\dim M_{-n}=1$.
If these equivalent conditions are fulfilled, then
the natural map $\g ^2\otimes M_{-n}\to M_{-n+2}$ is an isomorphism.
\endproclaim
\demo{Proof} Let $\lambda\in\h ^*$ be the lowest weight of $M$.
Then the lowest weight space $M^{-\lambda}$ is
contained in $M_{-n}$. Hence
$\dim M_{-n}=1$ is equivalent to $M^{-\lambda}=M_{-n}$. The latter is
equivalent to: $\g _{-2}+\g _0$ stabilizes  $M^{-\lambda}$, which in turn is
equivalent to $\lambda$ vanishing on $B_0$.

Suppose the two conditions satisfied. The fact that $\g _{-2}+\g _0$ is the
$\g$-stabilizer of $M_{-n}$ implies that $\g ^2\otimes M_{-n}\to M_{-n+2}$
is injective. Since $M$ is as a $U\g _2$-module generated by $M_{-n}$, it is
also surjective.
\enddemo

\label Let us describe the fundamental Jordan--Lefschetz representations
in the classical cases (i.e., those that are not of type $(E_7,E_6)$). We only
do this over the complex numbers.

\smallskip
{\it Case $(A_{2m-1},A_{m-1}+A_{m-1})$}: $(\sli (2m),\sli (m)\times\sli (m))$.
Let $V$ be a vector space of dimension $2m$, $V=V_{-1}\oplus V_{1}$ a direct
sum
decomposition into subspaces of dimension $m$. We take $\g :=\sli (V)$  and let
$h\in \sli (V)$ be $\pm 1$ on $V_{\pm 1}$.
Then $\g _0$ maps isomorphically onto $\sli (V_{-1})\times\sli (V _1)\times\C
h$
and  $\g _2\cong \Hom (V_{-1},V_1)$ resp.\  $\g _{-2}\cong \Hom (V_1,V_{-1})$.
The corresponding fundamental representation is $M:=\wedge ^m V$ with lowest
degree piece $M_{-m}=\wedge ^mV_{-1}$.

\smallskip
{\it Case $(B_m,B_{m-1})$ or $(D_m,D_{m-1})$}:  $(\so (n),\so (n-2))$ with
$m=2n+1$ resp.\ $m=2n$ ($n\ge 2$).
Let $V$ be a vector
space of dimension $n$ equipped with a nondegenerate symmetric bilinear form
and let $V_{\pm 2}$ be isotropic lines in $V$ such that $V_{-2}\oplus V_2$ is
nondegenerate. Let $V_0$ be the orthogonal complement of $V_{-2}\oplus V_2$ in
$V$. We take $\g =\so (V)$ and let $h\in\so (V)$ be the element with the eigen
space decomposition  $V_{-2}\oplus V_0\oplus V_2$. Then $\g _0=\so (V_0)\times
\gl (V_2)$ and $\g _{\pm 2}$ projects isomorphically to $\Hom (V_0,V_{\pm
2})$. We take $M=V$.

\smallskip
{\it Case $(C_m,A_{m-1})$}: $(\sy (2m),\sli (m))$ ($m\ge 2$). Let $V$ be a
vector
space of dimension
$2m$ equipped with a nondegenerate symplectic form and let $V=V_{-1}\oplus V_1$
be a decomposition of $V$ into totally isotropic subspaces of dimension $m$. We
take $\g =\sy (V)$ and let $h\in\sy (V)$ be the element with the eigen space
decomposition  $V_{-1}\oplus V_1$. Then $\g _0$ maps isomorphically to $\gl
(V_1)$ and $\g _{\pm 2}$ is naturally isomorphic to the space of symmetric
elements in $(V_{\pm 1})^{\otimes 2}$. We take for $M$ the primitive quotient
of
$\wedge ^mV$ (i.e., the quotient by $\wedge ^{m-2} V\wedge\omega$, where
$\omega\in \wedge ^2V$ is the dual of the symplectic form).

\smallskip
{\it Case $(D_{2m},A_{2m-1})$}: $(\so (4m),\sli (2m))$ ($m\ge 2$).
Let $V$ be a vector space of
dimension $4m$ equipped with a nondegenerate symmetric bilinear form and let
$V=V_{-1}\oplus V_1$ be a decomposition of $V$ into totally isotropic subspaces
of dimension $2m$. We take $\g =\so (V)$ and let $h\in\so (V)$ be the element
with the eigen space decomposition  $V_{-1}\oplus V_1$. Then $\g _0$ maps
isomorphically to $\gl (V_1)$ and $\g _{\pm 2}$ maps onto the skew elements in
$(V_{\pm 1})^{\otimes 2}$. Consider the spinor representation $\wedge
^{\bullet} V_1$. (We recall that this factors through the representation of the
Clifford algebra of $V$ on $\wedge ^{\bullet} V_1$ for which $v\in V_1$ acts as
wedging with $v$ and $v\in V_{-1}$ acts as the interior product under the
obvious isomorphism $V_{-1}\cong V_1^*$.) It splits into a direct sum of
subrepresentations $\wedge ^{\ev} V _1$ and $\wedge ^{\odd} V _1$. They are
irreducible and nonisomorphic and correspond to the case when $\beta$ is an end
of the Dynkin diagram connected with a branch point. For even $m$ both
are orthogonal and for odd $m$ both are symplectic. We take
$M=\wedge ^{\ev} V _1[m]$.

\smallskip
\label The Jordan--Lefschetz modules give rise to Frobenius
algebra's with remarkable properties. Let $(\g ,h)$  be a Jordan pair and let
$M$
be a Jordan--Lefschetz module of $(\g ,h)$ of depth $n$. Then $M$ is a monic
module
over the commutative algebra $U\g _2$. So if $I\subset U\g _2$ denote the
annihilator of $M$, then $I$ defines the origin in $\Spec (U\g _2)=\g _2^*$
with
local algebra $A:=U\g _2/I$. The latter is an evenly graded Lefschetz algebra
that
has $M$ as a free graded module of rank one. The next proposition shows that it
has a lot of automorphisms. Let $\g _0'$ be the Lie subalgebra of $\g _0$ that
kills $1\in A_0$ (or equivalently, kills $M_{-n}$); this Lie algebra is
complementary to the span of $h$ in $\g _0$.

\proclaim{\label Proposition}
The Lie algebra $\g '_0$ acts on $A$ as derivations and so the associated Lie
subgroup of $\GL (A)$ is a group of algebra automorphisms of $A$. Moreover, the
Lie subgroup $G_0\subset\GL (A)$ associated to $\g _0$ has a dense orbit in
$A_2$ consisting of Lefschetz elements.
\endproclaim
\demo{Proof}
If $u\in\g _0$ and $e\in\g _2$, then for all $x\in A$, $u(ex)=
[u,e]x+ e(ux)$. Since $A$ is generated by $\g _2$, it follows with  induction
that $u$ acts as a derivation if (and only if) $u$ kills $1$. The  last
assertion
follows from a well-known result \cite{Bourbaki}, Ch.~VIII, \S 11, Prop.\ 6.
\enddemo

\label
A decomposition of $\g$ into simple components corresponds to a decomposition
of
$A$ as a tensor product of algebra's, so little is lost in assuming that $\g$
is
simple. The form  $\int$ defined in \refer{1.4} makes of $A$ a Frobenius
algebra
with soccle $A_{2n}$; at the same time $\int$ defines a generator of $A_{2n}$
that
serves as the identity element for another algebra structure defined by the
action of $U\g _{-2}$. Here is a description of the
algebra's associated to the fundamental Jordan--Lefschetz modules in all
cases.

\smallskip
{\it Case $(A_{2m-1},A_{m-1}+A_{m-1})$}:
Let $W$, $W'$ be vector spaces of dimension $m$ and let $A:= \oplus _{k=0}^m
\wedge ^k W\otimes\wedge ^ kW'$. This is just the subalgebra of the graded
algebra $\wedge ^{\bullet}W\otimes\wedge ^{\bullet}W'$ generated by $W\otimes
W'$. It is the fundamental Frobenius algebra associated to $(\sli (W\oplus
W'),h=(1_W,-1_{W'}))$. (To see the relation with the description given in
\refer{2.9}, take $V_{-1}=(W')^*$ and $V_1=W$ and observe that a choice of a
generator of $\wedge ^mV_{-1}$ identifies $\wedge ^{m-k}V_{-1}$ with
$\wedge ^k W'$ and hence $\wedge ^m(V_1\oplus V_{-1})$ with $A$.)

Here is a presentation of this algebra.  Choose bases $(w_1,\dots
,w_m)$ of $W$ and $(w'_1,\dots ,w'_m)$ of $W'$. Then $x_{ij}:=w_i\otimes w'_j$,
($i,j=1,\dots ,m$) generate $A$ as an algebra and a set of defining relations
is
 $x_{ij}x_{kl}+x_{il}x_{kj}=0$.

\smallskip
{\it Case $(D_{2m},A_{2m-1})$}: Let $W$ be a
vector space of dimension $2m$, then let $A$ be the subalgebra of $\wedge
^{\bullet}W$ generated by $\wedge ^2W$. This is the
fundamental Frobenius algebra associated to $\so (W\oplus W^*),
h=(1_W,-1_{W^*}))$.

A presentation of $A$ is a  follows $(w_1,\dots ,w_m)$ is a basis of $W$, then
generators for $A$ are $\omega _{ij}:=w_i\wedge w_j$ ($1\le i,j\le m$), subject
to the linear relations  $\omega _{ij}+\omega _{ji}=0$ and the quadratic
relations $\omega _{ij}\omega _{kl}+\omega _{il}\omega _{kj}=0$. So this is a
quotient of the algebra of the previous case.

\smallskip
{\it Case $(C_m,A_{m-1})$}: Let $W$ be a vector space of dimension $m$. Then
a model for the fundamental Frobenius algebra associated to $\sy (W\oplus
W^*), h=(1_W,-1_W))$ is the subalgebra $A$ of $\wedge
^{\bullet}W\otimes\wedge ^{\bullet}W$ generated by the symmetric elements
$w\otimes w$, $w\in W$. If $(w_1,\dots ,w_m)$ is a
basis of $W$, then generators for $A$ are $u _{ij}:=w_i\otimes w_j +w_j\otimes
w_i$ ($1\le i,j\le m$). A defining set of relations is $u_{ij}=u_{ji}$ and
$u_{ij}u_{jk}+u_{jj}u_{ik}=0$.

\smallskip
{\it Cases $(B_n,B_{n-1})$ and $(D_n,D_{n-1})$:} Let $(W,q:W\to K)$ be a
vector space of dimension $m$ endowed with a nondegenerate quadratic form.
Consider the graded vector space $A:=K\oplus W\oplus K\mu$, where $K$ has
degree zero, $W$ has degree $2$ and $\mu$ is a generator of a
one-dimensional $K$-vector space that has degree $4$. Equip $A$ with
the quadratic form $\tilde q(x+w+x'\mu)=q(w)-xx'$.
Make it also a graded algebra by letting $1\in K$ be the identity element and
letting for $w,w'\in W$, $w.w'=b(w,w')\mu$, where $b$ is the bilinear form
associated to $q$. Then $A$ is the fundamental Frobenius algebra associated to
$(\so (A,\tilde q), h=(-2,0_W,2))$.
If $(w_1,\dots ,w_m)$ is a $q$-orthogonal basis of $W$ then a presentation of
$A$
has generators $w_1,\dots ,w_m$ and relations $w_iw_j=0$ ($i\not=j$) and
$q(w_i)^{-1}w_iw_i=q(w_j)^{-1}w_jw_j$.

\smallskip
{\it Case  $(E_7,E_6)$:}
Let $W$ be a finite dimensional $K$-vector space
and let be given a cubic form $c:\Sym ^3(W)\to K$ that does not factor through
proper linear quotient of $W$. The latter condition implies that the associated
homomorphism $\tilde c:\Sym ^2W\to W^*$ is surjective. Consider the graded
vector
space  $A:=K\oplus W\oplus W^*\oplus K\mu$ with respective summands in degree
$0$, $2$, $4$, $6$ (here $\mu$ is a generator of a one dimensional vector
space)
and use $\tilde c$ and the obvious pairing $W\times W^*\to
K\mu$ to give $A$ the structure of a commutative graded $K$-algebra. We
require that multiplication makes $A[3]$ a Jordan--Lefschetz module over
$A_2$.
This forces $\dim W=27$ and after possibly passing to an algebraic closure of
$K$, $c$ will be unique up to a linear transformation. The group $G_c$ of
$g\in GL(W)$ that leave $c$ invariant is of type $E_6$ and the
associated Lie algebra $\g$ is of type $E_7$. (This Lie algebra can be
characterized as the Lie algebra of linear transformations of $A$ that leaves
invariant a certain quartic form on $A$.) If we regard $A$ as a quotient of the
symmetric algebra of $W$, then we checked by means of the program \cite{LiE}
that the relations are again quadratic: the ideal $I\subset\Sym (W)$ that
defines $A$ is generated by $\ker (\tilde c)$ (the latter is an irreducible
representation of $G_c$ whose highest weight is twice that of $W$).

\medskip
The Jordan--Lefschetz algebra's of higher level can be expressed
in terms of a fundamental one:

\proclaim{\label Proposition}
Let $A$ be a fundamental Jordan--Lefschetz algebra
for a Jordan--Lefschetz pair $(\g ,h)$. Let $k$ be a positive integer and give
$\Sym
^k(A)$ the structure of an algebra by identifying it with the algebra of
symmetric
invariants in $A^{\otimes k}$. Then the subalgebra $A(k)$ of $\Sym ^k(A)$
generated
by $A_2$ is  a Jordan--Lefschetz algebra of level $k$ for $(\g ,h)$.
\endproclaim
\demo{Proof}
We regard $A$ as an fundamental representation of $\g$ with lowest weight space
spanned by its unit element $1\in A$. Then the irreducible representation of
$\g$
with lowest weight $k$ times the one of $A$ is contained in $A^{\otimes k}$ as
the $U\g$-submodule generated by the unit element $1\otimes\cdots\otimes 1$.
But this is also the $U\g_2$-submodule generated by this element. In other
words,
this is the subalgebra of $A^{\otimes k}$ generated by the elements
$a\otimes 1\otimes\cdots\otimes 1 + 1\otimes a\otimes\cdots\otimes 1 +\cdots +
1\otimes 1\otimes\cdots\otimes a$, with $a\in A_2$.
\enddemo

In section 4 we will be concerned with the algebra's of higher level in the
cases
$(B_n,B_{n-1})$ and $(D_n,D_{n-1})$, and that is why we want to describe them
here
in more explicit terms. Let $(W,q)$ and $(A,\tilde q)$ be as under the relevant
case above.

\proclaim{\label Proposition} Fix an integer $k\ge 1$ and let $I_k$ be the
ideal in $\Sym (W)$ generated by the $k+1$-st powers of $q$-isotropic vectors
(i.e., the $w^{k+1}$ for which $q(w)=0$).  Then $A(k)=\Sym (W)/I_k$. Moreover,
the
subalgebra of $\Sym (W)/I_k$ of $\so (W)$-invariants is $K(u)/(u^{k+1})$,
where $u\in \Sym ^2(W)$ represents the inverse form of
$q$ on the dual of $W$. The soccle of $A(k)$ is spanned by
the image of $u^k$ and if we identify $W$ with its image in $A(k)$, then
$x^2u^{k-1}=q(x)u^k$.

\endproclaim
\demo{Proof} Consider the graded algebra obtained by dividing $\Sym (W)$ out
by the
ideal generated by $u$,  $\Sym (W)/(u)$, and let $W(d)$ denote the image of
$\Sym
^d(W)$. Then $W(d)$ is an irreducible representation of $\so (W)$ that can be
identified with the subrepresentation of $\Sym ^d(W)$ which is linearly spanned
the
pure $d$th powers of $q$-isotropic vectors in $W$. In particular, $I_k$ is the
ideal
in $\Sym (W)$ generated by $W(k+1)$. It is easy to see that
$$
\Sym ^d(W)=\oplus _{i=0}^{\lfloor d/2\rfloor} u^iW(d-2i),
$$
as graded $\so (W)$-representations. A Clebsch--Gordan
rule asserts that for $p\ge q\ge 0$, the image of $W(p)\otimes W(q)\to \Sym
^{p+q}W$ under the multiplication map is $\oplus _{i=0}^q u^iW(p+q-2i)$. So the
cokernel can be identified with $u^{q+1}\Sym ^{p-q-2}(W)$. This remains so if
we replace $W(q)$ by $\Sym ^q(W)$. It follows that
$$\align
\Sym (W)/I_k =&\Sym ^0(W)\oplus\cdots \oplus\Sym ^{k-1}(W)\oplus\Sym ^k(W)\\
&\oplus u\Sym ^{k-1}(W)\oplus u^2\Sym ^{k-1}(W)\oplus \cdots\oplus u^k\Sym
^0(W).
\endalign
$$
If $\tilde u\in\Sym ^2(A)$ denotes the symmetric tensor dual to $\tilde q$,
then
$A(k)$ is the image of $\Sym ^k(A)$ in $\Sym (A)/(\tilde u)$. If we write
$A=Kt\oplus W\oplus K\mu$, with $\deg t=0$, then $\Sym (A)=K[t,\mu ]\Sym W$
with $\tilde u$
corresponding to $t\mu +u$. Hence we have the following identity of graded $\so
(W)$-representations:
$$
A(k)=\oplus _{i=0}^k t^{k-i}\Sym ^i(W)\,\oplus\,\oplus _{i=0}^k\mu
^{k-i}\Sym ^i(W).
$$
We know that there is a surjective graded algebra homomorphism
$\Sym (W)\to A(k)$. This homomorphism is $\so (W)$-equivariant and hence
contains $W(k+1)$ in its kernel. We therefore have a graded, $\so
(W)$-equivariant, algebra epimorphism $\Sym (W)/I_k\to A(k)$. The
$\so (W)$-decompositions that we found for both source and target show that
this must be an isomorphism.

The last two assertions are clear.
\enddemo

\head
\section Geometric examples of Jordan type I: complex tori
\endhead

\noindent
In this section we describe the total Lie algebra of a complex torus, the
K\"ahler Lie algebra of a complex torus, and the N\'eron--Severi Lie algebra of
an
abelian variety.

\medskip\label We begin with determining the total Lie algebra of a complex
torus.
As observed before this is independent of its complex structure.  We first
recall
the construction of the spinor representations.

\smallskip
Let $V$ be a real finite dimensional vector space.
For $\alpha\in V^*$, left exterior product with
$\alpha$ defines a map $e_{\alpha}:\wedge
^{\bullet}V^*\to \wedge ^{\bullet +1}V^*$ of degree one.  Dually, contraction
with $a\in V$ defines a map $i _a: \wedge^{\bullet}V^*\to \wedge
^{\bullet -1}V^*$ of degree minus one. Both are derivations of odd degree:
$i_a(\omega\wedge\omega ')= i_a(\omega )\wedge\omega '+(-)^{\deg\omega}\omega
\wedge i_a(\omega ')$ and similarly for $e_{\alpha}$.  The anti-commutator
$i_a e_{\alpha}+e_{\alpha}i_a$ is simply multiplication by $\alpha (a)$.
Iterated use of this identity shows that for $a,b\in V$ and $\alpha,\beta\in
V^*$, we have
$$
\align
[i_a i_b,e_{\alpha} e_{\beta}] =&
-\alpha (b) e_{\beta}i_a +\alpha (a) e_{\beta}i_b
+\beta (b) e_{\alpha}i_a - \beta (a) e_{\alpha} i_b\\
&+(-\alpha (a)\beta (b)+\beta (a)\alpha (b))\bold{1}_{\wedge ^{\bullet}V^*}
.\tag{**}
\endalign
$$
We rewrite this as follows: if we denote by $\sigma(a\wedge b,\alpha\wedge\beta
)\in\gl (V^*)$ the transformation
$$
\xi\mapsto
-\alpha (b)\xi (a)\beta +\alpha (a)\xi (b)\beta
+\beta (b)\xi (a)\alpha - \beta (a)\xi (b)\alpha ,
$$
and write $\tilde\sigma(a\wedge b,\alpha\wedge\beta )$ for its extension as a
degree
zero derivation in $\wedge ^{\bullet}V^*$, then the righthand side of
eq.~\refer{**}  is simply the sum of $\tilde\sigma(a\wedge b,\alpha\wedge\beta
)$
and the scalar operator that is multiplication by $-{1\over 2}\Tr
(\sigma(a\wedge
b,\alpha\wedge\beta ))$.

With the help of this formula we determine the Lie subalgebra $\g$ of
$\gl (\wedge ^{\bullet} V^*)$ generated by these operators. In fact,
by means of a Clifford construction we shall identify it with the Lie
subalgebra $\so (V\oplus V^*)$ of  $\gl (V\oplus V^*)$ of infinitesimal
automorphisms  of the quadratic form $q(x,\xi )=\xi (x)$.

The semi-simple element $u:=(-\bold{1}_V,+\bold{1}_{V^*})\in\so (V\oplus V^*)$
defines a grading of the latter with degrees $2$, $0$ and $-2$. We define a Lie
algebra isomorphism
$$
\psi _0: \gl (V ^*)\to \so (V\oplus V^*)_0;\quad
\psi _0(\sigma )(x,\xi )=(-\sigma ^*(x),\sigma (\xi )).
$$
We also have isomorphisms of abelian Lie algebra's
$$
\align
\psi _2: \wedge ^2V^*\to \so (V\oplus V^*)_2;&\quad
\psi _2(\alpha\wedge\beta )(x,\xi)=(0,\alpha (x)\beta -\beta (x)\alpha ),\\
\psi _{-2}: \wedge ^2V\to \so (V\oplus V^*)_{-2};&\quad
\psi _{-2}(a\wedge b)(x,\xi)=(\xi (a)b-\xi (b)a,0).
\endalign
$$

\proclaim{\label Proposition}
The maps $\psi _2(\alpha\wedge\beta )\mapsto e_{\alpha}e_{\beta}$,
$\psi _{-2}(a\wedge b)\mapsto i_ai_b$ extend to a graded Lie algebra
isomorphism of $\so (V\oplus V^*)$ onto  $\g$ that maps $\psi _0(\sigma )$ to
$\tilde\sigma -{1\over 2}\Tr (\sigma )\bold{1}_{\wedge ^{\bullet}V^*}$, where
$\tilde\sigma$ denotes the  extension of $\sigma$ as a degree zero derivation
of $\wedge ^{\bullet}V^*$.
\endproclaim
\demo{Proof} In view of formula \refer{**} it suffices to show that
$[\psi _{-2}(a\wedge b),\psi _2(\alpha\wedge\beta )]=
\psi _0(\sigma(a\wedge b,\alpha\wedge\beta ))$.
This is straightforward.
\enddemo

The pair $(\so (V\oplus V^*),u)$ is a Jordan--Lefschetz pair (it is a real
form  of
case $(D_{2n},A_{2n-1})$).

Now let $X$ be a real torus of even dimension $2n$. We identify the
universal cover of $X$  with $\Hm _1(X;\R )$.
We will write $V$ for this real
vector space (of dimension $2n$) so that $\Hm (X;\R )=\wedge
^{\bullet}V^*$. The rational homology defines a rational structure
$V_{\Q}\subset V$. Let $\kappa\in \wedge ^2V^*$ be nondegenerate. If $\alpha
_{\pm 1},\dots ,\alpha _{\pm n}$ is a basis of $V^*$ such that $\kappa =\sum
_{k=1}^n \alpha _k\wedge\alpha _{-k}$ then $e_{\kappa}=\sum _{k=1}^n e_{\alpha
_k}e_{\alpha _{-k}}$. If $a_{\pm 1},\dots a_{\pm _n}$ is the dual basis, then
we
see from eq.~\refer{**} that
$$
[e_{\kappa},\sum _{k=1}^n i_{a_{-k}}i_{a_k}]=-n+\sum _{k=1}^n
(e_{\alpha _k}i_{a_k}+e_{\alpha _{-k}}i_{a_{-k}}).
$$
Since this element acts on $\wedge
^lV^*$ as multiplication by $-n+l$, it follows that $f_{\kappa}$ is defined
and equal to $\sum _{k=1}^n i_{a_{-k}}i_{a_k}$. The nondegenerate $2$-forms
make up a nonempty open subset of $\wedge ^2V^*$ and therefore span that
space.
The corresponding $2$-vectors form
an open subset of $\wedge ^2V$ and so  $\g _{\tot}(X)$ is generated by $\g
_2\oplus\g _{-2}$ as a Lie algebra. Combining this with the above computation
gives:

\proclaim{\label Proposition}
There is a natural identification $(\g _{\tot}(X;\R ),h )\cong (\so (V^*\oplus
V),u)$; this is a real form of the case $(D_{2n},A_{2n-1})$.  Furthermore, $\Hm
^{\ev}(X)[n]$ is a semispinorial representation of
$g _{\tot}(X;\R )$ and a fundamental Jordan--Lefschetz module of $\Hm ^2(X,\R
)$.
\endproclaim

\medskip
We now assume that $X$ comes with a complex structure. We shall
determine its K\"ahler Lie algebra. Let $V^*$ resp.\ $\bar V^*$ denote the
$\R$-dual of $V$ equipped with the complex structure
$J^*$ resp.\ $-J^*$. The quadratic form $q$ defined above is invariant under
the
complex structure $(J,-J^*)$ and so extends to a Hermitian form on $V\oplus\bar
V^*$. Let $\su (V\oplus V^*)$ be the Lie algebra of the corresponding special
unitary group. Since $u\in\su (V\oplus \bar V^*)$, it inherits a grading with
degrees $-2$, $0$ and $2$. In fact, $(\su (V\oplus \bar V^*),u)$ is a
Jordan--Lefschetz pair. It is a real form of case $(A_{2n-1},
A_{n-1}+A_{n-1})$.

The $J^*$-invariant elements  of $\wedge _{\R}^2V^*$ are
precisely the real $(1,1)$-forms. So they make up the span of the K\"ahler
classes of $X$. Now a straightforward verification shows that $\psi$ maps
the $J$-invariants of $\wedge _{\R}^2V$ resp.~$J^*$-invariants of $\wedge _{\R}
^2V^*$ onto $\su (V\oplus \bar V^*)_{-2}$ resp.~$\su (V\oplus \bar V^*)_2$. We
have:

\proclaim{\label Proposition}
There is a natural identification $(\g _{K}(X;\R ),h )\cong (\su (V\oplus
\bar V^*),u)$; this is a real form of the case $(A_{2n-1},A_{n-1}+A_{n-1})$.
Furthermore, the subspace  $\oplus _k\Hm ^{k,k}(X)[n]$ is a  fundamental
Jordan--Lefschetz module of $\g _{K}(X)$.
\endproclaim
\demo{Proof}
The first assertion is clear from the preceding discussion.
The choice of a generator of $\wedge _{\C}^nV$
determines an isomorphism $\wedge^{n-k}_{\C}V\cong \wedge ^k_{\C} V^*$.
This yields a graded isomorphism
$$
\wedge ^n_{\C}(V\oplus \bar V^*)=\oplus _k \wedge ^{n-k}_{\C}V\otimes\wedge ^k
_{\C}\bar V^*\cong\oplus _k\wedge ^k _{\C}V^*\otimes\wedge ^k_{\C}\bar V^*[n]
\cong \oplus _k \Hm ^{k,k}(X)[n].
$$
The last assertion follows from this.
\enddemo

\medskip
We next determine the N\'eron--Severi Lie algebra (or rather the  Lie algebra
of
its rational points) of an abelian variety $X$. We adhere to the convention to
denote the $\Q$-algebra $\End (X)\otimes\Q$ by $\End ^0(X)$ and
we write $V_{\Q}$ for $\Hm _1(X;\Q )$ and $V$ for $\Hm _1(X; \R )$. We think of
$\End ^0(X)$ as a subalgebra of $\End (V_{\Q })$;
if $J:V\to V$, $J^2=-\bold{1}_V$, is the
complex structure determined by the one of $X$, then  $\End ^0(X)$ is
centralizer of $J$ in $\End (V)$ intersected with $\End
(V_{\Q})$. Likewise, $\NS (X)\otimes\Q$ may be identified with the
$J$-invariants in $\wedge ^2V^*$ intersected with $\wedge ^2V_{\Q}^*$.

If $\kappa\in \NS (X)$ is a polarization, then taking adjoints with respect to
this form defines an (anti-)involution ${}^{\dagger}$ in $\End (V)$:
$$
\kappa (\sigma v,w)=\kappa (v, \sigma ^{\dagger}w),
$$
which preserves $\End ^0(X)$; it is called the {\it Rosati involution}
defined by $\kappa $. A well-known fact can be stated as follows:

\proclaim{\label Proposition}
If $\lambda\in\NS (X)$ then the restriction of $[e_{\lambda},f_{\kappa}]\in\gl
(\wedge ^{\bullet} V^*)$ to $V^*$ is in $\End ^0(X)$ (acting on $V^*$
contragradiently) and is invariant under ${}^{\dagger}$; this defines an
isomorphism of $\NS (X)$ onto the ${}^{\dagger}$-invariants in $\End ^0(X)$.
\endproclaim

\demo{Proof} For any bilinear form $\lambda:V\times
V\to\R$ there is a unique $\sigma\in\End (V)$ such that  $\lambda (a,b)=\kappa
(\sigma _{\lambda} a, b)$. The condition that $\lambda$ be anti-symmetric is
equivalent to that $\sigma _{\lambda}$ be ${}^{\dagger}$-invariant; the
condition that $\lambda$ be $J$-invariant to that $\sigma _{\lambda}$  be
$J$-equivariant. If $\lambda$ is skew-symmetric and is regarded as an element
of $\wedge ^2V^*$, then eq.~ \refer{**} shows that
$[e_{\lambda},f_{\kappa}]|V^*$ is equal to $\pm\sigma ^*_{\lambda}$ plus a
scalar operator. The proposition follows.
\enddemo

\label Let us write $\End ^0(X)^{\pm }$ for the $\pm 1$-eigen space of
${}^{\dagger}$ in
$\End ^0(X)$. Since ${}^{\dagger}$ is an anti-involution, $\End ^0(X)^-$ is a
Lie subalgebra of $\End ^0(X)$ and $\End ^0(X)^+$ is a module of this Lie
algebra.
The group of units $(\End (X)\otimes\R)^{\times}$ acts on $\NS (X)\otimes\R$
and it is well-known that the polarisations are contained in a single orbit. So
the Rosati involutions are all conjugate under
$(\End (X)\otimes\R)^{\times}$. Let $\uf (X)$ denote the set of
elements in $\End ^0(X)$ that are  anti-invariant with respect to all
Rosati involutions. This is clearly a Lie ideal in $\End ^0(X)$.

\proclaim{\label Proposition}
The N\'eron--Severi Lie algebra of the abelian variety $X$ is of Jordan type.
Its degree $2$ summand is canonically isomorphic to $\NS (X)\otimes\Q$. Its
degree $0$ summand can be identified with the Lie ideal of $\End ^0(X)$
that is generated by $\End ^0(X)^+$ and this
isomorphism makes $h$ correspond to a scalar operator in $\End ^0 (X)$.
Moreover, $\End ^0(X) =\g _{NS}(X;\Q )_0\times \uf (X)$.
\endproclaim

\demo{Proof} We only prove the last two assertions, as the others just sum up
the preceding discussion. For a Jordan pair $(\g ,h)$, $[\g _2,\g
_{-2}]$ generates $\g _0$ and  so \refer{3.5} implies that $\g _{NS}(X)_0$ is
the Lie algebra generated by the elements invariant under some Rosati
involution.  As the Rosati involutions are dense in a single $(\End
(X)\otimes\R)^{\times}$-conjugacy class, a standard argument
shows that $\g _{NS}(X,\R )_0$ must be the Lie ideal in $\End (X)\otimes\R$
generated by  $(\End (X)\otimes\R )^+$.

If $\sigma\in \End ^0(X)^-$ and $\tau\in \End ^0(X)^+$, then
$\Tr (\sigma\tau )=\Tr ((\sigma\tau )^{\dagger})=\Tr (-\tau\sigma )=-\Tr
(\sigma\tau )$
(here $\Tr$ denotes the $\Q$-trace). This shows that $\End ^0(X)^-$ is the
the orthoplement of  $\End ^0(X)^+$ in $\End ^0(X)$ with respect to
the trace form. So $\uf (X)$ is the orthogonal complement of the span of the
elements  fixed by some Rosati involution. This orthoplement is an ideal as
well and hence equal to  $\g _{NS}(X)_0$.
\enddemo

\label
In order to convert this into a more explicit statement, we first make a
standard reduction.

The N\'eron--Severi Lie algebra of $X$ only depends on the isogeny type. So by
the
Poincar\'e's complete reducibility theorem we may without loss of generality
assume that the abelian variety is of the form
$$
X=X_1^{m_1}\times\cdots \times X_k^{m_k}
$$
with $X_1,\dots ,X_k$ simple, pairwise non-isogenous abelian varieties. Since
the N\'eron--Severi group of $X$ is just the direct sum of the
N\'eron--Severi groups of its isotypical factors, the same is true for the
N\'eron--Severi Lie algebra:
$$
\g _{NS}(X)=\g _{NS}(X_1^{m_1})\times\cdots \times\g _{NS}(X_k^{m_k}).
$$
We therefore assume that $X$ is a power $A^m$ of a simple abelian variety
$A$.

We first concentrate on $A$.
(For a discussion and the proofs of the properties that we are going to use
we refer to \cite{Lange-Birk}.) Since $A$ is simple, $\End ^0(A)$ is a skew
field over $\Q$. We shall write $F$ for it and denote the center of $F$ by $K$.
Then $\Hm _1(A;\Q )$ is in
a natural way a $K$-vector space. We fix a polarization $\kappa$ so that there
is defined a corresponding Rosati involution ${}^{\dagger}$. This involution is
positive
in the sense that for every nonzero $g\in F$, the action of $g^{\dagger}g$ on
the $\Hm _1(A;\Q)$  has positive trace over $\Q$.
The involution that $^{\dagger}$ induces in $K$ is independent of $^{\dagger}$:
for any embedding of $K$ in the complex field it is given by complex
conjugation. We therefore denote it by $\bar{}$.
The subfield $K_0\subset K$ fixed by $\bar{}$ is totally real; so if
 $\bold{e}_0$ stands for the set of embeddings of $K_0$ in $\R$, then
 $K_0\otimes _{\Q}\R $ is as a real vector space
canonically isomorphic to $\R ^{\bold{e}_0}$. The cardinality $e_0$ of
$\bold{e}_0$ is the degree of $K_0$ over $\Q$. The isomorphism between $\NS
(A)\otimes\Q$ and the ${}^{\dagger}$-invariants $F^+$ in $F$ gives $\NS
(A)\otimes\Q$ the
structure of a $K_0$-vector space as well. This is also independent of
$\dagger$. There are four cases:

\smallskip
{\it The case of totally real multiplication.} Then $F=K=K_0$, in particular,
the involution is trivial on $F$.

\smallskip
{\it The case of totally indefinite quaternion multiplication.} Here $K=K_0$
and $F$ is a $K_0$-form of $\End (2)$: there is an $\R$-algebra isomorphism
 $F\otimes _{\Q}\R\cong\End (2,\R )^{\bold{e}_0}$ such that the involution
 corresponds
to the transpose in every summand. In particular, $F^+$ is a $K_0$-form of the
space of binary quadratic forms.

\smallskip
{\it The case of totally definite quaternion multiplication.}
Here also $K=K_0$ and $F$
is a  $K_0$-form of the quaternion algebra $\quat$ over $K_0$: there is an
$\R$-algebra isomorphism $F\otimes _{\Q}\R\to \quat ^{\bold{e}_0}$ such that
${}^{\dagger}$ corresponds to quaternion
conjugation in every summand. Since the self-conjugate elements in $\quat$ are
the reals, it follows that $F^+=K_0$.

\smallskip
{\it The case of totally complex multiplication.}  The field $K$ has no real
embedding (so is a totally imaginary quadratic extension  of $K_0$) and $F$ is
a
$K$-form of $\End (d)$: there is an $\R$-algebra isomorphism $F\otimes
_{\Q}\R\to\End (d,\C )^{\bold{e}_0}$ such that the involution corresponds to
the
conjugate transpose in every summand. So $F^+$ is a $K_0$-form of the space of
Hermitian $d\times d$-matrices.

\medskip According to Albert all these cases occur. Recall that an
isogeny type of a polarized abelian variety is given by rational vector space
$V_{\Q}$, a nondegenerate symplectic form $\kappa$ on $V_{\Q}$ and
a complex structure $J$ on the realification $V$ of $V_{\Q}$ such  that
$\kappa$
is $J$-invariant and $\kappa (a,Ja)>0$ for all nonzero $a\in V$.
In order to realize the above cases, we fix a free finitely  generated left
$F$-module $W$ and a nondegenerate skew-symmetric Hermitian form $\phi :
W\times
W\to F$ (i.e., $\phi (b,a)=\phi (a,b)^{\dagger}$ and $\phi$ $F$-linear in the
first
variable). Such a $\phi$ can be brought into a standard form: if the
involution
is trivial ($F=K=K_0$), then  $\dim _FV$ must be even, say $2r$, and there
exists a basis $(e_{\pm 1},\dots ,e_{\pm r})$ such that  $\phi (a,b)=\sum
_{i=1}^r a_ib_{-i}-a_{-i}b_i$; if it is not, then there exists a basis
$(e_1,\dots ,e_r)$ of $W$ and nonzero $u_i\in F$ with $u_i^{\dagger}=-u_i$
($i=1,\dots
,r)$ such that $\phi (a,b)=\sum _{i=1}^r a_iu_ib_i^{\dagger}$. We take for
$V_{\Q}$ the
$\Q$-vector space underlying $W$ and let $\kappa :=\Tr _{F/\Q}\phi
:V_{\Q}\times
V_{\Q}\to\Q$. Then one can find an complex structure $J$ on $V$ which commutes
with $F$, preserves $\kappa$ and is such that $\kappa (a,Ja) >0$ for all
nonzero
$a\in V$. Under some mild restrictions (given in \cite{Shimura}, \S 4, Thm.~5
ff.), one can also arrange that the centralizer of $J$ in $\End (V_{\Q})$ is no
more than $F$; this means that $F$ appears as the endomorphism algebra
tensorized with $\Q$ of any abelian variety associated to $(V_{\Q},J)$.

\medskip
Let us define a  $K_0$-Lie subalgebra of $\End (2m,F)$ by:
$$
\sku (2m,F,\dagger )=\{
\pmatrix
A& B\\
C& -{}^t\! A^{\dagger}
\endmatrix\right) | A,B,C\in\End (m,F); B={}^{t}\! B^{\dagger}, C={}^{t}\!
C^{\dagger}\}.
$$
This is the Lie algebra of infinitesimal automorphisms of the skew-hermitian
form
$$
\sum _{k=1}^m(z_kw^{\dagger}_{-k}-z_{-k}w^{\dagger}_k).
$$
It is a reductive $K_0$-Lie algebra whose center is the space of
scalars $\lambda\in K$ with $\lambda ^{\dagger}=-\lambda$. So $\sku
(2m,F,\dagger )$ is
semisimple unless we are in the case of totally complex multiplication. We
grade
this Lie algebra by means of the semisimple element
$$
u_m:=
\pmatrix
-\bold{1}_m& 0\\
0&\bold{1}_m
\endmatrix\right) \in \sku (2m,F,\dagger )
$$
so that $A$, $B$ and $C$ parametrize the summands of degree $0$, $-2$ and $2$
respectively. Let $\g (2m,F,\dagger )$ denote the $K_0$-Lie subalgebra of
$\sku (2m,F,\dagger )$  generated by the summands of degree $2$ and $-2$; let
$\uf (m,F,\dagger )$
denote the union of $\GL (m,F)$-conjugacy classes in $\End (m,F)$ made up of
anti-invariants with respect to the involution $A\mapsto {}^tA^{\dagger}$ and
identify
$\uf (m,F,\dagger )$ with a subspace of $\sku (2m,F,\dagger )_0$ in an obvious
way.

The proof of the following lemma is left to the reader.

\proclaim{\label Lemma}
The pair $(\g (2m,F,\dagger ),u_m)$ is a Jordan pair.
The space $\uf (m,F,\dagger )$ is a Lie ideal in
$\sku (2m,F,\dagger )$ supplementary to $\g (2m,F,\dagger )$. It is
trivial except in the following cases:
\roster
\item $m=1$ and $F$ is totally definite quaternion: then $\g (2,F,\dagger
)\cong
\sli (2,K_0)$ and $\uf (1, F)$ can be identified with the
${}^{\dagger}$-anti-invariants in $F$ or
\item or $K$ is totally complex: then $\g (2m,F,\dagger )$ consists of
the matrices for which $A$ has its $K$-trace in $K_0$, whereas $\uf (m, F)$ can
be identified with the purely imaginary scalars in $K$ (i.e., the $\lambda\in
K$ with $\bar\lambda =-\lambda $).
\endroster
\endproclaim

Notice that in the exceptional cases the connected Lie subgroup of
$\GL(m,F\otimes _{\Q}\R )$ with Lie algebra $\uf (m,F,\dagger )\otimes _{\Q}\R
$ is a
product of $e_0$ copies of $U(1)$ resp.\ $SU(2)$ and hence compact.

\proclaim{\label Theorem}
The graded Lie algebra $\g _{NS}(A^m;\Q )\times\uf (A^m)$ is in a natural way a
product of graded $K_0$-Lie algebra's and as such it is isomorphic
(factor by factor) to  $\g (2m,F,\dagger )\times \uf (m,F,\dagger )=\sku
(2m,F,\dagger )$.
\endproclaim
\demo{Proof} We make the identification $\Hm _1(A^m;\Q)=\Hm _1(A;\Q )^m$. This
identifies the algebra $\End (A^m)\otimes\Q$ with $\End (m, F)$. We polarize
each summand by means of $\kappa$. Then the sum of these polarizations is a
polarization of $A^m$ and the correponding Rosati involution in $\End (m, F)$
is given by $\sigma\mapsto {}^{t}\sigma ^{\dagger}$. Hence $\NS (A^m)\otimes\Q$
can be
identified with the space of $^{\dagger}$-hermitian matrices in $\End (m, F)$.
This
identifies $\g _{NS}(A^m;\Q)_2$ resp.~$\g _{NS} (A^m;\Q)_{-2}$ with  $\g
(2m,F,\dagger )_2$ resp.~$\g (2m,F,\dagger )_{-2}$. Since the N\'eron--Severi
Lie algebra is
generated by these summands, it follows that this extends to an isomorphism
of $\g _{NS}(A^m;\Q )$ onto $\g (2m,F,\dagger )$. The rest is easy.
\enddemo

We describe the situation in each of the four cases:

\smallskip
{\it Totally real multiplication.} Then $\g _{NS}(A^m;\Q )$ is a
$K_0$-form of $\sy (2m)$ and we have $\g _{NS}(A^m;\Q )_0=\End
(A^m)\otimes\Q\cong\End (m,K_0)$.

This is a $K_0$-form of the case $(C_m,A_{m-1})$.

\smallskip
{\it Totally indefinite quaternion multiplication.} Then $\g
_{NS}(A^m;\Q )$ is a $K_0$-form of $\sy (4m)$ and  $\g _{NS}(A^m;\Q )_0=\End
(A^m)\otimes\Q\cong\gl (2m;K_0)$.

This is a $K_0$-form of the case $(C_{2m},A_{2m-1})$.

\smallskip
{\it Totally definite quaternion multiplication.}
The Lie algebra $\g _{NS}(A;\Q )$ is isomorphic to $\sli (2,K_0)$ and so
$\g _{NS}(A;\Q )_0\cong K_0$, whereas $\End (A)$ is a quaternion
$K_0$-algebra. For $m\ge 2$, $\g _{NS}(A ^m;\Q )$ is a $K_0$-form of $\so
(4m)$ and  $\g _{NS}(A;\Q )_0=\End (A^m)\otimes\Q\cong \End (2m,K_0)$.

This is a $K_0$-form of the case $(D_{2m},A_{2m-1})$.

\smallskip
{\it Totally complex multiplication.} The Lie algebra $\g _{NS}(A^m;\Q )$
is a $K_0$-form of $\sli (2md)$.
The inclusion $\g _{NS}(A^m;\Q )_0\subset \End (A^m)\otimes\Q$
corresponds to $\sli (md,K)\times K_0\bold{1}_m\subset\gl (md,K)$.

This is a $K_0$-form of the case $(A_{2md-1},A_{md-1}+A_{md-1})$.

\smallskip In all these cases,  the subalgebra of $\Hm (A)$ generated by the
N\'eron--Severi group is a Jordan--Lefschetz algebra, or rather a tensor
product
of such: if $A(k)$ denotes the Jordan--Lefschetz $K_0$-algebra of level $k$
associated to $K_0$-Jordan pair $(\g _{NS}(A ^m;\Q ),h)$, then the subalgebra
of
$\Hm (A;\R )$ generated by the N\'eron--Severi group can be identified with
the tensor product of the $\R$-algebra's $\R\otimes _{\sigma}A(k)$, where
$\sigma$ runs over $\bold{e}_0$ and $k$ can be calculated from the equality
$ke_0.\depth A(1)=\dim _{\C}A^m$. In the case of complex multiplication, $k$
will be divisible by $d$.

\medskip\label
It is clear that the Hodge algebra $\Hdg (X)\subset \Hm (X)$ (i.e., the complex
span  of the rational part of $\oplus _k\Hm ^{k,k}(X)$) is $\g
_{NS}(X)$-invariant and contains the subalgebra generated by the
N\'eron--Severi
classes as a $\g _{NS}(X)$-submodule. So the Hodge conjecture for $X$ is
basically concerned with the primitive subspace $\Prim (\Hdg (X))$ (in the
sense
of \refer{1.13}) in positive cohomological degree. In this way the
N\'eron--Severi Lie algebra neatly supplements the Mumford--Tate group in
helping us (at least in principle) to understand the Hodge ring: the latter
characterizes $\Hdg (X)$ as the ring of invariants of the  Mumford--Tate group,
whereas the decomposition of $\Hdg (X)$ into $\Q$-irreducible representations
of
$\g _{NS}(X)$ tells us among other things which classes do not come from
divisors. To illustrate the point, let us observe that $\uf (X)$ kills the
N\'eron--Severi group, hence kills the subalgebra of $\Hdg (X)$
generated by this group. So as soon as $\uf (X)$ acts nontrivially
on $\Hdg (X)$, $X$ will have  Hodge classes that are not in this subalgebra.
This often happens when $\uf (X)\not= 0$ \cite{Moonen-Zah}.

\head
\section Geometric examples of Jordan type II: hyperk\"ahlerian manifolds
\endhead

\noindent\label
In the previous section we saw  that complex tori and abelian varieties furnish
examples of the Jordan--Lefschetz algebra's that are of type $(A
_{2m-1},A_{m-1}+A_{m-1})$, $(C_m,A_{m-1})$  and $(D_{2m},A_{2m-1})$. The
classical
Jordan--Lefschetz algebra's that remain are those of type $(B_m,B_{m-1})$ and
$(D_m,D_{m-1})$ and the purpose of this section is provide geometric examples
of
them. One way to get such examples is to take a compact K\"ahler surface $X$:
then
$\Hm ^{\ev}(X;\R )[2]$ is a fundamental Jordan--Lefschetz module of $\Hm
^2(X;\R )$
of the desired type.  The corresponding  Lie algebra is $\so (\phi )$, where
$\phi$
is the form defined in \refer{1.9} and its degree zero part is $\so (\Hm
^2(X;\R))\times\R h$. A more interesting  class of examples was found by  M.\
Verbitsky, and this is what we will discuss in what follows.

\medskip
Let $\quat$ be a
quaternion algebra over $\R$. We denote its trace by $\Tr :\quat\to\R$ so that
${1\over 4}\Tr$ is the projection onto the real subfield. Elements of the
kernel of $\Tr$ are called {\it pure quaternions}; they make up a Lie algebra
that we denote by $\quat _0$. The {\it pure part} of $a\in\quat$ is its
projection in  $\quat _0$: $a_0:=a-{1\over 4}\Tr (a)$.

We have also defined the norm $\Nm :\quat\to \R$, $\Nm (a)=a\bar a$ (where
$\bar{}$ is the natural anti-involution that is $-1$ on $\quat _0$).
The set $\quat _1$ of elements of norm $1$ is a Lie subgroup of the group of
units $\quat ^{\times}$ of $\quat$ and has $\quat _0$ as its Lie algebra. It is
isomorphic to $SU(2)$. The unit sphere in
$\quat _0$, $\quat _0\cap \quat _1$, is precisely the set of square roots of
$-1$ in $\quat$ and so effectively parametrizes the field  homomorphisms $\C\to
\quat$. It is a conjugacy class of $\quat ^{\times}$.

Let $T$ be a left $\quat$-module of finite rank $m\ge 1$, equipped with
positive
definite real inner product $\la\, ,\,\ra :T\times T\to\R$ that is
$\quat$-invariant (this gives rise to $\quat$-Hermitian form). We write $V$ for
its
$\R$-dual $\Hom (T,\R)$ and we let $\quat$ act on the latter on the right.

Every $J\in\quat _0\cap \quat _1$ gives $T$ the structure of a complex vector
space of dimension $2m$. Since the inner product is $\quat$-invariant,
$H_J(x,y):=\la x,y\ra -\sqrt{-1}\la Jx,y\ra$ is a $J$-Hermitian form on $V$.
Its imaginary part is antisymmetric and thus determines a 2-form
$\kappa _J\in\wedge ^2V$. Wedging with $\kappa _J$ defines an operator in
$\wedge V [2m]$ that we denote by $e_J$. It has the Lefschetz property: the
corresponding degree $-2$ operator $f_J$ is characterized by
$f_J=\star e_J\star ^{-1}$. This makes sense for any nonzero element
$a\in\quat _0$: $e_a$ has the Lefschetz property and $f_a=\Nm (a)^{-1}
\star e_a\star ^{-1}$. It is clear that the $f_a$'s commute.
We denote the Lie algebra generated by these elements by $\g (V)$.
This applies in particular to the left $\quat$-module that underlies
$\quat$ itself (with inner product given by the norm). So there is defined
a Lie algebra $\g (\quat )$.
This is actually the universal case, because an orthogonal splitting
of $T$ into $m$ $\quat$-lines allows us to identify $T$ with an orthogonal
direct sum $\quat\oplus\cdots\oplus\quat$. This induces a graded algebra
isomorphism $\wedge ^{\bullet}V\cong (\wedge ^{\bullet}_{\R}
\quat)\otimes\cdots\otimes (\wedge ^{\bullet}_{\R}\quat)$ that is
compatible with the actions of $e_a$ and $h$. The Lie algebra $\g (V)$ now
appears as $\g (\quat )$ acting on the $m$-fold tensor power of its defining
representation. In particular, we find an isomorphism of graded Lie algebra's
$\g (\quat)\cong\g (V)$ that extends the identifications between the actions
of $e_a$, $a\in\quat _0$.

The following lemma gives more information about the nature of this Lie
algebra.

\proclaim{\label Lemma}
\roster
\item"{(i)}" $(\g (\quat ), h)$ is a Jordan--Lefschetz pair with
$\g (\quat )_2$ canonically isomorphic to the vector space underlying
$\quat _0$.
\item"{(ii)}" We have a natural isomorphism $\g (\quat
)_0\cong\quat _0\times\R h$, where $\quat _0$ is regarded as the  Lie algebra
of
$\quat _1$. The given action of $\quat _1$ on $\wedge ^{\bullet}V$ integrates
the action of this summand.
\item"{(iii)}" $(\g (\quat ),h)$ is isomorphic to the orthogonal Lie
algebra defined by the form $x_1x_5+x_2^2+x_3^2+x_4^2$ with $h$ corresponding
to $\text{diag}(-1,0,0,0,1)$.
\item"{(iv)}" The subalgebra $M\subset\wedge ^{\bullet}V$ generated by
the
$\kappa _J$'s is invariant under the star operator and $\g (\quat )$, and
$M[2m]$ becomes a Jordan--Lefschetz module of $(\g (\quat ), h)$ of level $m$.
\endroster
\endproclaim
\demo{Proof}  In view of the preceding discussion, we may assume that $m=1$.
Besides the established fact that the $e_a$'s and the $f_a$'s commute  among
each other,  one verifies that
$$
[e _a, f_b]= -(ab^{-1})_0 +{1\over 4}\Tr (ab^{-1})h,
$$
where $h$ defines the grading and $(ab^{-1})_0\in\quat _0$ operates on $\wedge
^{\bullet}V$ on the right via the $\quat$-module structure on $V$. One further
checks that $e_a$ and $f_a$ commute with the action of $\quat$.  The assertions
then follow in a straightforward manner, but let us
nevertheless make some remarks that make the verification rather simple and
give
a clearer picture as well. We keep assuming that $m=1$.

The star operator $\star$ in $\wedge ^{\bullet}V$ defines an involution
in $\wedge ^2V$ whose eigen spaces we denote $(\wedge ^2V)^{\pm}$.
There are two interesting symmetric bilinear forms on $\wedge ^{\bullet}V$.
One,
denoted $\phi$, is characterized by
$$
\phi (x,y):=\int x\wedge y \quad\text{ if } \deg (x)\ge 2,
$$
where $\int :\wedge ^{\bullet}V\to\wedge ^4V\cong\R$ is the obvious projection.
It has the
property that the $e_u$'s and $h$ leave this form infinitesimally
invariant, so that $\g\subset\so (\phi )$.
The other form is the natural extension of the inner product: $\la
x,y\ra =\int (x\wedge \star y)$. So on $(\wedge ^2V)^+$ both are equal, whereas
on $(\wedge ^2V)^-$ they are opposite, and these eigen spaces are perpendicular
for both forms. The $\kappa _J$'s make up the unit sphere in the space of
selfdual $2$-forms $(\wedge ^2V)^+$. Since the wedge product of a selfdual form
and an  antiselfdual form is zero, we find that each $e_u$
annihilates $(\wedge ^2V)^-$, so that $\g (\quat )$ acts trivially on this
space.  Hence $\g (\quat )$ acts on $\wedge ^{\ev}V$ via $M=\R\oplus(\wedge
^2V)^+\oplus\wedge ^4V$. One checks that the $e _J$'s and $f_J$'s
generate all of $\so (M, \phi )$, so that we have a surjective Lie homomorphism
$\g (\quat )\to \so (M, \phi )$ (a quick way to see this is to invoke
\refer{2.1}
and our classification of Jordan--Lefschetz pairs: the former  implies that the
image is a Jordan--Lefschetz pair with $3$-dimensional degree two summand and
latter
shows that such a pair must be of type $B_2$). This  homomorphism is in fact an
isomorphism. The action on $\wedge ^{\odd} V=V\oplus \star V$ can be
understood as follows:
$V$ is after complexification no longer irreducible as a $\quat$-module: we can
write $V\otimes\C=P\oplus \bar P$ where $P$ is a
$\quat$-invariant complex plane. Such a plane is totally isotropic with
respect to the complexified inner product. The group $\quat _1$ acts on $P$
faithfully with image $SU(P)$. Now $P\oplus\star\bar P$ is
$\g (\quat)\otimes\C$-invariant and comes with a natural symplectic form. This
symplectic form is infinitesimally preserved by the $e_a$'s so that we have a
Lie homomorphism $\g (\quat)\otimes\C\to\sy (P\oplus\star\bar P)$. This is an
isomorphism (defining a spinor representation of $\g (\quat)$).
\enddemo

\label Let $(X,g)$ be a compact connected Riemann manifold of dimension $4m$.
Assume that the holonomy group at $p\in M$, $G_p\subset GL(T_pM)$, is
isomorphic to $U(m,{\quat })\subset GL(4m)$. The group  $U(m,{\quat })$ is
the group of quaternionic transformations that preserve a positive definite
sesquilinear inner product; it is a maximal compact subgroup of the complex
symplectic group of genus $2m$ (and is often denoted by $Sp (m)$, but we
avoid that notation since it may lead to confusion with the way we refer to the
symplectic groups). The centralizer of $U(m,{\quat })$ in $\End (4m)$ is the
algebra of quaternions. So the centralizer of $G_p$ in $\End (T_pM)$ is a
quaternion algebra $\quat (p)$ that gives $T_pM$ the structure of a (left)
vector
space over the skew field $\quat (p)$. The subalgebra's $\quat (p)$, $p\in M$,
make up a subbundle of $\End (TM)$ that is flat for the Levi-Civita connection
and has trivial monodromy. This allows us to identify each of them with the
algebra of flat sections of this bundle. Although that subalgebra
depends on $(M,g)$, we will somewhat ambiguously denote it by
$\quat$. So $TM$ is naturally endowed with an action of $\quat$. The metric is
an eigen tensor of this action with character the norm.

Any  $J\in\quat _0\cap\quat _1$ (notation is as above) defines an
almost-complex structure on $X$. This almost-complex structure flat with
respect
to the Levi-Civita connection and (hence) integrable. It combines with the
given metric on $X$ to a K\"ahler structure on $X$; we denote the
K\"ahler form by $\kappa _J$. In particular $J$ preserves the harmonic forms;
if we identify $\Hm (X;\R )$ with the space of harmonic forms on $X$, then this
action is just the Weil operator (which on $\Hm ^{p,q}(X,J)$ is multiplication
by $\sqrt{-1}^{-p+q}$).

The assignment $J\mapsto \kappa _J$ extends linearly to $\quat _0$ and this
extension is an isomorphism of $\quat _0$ onto a space of harmonic
$2$-forms, of which the nonzero elements are K\"ahler classes. We denote the
image of this map by $\a$. We think of $\a$ as $3$-plane in $\Hm (X;\R )$ and
call it the {\it characteristic $3$-plane} of the metric. Its nonzero elements
have the Lefschetz property, so that is defined the Lie algebra $\g (\a ,\Hm
(X;\R ))$.

\proclaim{\label Proposition}
\roster
\item"{(i)}"
There is a unique isomorphism $\g (\quat)\cong \g (\a
,\Hm (X;\R ))$ (of graded Lie algebra's) that extends the identifications
between
the actions of $e_a$ $a\in\quat$ and the semisimple element $h$ defining the
grading.
\item"{(ii)}" Under the isomorphism of (i), the action of the Lie group
$\quat _1$ on the space of harmonic forms integrates the action of the
semisimple part $\g (\quat )'_0$ of $\g (\quat )_0$ on $\Hm (X)$. This action
preserves the algebra structure on $\Hm (X)$ so that $\g (\quat )'_0$ acts on
$\Hm (X)$ by derivations.
\item"{(iii)}" The subalgebra
$A_{\a}\subset \Hm (X;\R )$ generated by $\a$ is invariant under the star
operator and $\g (\quat )$ and $A_{\a}[2m]$ is a Jordan--Lefschetz module of
$\g
(\quat )$ of level $m$.
\endroster
\endproclaim
\demo{Proof} Part (i) follows from the observation that for every $x\in X$, the
harmonic forms define a subspace of the exterior algebra of the cotangent
space of $X$ at $x$ that  is invariant under both $\star$ and cupping with
alefschetz operator $\kappa _a$ with $a\in\quat _0$ nonzero.

For (ii) we remember that every $J\in \quat _1\cap\quat _0$ acts as a Weil
operator. So $\cos\theta +J\sin\theta $ acts on $\Hm ^{p,q}(X,J)$ as
multiplication by $\exp{(-p+q)\theta}$, hence acts as an algebra automorphism.
This proves that $\quat _1$ acts by algebra automorphism.

As for (iii), note that $A_{\a}$ is additively spanned by the subalgebra's $\R
[a]$, with $a\in\a$ nonzero. As
such a subalgebra is invariant under $\star$, so is $A_{\a}$. Hence $A_{\a}$ is
also invariant under $f_a$. The rest of the assertion is clear.
\enddemo

Most of the preceding statement is due to Verbitsky (he proved a weaker form of
(ii)).

\cite{Beauville} shows that once $X$ admits one Riemann metric with
$U(m,\quat )$ as holonomy group, then it admits many of them. As we have seen
above any such metric defines a characteristic $3$-plane in $\Hm ^2(X;\R)$.
Using
a theorem of S.T.\ Yau, he proves among other things the following:  \roster
\item"{(i)}" There is a nonempty open subset of the $3$-plane Grassmannian of
$\Hm ^2(X,\R )$ parametrizing characteristic planes.
\item"{(ii)}" There is a nonzero symmetric bilinear
form $q_0$ on $\Hm ^2(X,\R )$ with the property that for every characteristic
$3$-plane $H$ defined by the metric $g$, its $g$-orthogonal complement
$H^{\perp}$ coincides with its $q_0$-orthogonal complement and $q_0$ and $g$
define
proportional forms on $H$ (with positive ratio) and $H^{\perp}$
(with  negative ratio).
\endroster
So $q_0$ is unique up to positive factor and is nondegenerate of signature
$(3,b_2(X)-3)$.

This will imply the following analogue of \refer{4.4}, which is also  mostly
due
to \cite{Verbitsky 1995}.

\proclaim{\label Proposition} We then have:
\roster
\item"{(i)}" The pair $(\g _{\tot}(X;\R),h)$ is of Jordan--Lefschetz type of
type $(B,B)$ or $(D,D)$ with $\g _{\tot}(X;\R)$ isomorphic to $\so (4,
b_2(X)-2)$.
\item"{(ii)}"  We have natural identifications $\g _{\tot}(X;\R)_2\cong \Hm
^2(X;\R
)$ and $\g _{\tot}(X;\R)_0\cong \so (q_0)\times\R h$.
The semisimple part of $\g _{\tot}(X;\R)_0$, $\g _{\tot}(X;\R)'_0\cong\so
(q_0)$, acts on $\Hm (X;\R )$ by derivations.
\item"{(iii)}" The subalgebra $A$ of $\Hm (X;\R )$ generated by $\Hm ^2(X;\R )$
is
invariant under $\g _{\tot}(X;\R)$ and $A[2m]$ is a Jordan--Lefschetz module of
$\g _{\tot}(X;\R)$ of level $m$.
\endroster
\endproclaim

\demo{Proof} According to the previous proposition, the operators $f_a$ and
$f_b$
commute when $a$ and $b$ are nonzero elements of a characteristic $3$-plane. By
(i) this is therefore the case for $a$ and $b$ in a nonempty open subset of
$\Hm ^2(X;\R )$. Since the expression $[f_a,f_b]$ is rationally dependent on
its
arguments, it follows that $f_a$ and $f_b$ commute whenever both are defined.

Let $\g _2$ resp. $\g _{-2}$ denote the abelian Lie subalgebra's of
$\gl (\Hm (X;\R ))$ spanned by the $e_a$'s resp. $f_a$'s and let
let $\g _0$ be the Lie subalgebra generated by $[\g _2,\g _{-2}]$.

\smallskip
{\it Claim 1.} $\g _0 =\g _0'\times\R h$ (with $\g _0 '=[\g _0 ,\g _0 ]$) and
$\g _0'$ consists of derivations of $\Hm (X;\R )$. Moreover, the image of $\g
_0'$
in $\gl (\Hm ^2(X,\R ))$ is $\so (q_0)$.\par

Proof. Notice that by the above (rationality) argument,
$\g _0$ is already  generated by the brackets $[e_a,f _b]$ with $a,b$ nonzero
and contained in a characteristic $3$-plane. If we fix such a $3$-plane $\a$,
then  the Lie subalgebra of $\g _0$  generated by the
$[e_a,f_b]$, $a,b$ nonzero elements of $\a$, is $\g (\a ,\Hm (X;\R ))'_0\times
\R h$ with $\g (\a ,\Hm (X;\R ))'_0$ acting as infinitesimal algebra
automorphisms of $\Hm (X)$, that is, as derivations. Moreover
$\g (\a ,\Hm (X;\R ))'_0$ leaves $\q _0$ invariant. So $\g '_0$ acts as
derivations and maps naturally to $\so (q_0)$. This last homomorphism is
surjective, because $\so (q_0)$ is generated by its elements that kill the
$q_0$-orthogonal complement of a characteristic $3$-plane.

\smallskip
{\it Claim 2.} $\ad _{\g _0}$ leaves $\g _2$ and $\g _{-2}$ invariant.\par

Proof. Let $u\in \g '_0$. Since $u$ is a derivation,
we have for every $a\in \Hm ^2(X;\R )$ and $z\in \Hm (X;\R )$, that
$[u,e_a](z)=u(a.z)-a.u(z)= u(a).z= e_{u(a)}(z)$. So $\ad _u$ leaves the space
of operators $e_a$ invariant. If $G_0'\subset GL(\Hm (X;\R ))$ denote the
connected
closed subgroup with Lie algebra $\g '_0$, then for every $g\in G_0'$ we have
$ge_ag^{-1}=e_{g(a)}$ and $ghg^{-1}=h$. Hence if $f_a$ is defined, then
$gf_ag^{-1}=f_{g(a)}$. It follows that $\Ad _{G '_0}$ leaves $\g _{-2}$
invariant. The same assertion then holds for $\ad _{\g '_0}$.

\smallskip
We conclude that $\g :=\g _{-2}+\g _0 +\g _{-2}$ is a Lie subalgebra of
$\gl (\Hm (X;\R ))$ and hence equal to $\g _{\tot}(X;\R)$.
Clearly, $(\g ,h)$ is a Jordan--Lefschetz pair.

\smallskip
{\it Claim 3.} $A$ is an irreducible $\g$-submodule that has $1$ as lowest
weight
vector and $\g$ acts faithfully on $A$. Moreover, $\g _0'$ maps isomorphically
onto $\so (q_0)$.\par

Proof. From $U\g =U\g _2.U\g _0.U\g _{-2}$ and and the fact that $\g _0+\g
_{-2}$ stabilizes the unit element it follows that $A$ is $\g$-invariant and
has
$1$ as lowest weight vector.  It is clear that the kernel of the homomorphism
$\g\to\gl (M)$ is contained
 in $\g '_0$. This kernel acts trivially on $\Hm ^2(X;\R )$, hence has zero Lie
bracket with any $e_a$ (here we use that $\g _0'$ acts by derivations).
Applying the Jacobson--Morozov theorem to $e_a$ acting  on $\g$, we also find
that the kernel has zero Lie bracket with $f_a$, when defined. So the kernel
has
zero Lie bracket with $\g _2$, $\g _{-2}$ and hence also with $\g _0=[\g _2,\g
_{-2}]$. Since $\g$ is semisimple, this implies that the kernel is trivial.

We have already seen that the map $\g _0'\to\so (q_0)$ is surjective. Since
$\g _0'$ acts by derivations, its action on $A$ is completely determined by its
restriction to $\Hm ^2(M,\R )$. Hence it is also injective.

\smallskip
The theorem now follows easily.
\enddemo

\proclaim{\label Corollary}
The algebra structure on $\Hm (X;\Q)$ and the Hodge structure on $\Hm ^2(X)$
determine the Hodge  structure on all of $\Hm (X)$.
\endproclaim
\demo{Proof} Let $\tilde G$ be the closed connected $\Q$-subgroup of $GL (\Hm
(X))$ with Lie algebra $\g _{tot}(X)'_0$. Then its image $G$ in $GL(\Hm ^2(X))$
is an orthogonal group of rank $\ge 3$ and the projection $\tilde G\to G$ has
finite kernel. The Hodge structure on $\Hm ^2(X)$ is given by a real
representation of the Deligne torus $\bS$ on $\Hm ^2(X)$. This defines an
$\R$-morphism of algebraic groups $\bS\to G$. Similarly, the Hodge structure on
$\Hm (X)$ is given by an $\R$-morphism of algebraic groups $\bS\to\tilde G$.
The
latter must be a lift of the former. But clearly such a lift is unique.
\enddemo

\label{\it Example.} Let $T$ be a complex torus of complex dimension $2$ and
let $m$ be an positive integer. As Beauville explains, the $(m+1)$-fold
symmetric product $S^{m+1}(T)$ of $T$ admits a natural nonsingular resolution
$T^{[m+1]}\to S^{m+1}(T)$  (the Hilbert scheme paramatrizing finite subschemes
of $T$ of length $m+1$). Composition of this resolution with the  ``sum map''
$S^{m+1}(T)\to T$ gives a fibration  $T^{[m+1]}\to T$ which is locally trivial
in the \'etale sense. Let $K_m$ be the fiber over the origin. This fiber is
simply connected and if $m\ge 2$, then $\Hm ^2(K_m)$ (with its Hodge structure)
is canonically isomorphic to the direct sum of $\Hm ^2(T)$ and the span of the
class of the exceptional divisor restricted to $K_m$. So in that case, $\dim
\Hm ^2(K_m)=7$. Beauville shows that for a generic choice of $T$, $K_m$ admits
a quaternion K\"ahler metric. Hence \refer{4.5} implies that
$\g _{\tot}(K_m;\R )\cong \so (4,5)$.
The cohomology of $K_2$ is computed in \cite{Salamon}. We can interpret his
result as saying that as a $\g _{\tot}(K_m;\R )$-module, $\Hm (X)$ is the
orthogonal direct sum of the subalgebra generated by $\Hm ^2(K_2)$, a trivial
representation of dimension $80$ and the spinor representation (of dimension
$16$). Their Hodge polynomials are
$1+(s^2+5st+t^2)+(s^4+5s^3t+16s^2t^2+5st^3+t^4)+(s^4t^2+5s^3t^3+s^2t^4)+s^4t^4$,
$80s^2t^2$ and $(s^2t+st^2)+(s^4t+s^3t^2)$ respectively.

\medskip
We have already looked at algebra's such as $A$ in \refer{2.14}. Invoking
\refer{2.14} we find that if $u\in A_4\subset\Hm ^4(X)$ represents the dual of
the
Beauville-Bogomolov form $q_0$, then  $a\mapsto\int _Xa^2u ^{m-1}$ is
proportional to $q_0$. There is a natural choice for $u$:

\proclaim{\label Theorem}
Let $p\in \Hm ^4(X;\R )$ be the image of the first Pontryagin class of $X$
under
orthogonal projection of $\Hm (X;\R)$ onto $A$. Then the form  $r$ on $\Hm
^2(X;\R )$ defined by $r(a)=\int _X a^2p^{m-1}$ is a nonzero multiple of
Beauville-Bogomolov form.
\endproclaim

We first show:

\proclaim{\label Proposition}
The Lie algebra $\g _{\tot}(X;\R)'_0$ acts trivially on the
subalgebra $P(X)$ generated by the Pontryagin classes.
\endproclaim
\demo{Proof} Let $z$ be an element in this subalgebra of degree $4k$.
Choose a metric $g$ with holonomy group $\cong U(m,\quat )$. This determines an
action of $\quat$ on $\Hm (X;\R )$ and a characteristic $3$-plane $\a\subset
\Hm ^2(X;\R )$. If $J\in\quat _0\cap\quat _1$, then with respect to the complex
structure defined by $J$, $P^{4k}(X)$ consists of classes of bidegree
$(2k,2k)$.
Equivalently: the Weil operator $J$ leaves $P(X)$ invariant. As this is true
for all $J\in \quat _0\cap\quat _1$, it follows that $P(X)$ is invariant under
all of $\quat _1$. Hence $P(X)$ is killed by the Lie algebra $\quat _0$. Since
these Lie algebra's generate the semisimple part of
$\g _{\tot}(X;\R )_0$, the proposition follows.
\enddemo

\demo{Proof of \refer{4.8}}
Denote by $\pi :\Hm (X;\R )\to M$ the $\phi$-orthogonal projection onto $M$. So
if
$a\in \Hm ^2(X;\R )$,  then $\pi (a)\in M$ is characterized by the property
that
for all  $b\in \Hm (X;\R )$ we have $\int _X ab=\int _X \pi (a)b$. This is a
$\g
_{\tot}(A;\R )$-equivariant projection. Since $p_1(X)$ is $\g _{\tot}(M;\R
)$-invariant, $p=\pi (p_1(X))$ must be a multiple of $u$ and so all we need to
see
is that $p\not=0$. But this follows from a theorem of \cite{Chen-Ogiue} which
states that for some K\"ahler class $\kappa $, $\int _X \kappa ^{2m-2}p_1(X)$
is
positive.
\enddemo

\medskip
\label We have now seen that all the classical Jordan--Lefschetz algebra's
arise
geometrically. The question comes up whether the same is true for the
exceptional
case of type $E_7$ (that corresponds to a Jordan algebra  of dimension $27$).
In
the topological setting the answer is yes: if the $27$-dimensional vector space
$W$
in \refer{2.11} is equipped with an integral structure $W_{\Z}$ such that the
cubic
form only takes integral even values (this is indeed possible), then it follows
from theorems of Wall and Jupp that there is a simply connected closed oriented
$6$-manifold $X$ for which the integral cohomology ring is isomorphic to  the
corresponding integral algebra $A_{\Z}$ (see \cite{Okon-vdVen} for a general
discussion). Now $H^{\ev}(X;\Z )$, does not change if take a connected sum of
$X$
with a number of copies of $S^3\times S^3$. We wonder whether such a manifold
admits a complex structure. Since all our interesting examples have trivial
canonical bundle we are inclined to make this question more specific
by asking:

\smallskip{\it Question.} Does there exist a Calabi-Yau $3$-fold with Picard
group of rank $27$ such that its N\'eron--Severi Lie algebra is of type $E_7$?

\smallskip We notice that the mirror dual family of such a Calabi-Yau manifold
will, if it exists,  have its period mapping take values in a Hermitian domain
of type $E_7$.

\head
\section Filtered Lefschetz modules
\endhead

\noindent\label We recall a version of the
Jacobson--Morozov lemma. If $e$ is a nilpotent transformation in a vector space
$M$, then there is a unique
nonincreasing filtration $W^{\bullet}$ preserved by $e$ such that $e$ has
the Lefschetz property in $\Gr _W^{\bullet}(M)$. Any $\slt$-triple $(e,h,f)$
containing $e$ descends to an $\slt$-triple in $\Gr _W(M)$ and splits the
filtration (so the $k$-eigen space of $h$ is a supplement of $W^{k+1}$ in
$W^k$). We shall refer to $W^{\bullet}$ as the {\it Lefschetz filtration} of
$e$.

\proclaim{\label Lemma}
Let $\g$ be a reductive Lie algebra, $\s\subset\g$ a commutative subalgebra
consisting of semisimple elements and $\chi\in\s ^*$ a character of
$\s$ in $\g$. Then for every nilpotent $e\in\g ^{\chi}$ there exists a $f\in\g
^{-\chi}$ such that $(e, [e,f],f)$ is an $\sli (2)$-triple.
\endproclaim
\demo{Proof}
Choose an $\slt$-triple $(e,h',f')$ containing $e$. Let $f''$ be the
$\g ^{-\chi}$-component of $f'$. Then $h:=[e,f'']$ is the $\g ^0$-component of
$h'$ and so $[h,e]$ is the $\g ^{\chi}$-component of $[h',e]=2e$ and hence
equal
to $2e$. Since $h$ is in the image of $\ad (e)$,  the pair $(e,h)$ is
by \cite{Bourbaki}, Ch.~VIII, \S 11, Lemme $6$, extendable to a $\slt$-triple
$(e,h,f)$. This remains an $\slt$-triple if we replace $f$ by its $\g
^{-\chi}$-component and so the lemma follows.
\enddemo

Let $(\a ,M)$ be a Lefschetz module and $\hor ^{\bullet}M$ a nonincreasing
filtration on the graded vector space underlying $M$ (so $\hor ^kM=\oplus
_l\hor ^kM_l$) which is preserved by $\a$.  We shall refer to this filtration
as the {\it horizontal filtration}.  Then the associated {\it
vertical filtration} is defined by  $\ver _kM:=\sum _r \hor ^{r-k}M_r$. This
filtration is nondecreasing; we call the corresponding grading of $\Gr
_{\ver }M$ the {\it vertical grading}. The notations $\Gr _{\hor}M$ and  $\Gr
^{\ver}M$ refer to the same vector  space but with different gradings:
$\Gr _{\hor}^kM_r=\Gr ^{\ver}_{r-k}M_r$ has horizontal degree $k$ and vertical
degree $r-k$. The notation $\Gr M$ refers to their common bigraded structure.

Suppose now that some $a\in\a$ preserves the vertical
grading (i.e., $e_a (\hor ^kM)\subset \hor ^{k+2}M$ for all
$k$) and has the Lefschetz property in $\Gr _{\hor}M$:
$e_a ^k$ sends $\Gr _{\hor}^{-k}M$ isomorphically onto $\Gr _{\hor}^kM$.
It is then immediate that
the horizontal filtration is the Lefschetz filtration of the transformation
$e_a$ in $M$. If we apply \refer{5.2} to $e:=e_a$ and $\s:=\C h$, we find an
$\slt$-triple $(e_a,h_{\hor},f_a)$ in $\g (\a ,M)$ with $f_a$ of total degree
$-2$ and $h_{\hor}$ of total degree $0$. So $f_a$ will map $\hor ^kM_r$ to
$\hor
^{k-2}M_{r-2}$. This shows that $f_a$ preserves the vertical filtration. It
also
follows that $h_{\hor}=[e_a,f_a]$ has this property. It is clear that the eigen
spaces of $h_{\hor}$ split the horizontal filtration. This element commutes
with $h$, so if we put $h_{\ver}:=h-h_{\hor}$, then the eigen spaces of the
commuting pair $(h_{\hor},h_{\ver})$ define a bigrading of $M$ that identifies
$M$ with $\Gr M$. The eigen spaces of $(h_{\hor},h_{\ver})$ under the adjoint
representation also define a bigrading of $\g (\a ,M)$. The image of $\a$ in
$\g
(\a ,M)$ need not be bigraded.

\proclaim{\label Proposition} Let $\a _{\hor}\subset\a$ be the set of $a\in\a $
that preserve the vertical grading and suppose that $\Gr _{\hor}M$ is a
Lefschetz module of $\a _{\hor}$. Then we can write $h=h_{\hor}+h_{\ver}$
with $h_{\hor}$ and $h_{\ver}$ semisimple elements of $\g (\a, M)$ that have
integral eigen values and commute with each other (so for the resulting
bigrading of $M$, $M_{k,l}$ gets identified with $\Gr _{\hor}^kM_{k+l}$).

The span of the components of $\a _{\hor}$ of lowest horizontal degree $2$
make up an abelian subalgebra $\a _{2,0}$ of $\g (\a ,M)_{2,0}$ that has the
Lefschetz property in $M$ with respect to the horizontal grading. Moreover,
$\g (\a _{2,0},M_{\hor})$ is a subalgebra of $\g (\a ,M)_{\bullet ,0}$ that
maps isomorphically onto $\g (\a _{\hor},\Gr _{\hor}M)$.

If in addition, $\Gr ^{\ver}M$ is a Lefschetz module of $\a$, then
the span of the components of $\a$ of highest vertical degree $2$ make up an
abelian subalgebra $\a _{0,2}$ of $\g (\a ,M)_{0,2}$ that has the Lefschetz
property in $M$ with respect to the vertical grading. Moreover,  $\g (\a
_{0,2},M_{\ver})$ is a subalgebra of $\g (\a ,M)_{0,\bullet}$ that maps
isomorphically onto $\g (\a ,\Gr ^{\ver}M)$. The obvious map $\g (\a
_{2,0},M_{\hor})\times \g (\a _{0,2},M_{\ver})\to\g (\a ,M)$ is then an
injective homomorphism of Lie algebra's.
\endproclaim

\demo{Proof}
Everything follows from the preceding or is obvious except the very last
statement.
We claim that any ``horizontal'' $\slt$-triple $(e',h',e')$ acting in $\Gr M$
commutes
with any ``vertical'' $\slt$-triple $(e'',h'',e'')$ acting in $\Gr M$. This
just
follows from the fact that $(e',h')$ commutes with $(e'',h'')$ and the fact
that
in either case the last member is a rational expression in the first two. It
follows that $\g (\a _{\hor},\Gr _{\hor}M)$ and $\g (\a ,\Gr ^{\ver}M)$
commute.
The same is therefore true for their bigraded lifts $\g (\a _{2,0},M)$ and
$g (\a _{0,2},M)$ in $\g (\a ,M)$. To see that $\g (\a _{2,0},M)\cap \g (\a
_{0,2},M)=0$, note that this intersection has bidegree $(0,0)$ and is normal
in either of them. If it were nonzero, then it would contain a simple factor of
$\g (\a _{2,0},M)$ of bidegree $(0,0)$. But this is impossible since $\g (\a
_{2,0},M)$ is (as a Lie algebra) generated by its degree $\pm 2$ summands. The
proposition follows.
\enddemo

Let $f:X\to Y$ be a fibration of projective manifolds which is topologically
locally trivial and let $n$ and $m$ be the complex dimensions of $X$ and $Y$
repectively, so that $d:=n-m$ is the complex fiber dimension. Following
\cite{Deligne 1968} the Leray spectral sequence of $f$ degenerates: if
$L^
{\bullet}$ denotes the Leray filtration of $\Hm ^{\bullet}(X)$, then $\Gr
^k_L\Hm ^r(X)\cong \Hm ^k(Y,R^{r-k}f_*\C )$.

\proclaim{\label Proposition}
If $h$ denotes the basic semisimple element of $\g _{NS}(X)$, then we can write
$h=h_{\hor}+h_{\ver}$ with $h_{\hor}$ and $h_{\ver}$
semisimple elements of $\g _{NS}(X)$ that have integral eigen values and
commute
with each other so that for the resulting bigradings of $\Hm (X)$ and $\g
_{NS}(X)$
have the following properties: $\Hm (X)_{k,l}$ gets identified with
$\Hm ^{k+m}(Y,R^{l+d}f_*\C
)$ and
$$
\g (\NS (Y),\Hm ^{\bullet}(Y,Rf_*\C )[m])\times\g (\NS
(X/Y),\Hm (Y,R^{\bullet}f_*\C [d]))
$$
lifts (uniquely) to a
bigraded Lie subalgebra of $\g _{NS}(X)$.
\endproclaim
\demo{Proof} We prove the proposition by verifying
the hypotheses of the previous proposition for $M:=\Hm ^{\bullet}(X)[n]$ with
as
horizontal filtration the Leray filtration appropriately shifted: $\hor
^kM_r:=L^{d+k}\Hm ^{n+r}(M)$.

First recall that $R^lf_*\C$ underlies a variation of Hodge structure of weight
$l$. If $\xi \in \NS (X)$ is
ample relative $f$, then its image in $\Hm ^0(Y,R^2f_*\C)$ has the
Lefschetz property in the graded local system $R^{\bullet}f_*\C$ and induces
a polarization in each summand. This implies that $\xi$ has the
Lefschetz property in $\Hm ^k(Y,R^{\bullet}f_*\C [d])$ and satisfies the
hypotheses of \refer{1.6}.
So $\Hm ^k(Y,R^{\bullet}f_*\C )[d]$ is a Lefschetz module over $\NS (X/Y)$.

If $\eta\in\NS (Y)$ is a polarization, then cupping with $\eta ^k$ defines an
isomorphism $\Hm ^{m-k}(Y,R^lf_*\C )\to \Hm ^{m+k}(Y,R^lf_*\C )$ \cite{Zucker},
\cite{Saito} and $\eta\cup$ has the Lefschetz property in $\Hm
^{\bullet}(Y,R^lf_*\C )[m]$. We apply \refer{1.6} again and find that all the
hypotheses of \refer{5.3} are fulfilled. \enddemo

{\it Remarks.} This splitting of the Leray filtration is certainly a splitting
that is invariant under the action of $\NS (Y)$. We do not know however whether
it can be chosen to be a splitting of $\Hm (Y)$-modules.

If the graded local system $R^{\bullet}f_*\C$ is trivial (which is the case
when $Y$ is simply connected), then by the universal coefficient theorem, $\Hm
^k(Y,R^{\bullet}f_*\C )\cong \Hm ^k(Y)\otimes \Hm ^{\bullet}(X_y)$, where $X_y$
is a fiber. If we take $y$ sufficiently general, then $\NS (X/Y)$ restricts
isomorphically to $\NS (X_y)$ and then $\g (\NS (X/Y),\Hm (Y,R^{\bullet}f_*\C
[d]))$ can be identified (via the preceding  isomorphism) with $\g
_{NS}(X_y)$.

\medskip\label
The conditions $X$ and $Y$  nonsingular and $f$ topologically
locally trivial can all be eliminated in the context
of Hodge modules: if $f:X\to Y$ is a morphism of projective
varieties, and $E$ is a polarized Hodge module on $X$ of pure weight, then
according to \cite{Saito} the Leray spectral sequence for $f_*E$  degenerates
and the Leray filtration satisfies all the hypotheses of \refer{5.3} (here no
shifting is necessary). If the map from a space $Z$ to a  fixed singleton is
denoted $a_Z$, then we find that $\g (\NS (X),a_{X *}(E))$ contains graded Lie
subalgebra's isomorphic to $\g (\NS (Y),a_{Y*}\Hm ^{\bullet}(f_*E))$ and $\g
(\NS
(X/Y),a_{Y*}\Hm ^{\bullet}(f_*E))$ that centralize each other. (Here the
cohomology
is taken in the sense of Hodge modules; for a representing complex of
constructible sheaves, this amounts to taking perverse cohomology.)

\proclaim{\label Proposition}
Let $f:X\to Y$ be a (topologically locally trivial) fibration of projective
manifolds with fiber $\PP ^d$. Then $\g _{NS}(X)$ contains a Lie subalgebra
isomorphic to $\g _{NS}(Y)\times\slt$. If this subalgebra is  equal to
$\g _{NS}(X)$, then the characteristic classes of this bundle are trivial,
i.e.,
$\Hm (X)=\Hm (Y)\otimes \Hm (\PP ^n)$ as $\Hm (Y)$-algebra's.
\endproclaim
\demo{Proof} As an algebra, $\Hm (X)$ is a simple integral extension
of $\Hm (Y)$: $\Hm (X)=\Hm (Y)[\xi ]/(P)$ with $P$ a monic polynomial of degre
$d+1$:
$P=\xi ^{d+1} +c_1\xi  ^d+\cdots +c_{d+1}$. We make $P$
unique by requiring that $c_1=0$; then $c_2,\dots ,c_{d+1}$ are the
characteristic classes of our bundle. Here $\xi $ is any element of $\NS
(X)\otimes\Q $ that spans a supplement of $\NS (Y)\otimes\Q$ in $\NS
(X)\otimes\Q $. Take a bigrading on $\Hm (X)$ as in the previous proposition so
that we get an embedding of $\g _{NS}(Y)\times\slt$ in
$\g _{NS}(X)$. If it is surjective, then take for $\xi $ a nonzero element
of $\NS (X)$ that projects in $\Hm ^2(X)_{(-\dim Y,-d+2)}$. Cupping with $\xi $
then corresponds to a nonzero element of the $\slt$ factor. So we have $\xi
^{d+1}=0$. Hence all the $c_i$'s are zero. \enddemo

\proclaim{\label Theorem}
Let $f:X\to Y$ be a $\PP ^d$-bundle involving projective manifolds. Assume that
the N\'eron--Severi Lie algebra of $Y$ is maximal: $\g _{NS}(Y)=\aut (\Hm
(Y))$,
and
that $\Hm (Y)$ is not an inner product space of dimension $4$. Then either the
characteristic classes of $f$ are trivial and  $\g _{NS}(X)\cong \g
_{NS}(Y)\times\slt$ or the N\'eron--Severi Lie algebra of $X$ is maximal.
\endproclaim
\demo{Proof} Suppose not all the characteristic classes are
trivial and let $\xi \in\NS (X)\otimes\Q$ be such that
$\Hm (X)=\Hm (Y)[\xi  ]/(P)$ with $P=\xi ^{d+1} +c_2\xi ^{d-1}+\cdots +c_{d+1}$
as in the proof of \refer{5.6}. Then
$\g _{NS}(X)$ contains a copy of $\g _{NS}(Y)\times\slt$, but $\xi $,
viewed as an element of $\g _{NS}(Y)$, commutes with neither factor. So the
simple component of $\g _{NS}(X)$ that contains $\xi $ contains $\g
_{NS}(Y)\times\slt$ as well. It remains to apply the theorem of the appendix.
\enddemo

\proclaim{\label Theorem}
Let $X$ be the flag space of a simple complex algebraic group. Then its
N\'eron--Severi Lie algebra is equal to $\aut (\Hm (X))$.
\endproclaim
\demo{Proof} In case the group is of rank one, then $X=\PP ^1$ and then the
assertion is clear. So assume that the rank is $\ge 2$. Then $X$ admits an
iterated  fibered structure
$$
X=X_s\to X_{s-1}\to\cdots \to X_0
$$
with $X_0$ a singleton, $X_t\to X_{t-1}$ a projective space  bundle with
positive fiber dimension, and $s\ge 2$. Iterated application of \refer{5.6}
yields an embedding of $(\slt )^s$ in $\g _{NS}(X)$. In view of the theorem of
the appendix it is therefore enough to show that $\g _{NS}(X)$ is simple. Were
that not the case then we could write nontrivially $\g _{NS}(X)=\s _1\times \s
_2$  with $\s _i$ nonzero semisimple. In that case, let $H_i$ be the $\s
_i$-submodule of $\Hm (X)$ generated by the unit element. As $\Hm (X)$ is
generated
by $\NS (X)$ it follows that  $\Hm (X)=H_1\otimes H_2$ as algebra's. However,
such
a decomposition is precluded by Borel's description of $\Hm (X)$. According to
this theory, there is up to scalar a unique quadratic form on $\NS (X)\otimes\C
$
that becomes a relation for $\Hm ^4(X)$. This form is nondegenerate, and  so
$\Hm (X)$ cannot be the tensor product of two graded subalgebra's (see also
\refer{1.2}).
\enddemo

\label The Borel description of $\Hm (X)$ actually shows that as an algebra,
$\Hm (X)$ is isomorphic to $\Sym (V)/I$, where $V$ is the complexified weight
lattice of $G$ and $I$ is the ideal generated by the Weyl group invariant
homogeneous  forms
of positive degree. In other words, we are in the case of the example
discussed
in \refer{1.10} and the above theorem gives a complete description of $\g
(V,\Sym (V)/I)$ in case $W$ is an irreducible Weyl group.

Another interesting class of projective manifolds with the property that their
rational cohohomology is generated by the N\'eron--Severi group are the
Knudsen--Mumford moduli spaces $\Cal{M}_0^n$ of stable $n$-pointed curves of
genus zero. It it likely that here also the N\'eron--Severi Lie algebra of
$\Cal{M}_0^n$ equals $\aut (\Hm (\Cal{M}_0^n))$ (compare theorem \refer{6.8}
below).

\head
\section Frobenius--Lefschetz modules
\endhead

\noindent\label
We say that a Lefschetz module $(M,\a )$ of depth $n$ is {\it Frobenius} if it
satisfies the following three properties:
\roster
\item"{(1)}" $\Prim M_{-n}=1$ is of dimension one (and so $M$ is
irreducible),
\item"{(2)}" the map $\a\otimes M_{-n}\to M_{-n+2}$ is an
isomorphism,
\item"{(3)}" $M$ is generated as a $U\a$-module by $M_{-n}$.
\endroster

If only the first two conditions are satisfied we say that $M$ is a {\it
quasi-Frobenius} of depth $n$ and if instead of (3), we have
\roster
\item"{(3$'$)}" the $U\a$-module generated by $M_{-n}$ contains $M_{-n+2k}$ for
$k\le d$,
\endroster
then we say that $M$ is {\it Frobenius up to order $d$}.

\smallskip
Observe that if $A$ is Lefschetz algebra of depth $n$, then $A[n]$ is a
Frobenius--Lefschetz module of $A_2$ if and only if $A$ is generated by $A_2$.
Moreover, any Frobenius--Lefschetz module is of this form.

We also note that a Jordan--Lefschetz module is Frobenius. The property
of being quasi-Frobenius is a useful one, as it turns out to be a rather
strong  approximation to being Frobenius, that (somewhat in contrast to the
latter) is generally easy to verify in practice. Another reason of our interest
in this notion is that it occurs naturally in geometric examples:

\proclaim{\label Proposition} Let $X$ be a connected compact K\"ahler
(resp. complex projective) manifold. Then the $\g _K(X)$-submodule
(resp. $\g _{NS}(X)$-submodule) of $\Hm (X)$ generated by $\Hm ^0(X)$ is
quasi-Frobenius
 as a Lefschetz module of $\Hm ^{1,1}(X)$ (resp. $\NS (X)$).
\endproclaim
\demo{Proof}
In the K\"ahler case, it is clear that this submodule is contained in  $\oplus
_k\Hm ^{k,k}(X)$. So its intersection with $\Hm ^2(X)$ is $\Hm ^{1,1}(X)$ and
the
proposition follows. The N\'eron--Severi case is proved similarly.
\enddemo

The following proposition explains our terminology for it shows that a
Lefschetz module $(M,\a )$ is Frobenius if and only if $M$ is
Frobenius ($=$ Gorenstein) as a $U\a$-module.

\proclaim{\label Proposition}
Let $(\a ,M)$ be a Frobenius--Lefschetz module of depth $n$.\break
Then there exists
a nondegenerate $\g (\a ,M)$-invariant $(-)^n$-symmetric bilinear form on $M$.
\endproclaim
\demo{Proof} Since $\dim M_{-n}=1$, we also have $\dim M_n=1$.
The choice of a generator $u\in M_{-n}$ identifies $M$ with a
graded quotient $R$ of the symmetric algebra $U\a$ (with a shift of degree).
Pick a nonzero linear form $\int :R_{2n}\to\C$, and define a graded bilinear
form
$\langle\, ,\,\rangle :M\times M\to\C$ by $\langle au, bu\rangle =(-1)^k\int
(abu)$ if $a\in R^{2k}$ $b\in R^{2n-2k}$. This form
is symmetric or skew according to whether $n$ is even or odd.
We claim that it is $\g (\a, M)$-invariant. For if we regard
$\langle\, ,\,\rangle$ as
an element of the Lefschetz module $(M\otimes M)^*$, then it is of degree
zero and annihilated by $\a$. So if $f_a$ ($a\in \a$) is
defined, then $\langle\, ,\,\rangle$ is primitive for the  $\slt$-triple
$(e_a,h,f_a)$ and hence annihilated by $f_a$. This proves our claim.
Since $M$ is irreducible, $\langle\, ,\,\rangle$ must be nondegenerate.
\enddemo

\medskip
\label {\it Example.} The following example is not just instructive, it also
will be used in  a proof (of theorem \refer{6.8}).
We give $V(2l)=K[e]/(e^{2l+1})[2l]$ the unique $\slt$-invariant quadratic form
for which the inner product of $1$ and $e^{2l}$ is $(-1)^l$, so the inner
product of $e^p$ and $e^q$ is $(-1)^{l+p}$ if $q=l-p\in\{ 0,\dots ,l\}$ and
zero otherwise.  Let $V:=V(2k)\oplus V(2k-2)$, $k\ge 2$, regarded as the
orthogonal direct sum of $K[e]$-modules.  Any endomorphism of $V$ that
commutes with $e$ is  $K[e]$-linear and so representable by a $2\times 2$
matrix with coefficients in $K[e]$. We readily compute that the intersection
of the centralizer of $e$ with  $\so (V)_2$ is the set of matrices of the form
$$
\pmatrix
ae & be^2\\
b & ce\\
\endpmatrix
$$
with $a,b,c$ scalars. These do not mutually commute, so any
subspace $\a\subset\so (V)_2$ that defines Lefschetz module structure on $V$
and contains $e$ is of $\dim \le 2$. Hence $\a$ will be spanned
by $e$ and an element $e'$ that can be represented
by a matrix of the above type with $c=-a$. If $a\not=0=b$, then it is easy to
see $\g (\a ,V)\cong\sli (2)\times\sli (2)$. Suppose therefore that $b\not=0$.
We normalize $e'$ to make $b=1$.
Then $e$ and $e'$ satisfy: $(e')^2= (a^2+1)e^2$, $e^{2k+1}=0$ and
$e^{2k-1}(e'-ae)=0$. This suggests to take as a new basis for $\a$
$a_{\pm}:=e'\pm\sqrt{a^2+1}e$. Then the relations become:
$a_+a_-=a_+^{2k+1}=a_-^{2k+1}=0$ plus a relation of degree $2k$. It is then
clear that as a $\C [a_-,a_+]$-module, $V$ is generated by $V_{-2k}$.
(One can now show that $\g (\a ,V)=\so (V)$, but we will not need this in what
follows.)

Now assume that $k\ge 3$ and consider the semispinorial representation $W$ of
$\so (V)$. Let us first recall how $W$ is obtained. The inner product is
nondegenerate on the plane $V_0$. Let $F_0$ and $F_0'$ be complementary
isotropic lines in $V_0$ and put $F:=F_0 +\sum _{i>0} V_{2i}$ and
$F':=F_0'+\sum
_{i<0} V_{2i}$. These are complementary isotropic subspaces of $V$. One of the
semispinorial representations of $\so (V)$ can be realized on $W:=\wedge^{\ev}
F[n]$, with $n=k^2$. Now $W$ has as its lowest degree summand a line in degree
$-n$ spanned by $1\in \wedge^0 F$. The action of $\a$ on $\wedge^{\ev} F$ is as
follows: if $f\in F_0$ and $f'\in F_0'$ are such that their inner product is
$1$, then  we have $$ a\cdot (x_1\wedge\cdots \wedge x_{2l})=
-f\wedge a(f')\wedge x_1\wedge\cdots \wedge x_{2l}+
\sum _{i=1}^{2l} x_1\wedge\cdots\wedge a(x_i)\wedge\cdots \wedge x_{2l},
$$
where $a\in\a$ and $x_1,\dots ,x_{2l}\in F$.
A calculation shows that $a_{\pm}(f')=\pm ca_{\pm}(f)$, with $c$ a nonzero
constant. This enables us to determine the  $\C [a_-,a_+]$-submodule of $W$
generated by $1\in W_{-n}$: we find that in degree
$\le -n+6$ it is the span of $f,\, f\wedge a_{\pm}(f),\, f\wedge
a_{\pm}^2(f),\, a_+(f)\wedge a_-(f),\, f\wedge a_{\pm}^3(f)+a_{\pm}(f)\wedge
a_{\pm}^2(f),\, a_{\pm}(f)\wedge a_{\mp}^2(f)$. This contains the summands of
degree $\le -n+4$ of $W$, but not $W_{-n+6}$: since $n=k^2\ge 9$, $f\wedge
a_{\pm}^3(f)$ is not in this span.
So $W$ is Frobenius up to order $2$, but not up to order $3$.

\bigskip\label
Let $(\g ,h,\a )$ be a Lefschetz triple, $\h$ a Cartan subalgebra of $\g$
containing $h$ and $B$ a root basis of $R=R(\g ,\h )$ whose members are $\ge
0$ on $h$. Let $M$ be an irreducible representation of $\g$ with
$-\lambda\in\h ^*$ as lowest weight (with respect to $B$). We give a simple
criterion for $M$ to be quasi-Frobenius--Lefschetz module.

The lowest weight subspace of $M$, $M^{-\lambda}$, is a line. The
$\g$-stabilizer of this line is the standard parabolic subalgebra $\p =\p _X$
of
$\g$, where $X$ is the set of $\alpha\in B$ with $\lambda (\alpha ^{\vee})=0$.
In other words, $\p $ is the direct sum of $\h$ and the root spaces of
roots $\beta$ with $\lambda (\beta ^{\vee})\le 0$.
We endow $M$ with the grading defined by $h$ and denote the grade by a
subscript. It is clear that the depth $n$ of $M$ is
equal to $\lambda (h)$. So $M^{-\lambda}\subset M_{-n}$. The map that assigns
to
$x\in\g$ the homomorphism $x|M^{-\lambda}\in\Hom (M^{-\lambda},M)$
induces an injective map $\g /\p\to \Hom (M^{-\lambda},M/M^{-\lambda})$. In
particular, if $\p _2$  denotes the degree two part of $\p$:
$$
\p _2:=\sum _{\alpha\in R_2, \lambda (\alpha ^{\vee})=0} \g ^{\alpha},
$$
then we have an injective map $\g _2/\p _2\to \Hom (M^{-\lambda},M_{-n+2})$.
We conclude:

\proclaim{\label Lemma}
Suppose $\dim M_n=1$. Then $M$ is quasi-Frobenius if and only if $\a$
is a supplement of $\p _2$ in $\g _2$.
\endproclaim

\label If $\g$ is a Lie algebra, then call a  representation $M$ of $\g$ {\it
tautological} if $\g$ maps isomorphically onto  the Lie algebra of
infinitesimal isometries of a nondegenerate ($\pm$)-symmetric form on $M$.
So $\g$ is then a symplectic or orthogonal Lie algebra and
if $\g\not\cong \so (4)$, then $M$ is a fundamental representation of $\g$
with highest weight at an end of the Dynkin diagram.

We have not found an example of a simple Frobenius--Lefschetz module that is
not of Jordan--Lefschetz type or tautological. The following theorem says that
such an example must involve an exceptional Lie algebra.

\medskip
\proclaim{\label Theorem}
A simple Frobenius--Lefschetz module that is not a Jordan-Lef-schetz module and
whose Lie algebra is simple and of classical type, is tautological.
\endproclaim

\label
Before we begin the proof, we derive some general properties of
quasi-Frobenius--Lefschetz modules. So from now on, $M$ is a
quasi-Frobenius--Lefschetz module and we retain the notation introduced above.
Since $\p$ contains $\g _0$, we have $X\supset B_0$. We denote $X\cap B_2$ by
$B_2^{\p}$ and we let $B_2^{\a }:=B_2\setminus X=B_2-B_2^{\p}$.

For $\beta\in B_2$, we denote by $R_2(\beta )$ the set of roots in $R_2$ that
have coefficient one on $\beta$ and put
$$
\g _2(\beta):=\sum _{\gamma\in R_2(\beta )}\g ^{\gamma}.
$$
Since $R_2$ is the disjoint union of the $R_2(\beta )$'s,
$\g _2$ is the direct sum of the $\g _2(\beta)$'s.
Each $\g _2(\beta)$ is a $\g _0$-invariant subspace of $\g$. Notice that
$\p _2=\sum _{\beta\in B_2^{\p }} \g _2(\beta )$ so that
$$
\bar\a :=\sum _{\beta\in B_2^{\a }} \g _2(\beta ).
$$
is a $\g _0$-invariant supplement of $\p _2$ in $\g _2$. By
\refer{6.6}, $\a$ is the graph of a linear map $\phi :\bar\a\to\p _2$.

\proclaim{\label Lemma}
For every $\beta '\in B_2^{\p}$ there exists a $\beta\in B_2^{\a}$ such that
the map $\g _2(\beta )\to\g _2(\beta ')$ induced by $\phi$ is nonzero.
\endproclaim
\demo{Proof}
Suppose not. Then $\a$ centralizes the fundamental coweight
$p _{\beta}$
corresponding to $\beta$. We show that then $[f(e),p_{\beta}]=0$,  whenever
$f(e)$ is defined. This will imply that $\a\cup f(\a)$ is in the centralizer of
$p_{\beta}$ and thus contradict the fact that $\a\cup f(\a )$ generates $\g$ as
a Lie algebra.

Write $f(e)=\sum _k f_k$ with $[h_{\beta},f_k]=kf_k$. By homogeneity,
$(e,h,f_0)$ is then also an $\slt$-triple and by uniqueness of $f(e)$, we then
have $f(e)=f_0$.
\enddemo

Virtually all the properties that we shall derive about
quasi-Frobenius--Lefschetz modules come from the following lemma.

\proclaim
{\label Lemma}
The spaces $\bar\a$ and $\phi (\bar\a)$ are  abelian subalgebras
and $[X,\phi (Y)]=[Y,\phi (X)]$ for all $X,Y\in\bar\a$.
(In particular, $[X,\phi (Y)]=0$ if $X\in\g _2(\beta )$, $Y\in\g _2(\beta ')$
with $\beta ,\beta '\in B_2^{\a}$ distinct.)
\endproclaim
\demo{Proof} Let $\gamma$ and $\gamma '$ be distinct roots of $R_2$ that have a
$B_2^{\a}$-coefficient equal to one and let $X\in \g^{\gamma}$, $Y\in\g
^{\gamma '}$. Since $\a$ is abelian, we have
$$\align
0&=[X+\phi (X),Y+\phi (Y)]\\
 &=[X,Y]+([X,\phi (Y)]+[\phi (X),Y])+[\phi (X),\phi (Y)].
\endalign
$$
The three groups of terms belong to direct sums of root spaces that do not
intersect and so each of them must be zero.
\enddemo

\proclaim{\label Corollary}
Let $\beta\in B_2^{\a}$, $C$ its
connected component in $B-B_2^{\p}$. Then:
\roster
\item"{(i)}" $\beta$ is the unique element of
$B_2\cap C$ and the coefficient of $\beta$ in the highest root with support in
$C$ is $1$,
\item"{(ii)}" if $\beta '\in B_2^{\p}$ is connected with $C$ and $\beta _1\in
B_2^{\a}$, then $\phi $ induces a nonzero map $\g _2(\beta _1)\to \g
_2(\beta ')$ if and only if $\beta _1=\beta$,
\item"{(iii)}" if $B_2^{\p}\not=\emptyset$,
then $\beta$ is an end of $B$ or not connected with $B_0$.
\endroster
\endproclaim
\demo{Proof} (i) If $C\cap B_2$ contains an element distinct from $\beta$, then
let $\beta '$ be such an element that is not separated from $\beta$ by
another member $B_2$. Then $\beta '\in B_2^{\a}$. We then can find a $\gamma\in
R_2(\beta )$ such that $\beta '+\gamma$ is a root. But this contradicts the
fact
that $\g_2 (\beta )$ and $\g _2(\beta ')$ commute.

The fact that $\g _2(\beta )$ is abelian implies that the sum of no two
elements
of $R_2(\beta )$ is a root. So the coefficient of $\beta$ in the
highest root with support in $C$ is $1$ (see \refer{2.5}).

(ii) Let $\beta '\in B_2^{\p}$ and $\beta _1\in B_2^{\a}$ be as in the
statement. Choose $\gamma _1\in R_2(\beta _1)$ and $\gamma '\in R_2(\beta ')$
such that the map $\g ^{\gamma _1}\to \g ^{\gamma '}$ induced by $\phi$ is
nonzero and choose $\gamma\in R_2(\beta )$ such that $\gamma +\gamma '$ is a
root. So if $X^{\gamma}\in \g ^{\gamma}$ and $X^{\gamma _1}\in \g ^{\gamma _1}$
are generators, then the $\g ^{\gamma +\gamma '}$ component of
$[X_{\gamma},\phi (X_{\gamma _1})]$ is nonzero. By \refer{6.11}, this implies
that $\beta _1=\beta$.

(iii) For this assertion we prove that if $C-\{\beta\}$ is connected with
$\beta '\in B_2^{\p}$, then  $C-\{\beta\}$ is connected and separates $\beta '$
from $\beta$. Since the Dynkin diagram is a tree this amounts to showing that
the union $Y$ of connected components of
$C-\{\beta\}$ that do not separate $\beta$ and  $\beta '$ is empty.

Suppose this is not the case: $Y\not=\emptyset$. Since $\beta '$ is separated
from $Y$ by $\beta$, all the roots of $R_2(\beta ')$ will have zero coefficient
on $Y$. According to (ii) there exist
a $\gamma _1\in R_2(\beta )$ and a $\gamma '\in R_2(\beta ')$
such that $\phi (X^{\gamma _1})$ (with $X^{\gamma _1}\in\g ^{\gamma _1}$) has
nonzero $\g^{\gamma '}$-component. Choose $\gamma\in R_2(\beta )$ such that
$\gamma +\gamma '$ is a root and $\gamma-\gamma _1$ has a nonzero coefficient
on
an element of $Y$. Let $X^{\gamma}\in\g ^{\gamma}$ be a generator and  consider
the identity  $$
[X^{\gamma},\phi (X^{\gamma _1})]=[X^{\gamma _1},\phi (X^{\gamma})].
$$
The $\g^{\gamma +\gamma'}$-component of the lefthand side is clearly nonzero.
So this is also the case for the righthand side. This implies that
$-\gamma _1+\gamma +\gamma '\in R_2(\beta ')$. This is therefore a root whose
$Y$-coefficients are zero. However, these coefficients are those of
$-\gamma _1+\gamma$ and so we arrive at a contradiction.
\enddemo

\proclaim{\label Corollary}
If $B_2^{\p}=\emptyset$, then $B_2$ is a singleton, $\a =\bar\a =\g _2$ and $M$
is a Jordan--Lefschetz module.
\endproclaim
\demo{Proof} The first part of previous corollary implies that $B_2=B_2^{\a}$
is
a singleton, $\{\beta\}$, say. Then $\bar\a$ is the sum of the root spaces $\g
^{\gamma}$ with $\gamma\in R_2$. Since $\bar\a$ is abelian, the sum of two
elements of $R_2$ is never a root. Hence $R_2\cup R_0$ contains all the
positive
roots. This implies that $\bar\a=\g _2$. Since the lowest weight of $M$
is a negative multiple of the fundamental weight at the vertex labeled by
$\beta$, it is a Jordan--Lefschetz module.
\enddemo

Corollary \refer{6.12} does not exploit \refer{6.11} to the fullest. For
instance, we have:

\proclaim{\label Lemma}
Suppose that $C$ is connected with some $\beta '\in B_2^{\p}$ and assume that
$C\cup\{\beta '\}-\{\beta\}$ is of type $A_l$, $l\ge 2$. Then
$C$ is a string and $\beta$ has no greater root length then $\beta '$.
\endproclaim
\demo{Proof}  Let us number the roots of $C\cup\{\beta '\}-\{\beta\}$ in
order: $\beta ',\alpha _1,\dots ,\alpha _{l-1}$ and let $\beta$ be connected
with $\alpha _k$, $1\le k\le l-1$.

According to \refer{6.12-ii}, $\phi$ induces a nonzero map
$\phi _{\beta ',\beta}:\g _2(\beta )\to \g _2(\beta ')$ such that
$[X,\phi _{\beta ,\beta '} (Y)]$ is symmetric in $(X,Y)\in \g (\beta )\times\g
(\beta )$. The simple roots $\{\alpha _1,\dots ,\alpha _{l-1}\}$ define a Lie
subalgebra $\s\subset\g$ isomorphic to $\sli (l)$. We denote its fundamental
weights by $\varpi _i\in (\h \cap\s)^*$ ($i=1,\dots ,l-1$).
Any root $\gamma '$ of $R_2(\beta ')$
will have the form $\gamma +\alpha _1\cdots +\alpha _i$ (including the case
$i=0$: $\gamma =\gamma '$), such that $\gamma \in R_2(\beta ')$ has all its
$C$-coefficients zero. We take $\gamma '\in R_2(\beta ')$ such that $\phi
_{\beta ',\beta}$ has nonzero component on $\g ^{\gamma '}$. Then
$$
 V:=\g ^{\gamma}+\g ^{\gamma +\alpha _1}+\cdots + \g ^{\gamma +\alpha
_1+\cdots +\alpha _{l-1}}.
$$
is the $\s$-subrepresentation of $\g _2(\beta ')$ generated by $\g
^{\gamma '}$. It is irreducible with highest weight $\varpi _1$ and is
therefore
a standard representation of $\s$.

Suppose that $k\notin\{1, l-1\}$. In view of the classification of Dynkin
diagrams, $\beta$ has then the same root length as $\beta '$.

\smallskip {\it Claim 1:} The set of roots $R_2(\beta )$ is the orbit of
$\beta$
under the action of the Weyl subgroup $W$ defined by the subroot system
$\{\alpha
_1\dots ,\alpha _{l-1}\}$ of $R$.

\smallskip
Proof. Notice that $R_2(\beta )$ consists of positive roots that are linear
combinations of $\beta ,\alpha _1,\dots ,\alpha _{l-1}$ with the coefficient of
$\beta$ being $1$. Each $W$-orbit in $R_2(\beta )$
contains a root $\delta =\beta +r_1\alpha _1+\cdots +r_{l-1}\alpha _{l-1}$ in
the closed $W$-chamber opposite the fundamental one:
$\delta (\alpha _i^{\vee})\le 0$ for $i=1,\dots ,l-1$. This means
that $2r_i\le r_{i-1}+r_{i+1}$ for $i\not=k$ (where we put $r_0=r_l=0$) and
$2r_k\le r_{k-1}+r_{k+1}+1$. It is not difficult
to verify that this can only happen when all $r_i$'s are zero, i.e., when
$\delta =\beta$. So $R_2(\beta )=W\beta$.

\smallskip {\it Claim 2}: $\g (\beta)$ is as a $\s$-representation isomorphic
to  $\wedge ^kV$.
\smallskip

Proof. According to the previous claim the weights of $\g (\beta)$ with
respect to $\h$ are in a single $W$-orbit. This implies that $\g (\beta)$
is irreducible as a $\s$-representation. Since
$\beta$ defines the minus the fundamental weight $-\varpi _k$ of the root
system generated by $\{\alpha _1,\dots ,\alpha _{l-1}\}$, this representation
must be equivalent to $\wedge ^kV$.

\smallskip
We identify $\g (\beta)$ with $\wedge ^kV$ so that $\phi$ induces a nonzero
linear map $\psi :\wedge ^kV\to V$. The Lie bracket defines a bilinear
$\s$-equivariant map from $V\times\wedge ^kV$ to a representation space
of $\s$ and so we may think of it as a projection onto $\s$-subrepresentation
of $V\otimes\wedge ^kV$. This subrepresentation is nonzero since for instance
the root spaces $\g^{\gamma +\alpha _1+\cdots +\alpha _k}$ and $\g^{\beta}$
have nonzero Lie bracket. Moreover, $[\psi (x),y ]$ is symmetric in $(x,y)$. We
show that this is impossible.

We observe that $\wedge
^kV\otimes V$ decomposes into two irreducible representations: one is $\wedge
^{k+1}V$ and the other is the space $U$ spanned by the elements $x_1\otimes
x_1\wedge\cdots \wedge x_k$.
Suppose first that the image of the Lie bracket
has a nonzero projection on $\wedge ^{k+1}V$. Then
$$
\psi (x_1\wedge\cdots \wedge x_k)\wedge y_1\wedge\cdots\wedge y_k=
\psi (y_1\wedge\cdots \wedge y_k)\wedge x_1\wedge\cdots\wedge x_k
$$
for all $x_i,y_i\in V$. This identity shows that if
$\psi (y_1\wedge\cdots \wedge y_k)\not= 0$, then
each $x_i$ must be in the span of $\psi (x_1\wedge\cdots \wedge x_k),
y_1,\cdots,
y_k$. Since $\psi$ is not identically zero, we deduce (by letting $y_1,\dots
,y_k$ vary) that each $x_i$ is proportional to $\psi (x_1\wedge\cdots
\wedge x_k)$ or that $k=l-1$. The last case is excluded and the first case
implies that $k=1$, which is excluded as well. So the bracket map has image in
$U$. The $\s$-equivariant projection $\pi _U: V\otimes \wedge ^kV\to U$ is
given by
$$
\align
\pi &_U(x_0\otimes (x_1\wedge\cdots \wedge x_k))=\\
&={k\over k+1} x_0\otimes (x_1\wedge\cdots \wedge x_k)-{1\over k+1}
\sum _{i=1}^k (-1)^ix_i\otimes (x_0\wedge x_1\wedge\cdots \wedge\widehat
x_i\wedge\cdots\wedge x_k).
\endalign
$$
The same reasoning as above shows that if $\pi _U(\psi (x)\otimes y)$ is
symmetric in $(x,y)$ arguments, then $\psi =0$. This proves the first assertion
of the lemma.

\smallskip
The proof of the second assertion uses a similar argument. Suppose
that the root length of $\beta$ is greater than that of $\beta '$. In view of
the classification this can only happen when $\beta$ is connected  with $\alpha
_{l-1}$ so that $C\cup\{\beta '\}$ is of type $C_{l+1}$.

\smallskip {\it Claim 3}: $\g (\beta)$ is as a $\s$-representation isomorphic
to the space $\Sym ^2(V^*)$  of quadratic forms on $V$.
\smallskip
Proof. Note that $W$ has two orbits in $R _2(\beta )$: the orbit of the long
root $\beta$ and the orbit of the short root $\alpha _{l-1}+\beta$. The
$\s$-representation generated by $\g ^{\beta}$ has $\g ^{\beta}$ as lowest
weight space. The  lowest weight
is $-2\varpi _{l-1}$ and the corresponding representation is therefore $\Sym
^2(V^*)$. The weights of this representation are just the elements of $R
_2(\beta )$ and so the claim follows.

\smallskip
We finish the argument as before. The contraction mapping $V\otimes \Sym
^2(V^*)\to V^*$ is equivariant and surjective and its kernel $U'$ is an
irreducible representation of $\s$.

We first show that the Lie bracket cannot induce a nonzero mapping $V\otimes
\Sym ^2V^*\to V^*$. For then
$\psi : \Sym ^2(V^*)\to V$ is a nonzero linear map such that the expression
$\xi (\psi (\eta ^2))\xi$ is symmetric in $\xi,\eta\in V^*$. Since $\Sym
^2(V^*)$ is spanned by the squares, there exists a $\eta\in V^*$ with $\psi
(\eta ^2)\not= 0$. It then follows that $\eta$ is proportional to $\xi$ for
almost all $\xi$, which is absurd since $\dim V\ge 2$.

To finish the argument we now suppose that the Lie bracket induces a nonzero
mapping $V\otimes \Sym ^2V^*\to U'$. The equivariant  section of the above
contraction map assigns to $\xi\in V^*$ the symmetrization
$\sym (\bold{1}_V\otimes \xi)\in V\otimes \Sym ^2(V^*)$. So the equivariant
projection $V\otimes \Sym ^2V^*\to U'$ is given by
$v\otimes\xi ^2\mapsto v\otimes\xi ^2 -\sym (\xi (v)\bold{1}_V\otimes \xi)$.
This means that the expression $\xi (x)^2\psi (\eta ^2) -\xi (\psi (\eta
^2))\xi (x)x$
(with $\xi,\eta\in V^*$ and $x\in V$) is symmetric in $\xi$ and $\eta$. So if
$\eta (x)=0$, then $\xi (x)^2\psi (\eta ^2) -\xi (\psi (\eta ^2))\xi (x)x=0$.
By taking $\xi (x)\not= 0$, we see that $\psi (\eta ^2)$ and $x$ must be
proportional. As this is true for all $x\in V$ with $\eta (x)=0$, it follows
that $\eta (\psi (\eta ^2))=0$. If we substitute $\eta =t_1\eta _1+t_2\eta _2$,
and take the $t_1^2t_2$-coefficient, we find that $\eta _1\psi (\eta _1\eta
_2)=0$. This means that $\psi$ is identically zero, which contradicts our
assumption.
\enddemo

\proclaim{\label Proposition}
Let $(\a ,M)$ be an irreducible Lefschetz module of depth $n$ with
$\dim M_{-n}=1$. Suppose that an irreducible representation
of $\g (\a ,M)$ whose highest weight is $k\ge 1$ times that of $M$,
is Frobenius up to order $l$, with $1\le l\le k$. Then $(M ,\a )$ is
quasi-Frobenius and $M_{-n+2i}=\bar\a ^iM_{-n}$ for $i\le l$.
\endproclaim

\demo{Proof}
Let $u\in M_{-n}$  be nonzero and let $M(k)$ be the $\g (\a
,M)$-subrepresentation of  $\Sym ^k(M)$ generated by $u^k$. Then
 $M(k)$ is irreducible and has highest weight $k$ times that of $M$.
It is also a Lefschetz $\a$-module of depth $kn$ with $u^k$ spanning
$M(k)_{-kn}$. By assumption, $M(k)$ is Frobenius up to order $l$.
In particular, $\a (u^k)=\bar\a (u^k)=M(k)_{-kn +2}
$. In view of the fact that $M(k)_{-kn+2}=u^{k-1}M_{-n+2}$, it follows
that $\a M_{-n}=\bar\a M_{-n}=M_{-n+2}$, so that $M$ is quasi-Frobenius also.
Hence for $i\le l$, the map
$$
\Sym ^i(\bar\a )\otimes u^k\to \Sym ^i(M_{-n+2})u^{k-i},\quad
a_1a_2\cdots a_i\otimes u^k\mapsto a_1(u)\cdots a_i(u) u^{k-i}
$$
is an isomorphism. There is an  obvious projection of $\Sym ^k(M)_{-kn+2i}$
onto
its subspace $\Sym ^i(M_{-n+2})u^{k-i}$ which makes the above map factor as
$$
\Sym ^i(\bar\a )\otimes u^k\to M(k)_{-kn+2i}\to \Sym ^i(M_{-n+2})u^{k-i},
$$
with the first map given by the $\bar\a$-action on $M(k)$.
The image of that first map is just $\bar\a ^i(u^k)$ and hence
$\dim\, \bar\a ^i (u^k)=\dim \Sym ^i(\bar\a )=\dim \Sym ^i(\a )$.
Since $M(k)$ is
quasi-Frobenius up to order $l\ge i$, we also have $\dim M(k)_{-kn+2i}\le \dim
\Sym ^i(\a )$.  It follows that $M(k)_{-kn+2i}=\bar\a ^i(u^k)$. Taking the
projection in the summand $M_{-n+2i}u^{k-i}$, then gives that $\bar\a
^i(u)=M_{-n+2i}$.
\enddemo

\medskip\label {\it Example \refer{6.4} continued.} In this case
$\bar\a$ is the intersection of $\so (V)$ with $\Hom (V_{-2},F_0)+\Hom
(F_0',V_2)$. We may identify $\bar\a$ with $F_0\wedge V_2$ and if do so, then
its action on $W$ is given by the wedge product. From this it is immediate
that $\bar\a ^2$ acts trivially on $W$. So the previous proposition
implies that any irreducible representation of $\so (V)$ with lowest weight a
multiple of that of $W$ is not Frobenius.

\demo{Proof of \refer{6.8}}
In view of \refer{6.13}, the assumption that $(M,\a )$ is not a
Jordan--Lefschetz module implies that $B_2^{\p}$ is nonempty. According to
\refer{1.17}, $B_2$ is totally disconnected. So it follows  from \refer{6.12}
that the elements of $B_2^{\a}$ are ends of the Dynkin diagram.  Let the
numbers
$n$, $k$, $r$ and $1\le d_0<d_1<\cdots <d_r=d_{r+1}=\cdots =d_k$  have the same
meaning as in \refer{1.16}.

\smallskip
{\it Case $A_l$.} Since $B_2^{\a}$ consists of ends, we must have
$d_0=1$. Hence $M$ has lowest weight of the form $-p\varpi _1-q\varpi _l$,
with $p,q$ nonnegative integers. According to \refer{6.3}, $M$
is self-dual, so that $p=q$ and $B_2^{\a}=\{\alpha _1,\alpha _l\}$.

Notice  that for  $p=q=1$ we get the adjoint representation of $\g =\sli (V)$.
This representation is not quasi-Frobenius: its lowest degree summand is
$\sli (V)_{-2n}=\Hom (V_n,V_{-n})$ and so every element of $\a (\sli
(V)_{-2n})$
must be in $\Hom (V_{n-2},V_{-n})$, which is a proper subspace of $\sli
(V)_{-2n+2}=\Hom (V_n,V_{-n+2})\oplus \Hom (V_{n-2},V_{-n})$.

So by \refer{6.15}, $M$ cannot be quasi-Frobenius either and therefore this
case does not occur.

\smallskip
{\it Case $B_l$.} Then $n$ is even, $d_k$ is odd and $l=d_0+\cdots
+d_{k-1}+{1\over 2}(d_k-1)$. The elements of $B_2$ are in position
$d_0,d_0+d_1,\dots ,d_0+\cdots +d_{k-1}$. So in this case $d_0=1$ and
$B_2^{\a}=\{\alpha _1\}$, the simple root at the tautological vertex of $B$.

\smallskip
{\it Case $C_l$, odd parity.}
Then $n$ is odd and $l=d_0+\cdots +d_k$. The elements of $B_2$ are in position
$d_0,d_0+d_1,\dots ,d_0+\cdots +d_k=l$. The last element is the large simple
root, so by \refer{6.14} cannot belong to $B_2^{\a}$. Hence $d_0=1$ and
$B_2^{\a}=\{\alpha _1\}$, the simple root at the standard tautological vertex
of
$B$.

\smallskip
{\it Case $C_l$, even parity.}
Then  $n$ and $d_0,\dots ,d_k$ are even and $l=d_0+\cdots
+d_{k-1}+{1\over 2}d_k$. The elements of $B_2$ are in position
$d_0,d_0+d_1,\dots ,d_0+\cdots +d_{k-1}$. Neither the first nor the last root
are among them, so this case cannot occur.

\smallskip
{\it Case $D_l$, odd parity.} Then $n$ is odd, all $d_i$'s are even and
$l= d_0+\cdots +d_{k-1}+d_k$. The elements of $B_2$ are in position
$d_0,d_0+d_1,\dots ,d_0+\cdots +d_k$ and $d_k\ge 4$.
Now $\alpha _{l-d_k},\alpha _l\in B_2$, whereas $\alpha _i\in B_0$ for the
intermediate indices $i=l-d_k+1,\dots ,l-1$. Since $d_k\ge 4$, it follows from
\refer{6.14}, that we cannot have $\alpha _l\in B_2^{\a}$. So $d_0=1$ and
$B_2^{\a}=\{\alpha _1\}$, the simple root at the tautological vertex of $B$.

\smallskip
{\it Case $D_l$, even parity.} Then $n$ and $d_k$ are even and $l=
d_0+\cdots +d_{k-1}+{1\over 2}d_k$. If $d_k\ge 4$, then the elements of $B_2$
are in position $d_0,d_0+d_1,\dots ,d_0+\cdots +d_{k-1}$. So $\alpha _{l-1}$
and $\alpha _l$ do not belong to this set. It follows from \refer{6.12} that
$d_0=1$ and $B_2^{\a}=\{\alpha _1\}$, the simple root at the tautological
vertex of $B$.

Suppose now $d_k=2$. Then according to \refer{1.16}
$d_0=1$ and $d_1=\cdots =d_k=2$, in other words, $V\cong V(2k)\oplus V(2k-2)$
($k\ge 2)$. If $k=2$, then we are in the $D_4$-case with all ends belonging
to $B_2$ and the center in $B_0$. According to \refer{6.12} this can only be if
$B_2^{\a}$ is a singleton. So let us assume that $k\ge 3$; we are then in the
case of example \refer{6.4} and its continuation \refer{6.16}.
According to \refer{6.12},
$B_2^{\a}$ cannot contain both $\alpha _{l-1}$ and $\alpha _l$. Suppose it
contains one of them, say $\alpha _l$.
If $B_2^{\a}$ also contains $\alpha _1$, then
$\bar\a$ is the sum  of the root spaces corresponding to the four roots $\alpha
_1,\alpha _1+\alpha _2, \alpha _l,\alpha _l+\alpha _{l-2}$. But this
contradicts
our finding in \refer{6.4} that $\dim\a\le 2$. We also cannot have
$B_2^{\a}=\{\alpha _l\}$: if that were the case,
then $M$ has lowest weight $-p\varpi _l$ ($p>0$). The representation with
lowest weight $-\varpi _l$ is a semispinorial representation $W$ as discussed
\refer{6.4} and so this is excluded by \refer{6.16}.

\smallskip
We conclude that we are in the orthogonal
or symplectic case and that $B_2^{\a}=\{\alpha _1\}$ always. So $M$ has
lowest weight $-p\varpi _1$, with $p$ a positive integer. To finish the
argument, we must show that $p=1$. For this we invoke \refer{6.15}, with $V$
taking the r\^ole of $M$. Notice that in all cases the assumption
$B_2^{\p}\not=\emptyset$ implies that $n>2$.  The lowest degree part of $V$,
$V_{-n}$, is one dimensional and $\bar\a\subset\Hom (V_{-n},V_{-n+2})+\Hom
(V_{n-2},V_n)$. So $\bar\a V_n=V_{-n+2}$, but $\bar\a V_{-n+2}=0$, whereas
$V_{-n+4}\not= 0$ (here we use that $n>2$). So \refer{6.15} implies that
$p=1$.
\enddemo

\head
\section Appendix: a property of the orthogonal and symplectic Lie algebra's
\endhead

\noindent The purpose of this appendix is to prove the following theorem.

\proclaim{Theorem}
Let $U_i$ ($i=1,\dots ,k$) be finite dimensional complex vector spaces ($k\ge
2$) endowed with a nondegenerate form (symmetric or skew) and assume that
no $U_i$ is an inner product space of dimension $2$. If we give
$U_1\otimes\cdots\otimes U_k$  the product form
$$
\la u_1\otimes\cdots\otimes u_k,v_1\otimes\cdots\otimes v_k\ra :=
\la u_1,v_1\ra\cdots\la u_k,v_k\ra ,
$$
then every simple Lie subalgebra of  $\aut (U_1\otimes\cdots\otimes U_k)$ that
contains $\aut (U_i)$ or contains a copy of $\sli (2)$ in $\aut (U_i)$  acting
irreducibly on $U_i$ ($i=1,\dots k$) is equal to $\aut (U_1\otimes\cdots\otimes
U_k)$.
\endproclaim

Before we begin the proof we make some preliminary observations and
recall two results of Dynkin.

As before, $V(d)$ denotes the standard irreducible $\slt$-module of dimension
$d+1$, which we regard as the $d$-fold symmetric product of the tautological
representation of $\slt$.

\proclaim{\label Lemma}
The decomposition of $\gl (V(d))$ into irreducible $\slt$-submodules is
$$
\gl (V(d))=\oplus _{i=0}^{d-1} \gl ^{(i)}(V(d)),
$$
where $\gl ^{(i)}(V(d))$ is the $\slt$-submodule generated by $e^i$.
Furthermore,  $\gl ^{(0)}(V(d))$ consists of the scalars, $\gl ^{(1)}(V(d))$
can be identified with the image of $\slt$ in $\gl (V(d))$ and
$\gl ^{(\odd)}(V(d))=\aut (V(d))$.
\endproclaim
\demo{Proof}
Let $W$ be an irreducible $\slt$-submodule of $\gl (V(d))$ of dimension $m+1$.
If $T\in W$ is a highest weight vector, then $[e,T]=0$ and $[h,T]=mT$. Since
$V(d)$ is a monic $\C [e]$-module, it follows that $T$ is a polynomial in $e$
(of degree $\le d$, of course). Since $[h,e^i]=2i e^i$ it follows that
$m$ is even and that $T$ is proportional to $e^{{1\over 2}m}$. On the other
hand
is clear that for $i=0,\dots ,m$, $e^i$ is coprimitive of weight $2i$ and hence
generates an irreducible $\slt$-submodule of $\gl (V(d))$ of dimension $2i+1$.
This proves the first part of the lemma. The identity $\langle e^ix,y\rangle
=(-)^i \langle x,e^iy\rangle$ shows that $e^i\in\aut (V(d))$ if and only if $i$
is odd. \enddemo

\proclaim{\label Lemma}
In the situation of the previous lemma, we have for $i=0,\dots ,m$ that
$f^i\in \gl ^{(i)}(V(d))$. Furthermore, $u_i:=\ad _f^ie^i$ resp.\
$h_i:=[e^i,f^i]$ is a semisimple element in $\gl ^{(i)}(V))$ resp.\ $\aut
(V(d))$ which commutes with $h=[e,f]$. For $d\ge 3$ and $2\le i\le d$, we have
$h_i\notin\slt$. \endproclaim
\demo{Proof} There is an inner automorphism of $\slt$ that sends $e$ to $f$ and
so the first statement follows. For the second, regard $V(d)$ as the space of
homogeneous polynomials of degree $d$ in two variables $x,y$ and let $e$
resp.\ $f$ act as $x\partial /\partial y$ resp.\ $y\partial /\partial x$. The
calculation is then straightforward: $x^ky^l$ is an eigen vector of
$h_i$  with eigen value $c_{k,l}:=k(k-1)\cdots (k-1+d)(l+1)(l+2)\cdots (l+d)-
(k+1)(k+2)\cdots (k+d)l(l-1)\cdots (l-1+d)$. We have $c_{k,l}=-c_{k,l}$ and so
$h_i\in \aut (V(d))$. A simple verification shows that
under the given constraints, $h_i$ is not proportional to $h$ and hence not in
$\slt$. It is clear that $[h,u_i]=0$. Hence $u_i$ preserves each eigen space of
$h$ and so $u_i$ is semisimple.
\enddemo

We will also need two theorems due to Dynkin \cite{Dynkin 1952b}.

\proclaim{\label Theorem}
(Dynkin) Let $\s$ be a Lie subalgebra of $\aut (V(d))$
that contains $\sli (2)$ as a proper subalgebra. Assume $d\not= 6$. Then $\s
=\aut (V(d))$.
\endproclaim

\proclaim{\label Theorem}
(Dynkin) Let $U$ and $V$ be vector spaces with a nondegenerate
form (symmetric or skew), neither of which is an
inner product spaces of dimension $4$. Then every semisimple Lie
subalgebra $\g$ of $\aut (U\otimes V)$ that strictly contains $\aut (U)\times
\aut (V)$ is equal to $\aut (U\otimes V)$.
\endproclaim

The reason for excluding $4$-dimensional inner product spaces is that such a
space is of the form $V=W_1\otimes W_2$ with $W_i$ a
symplectic plane and $\aut (V)=\aut (W_1)\times \aut (W_2)$. In that case $\aut
(U\otimes W_1)\times\aut (W_2)$ is a semisimple Lie subalgebra of $\aut
(U\otimes V)$ that strictly contains $\aut (U)\times \aut (V)$ (at least, if
$\dim U\ge 3$). However this algebra is not simple and as we shall see the
exceptions  disappear if $\g$ is assumed to be simple.

\medskip
Let $U$ be vector space of finite dimension $\ge 2$ with a nondegenerate
$\epsilon$-symmetric form and denote its Lie algebra of infinitesimal
automorphisms by $\aut (U)$.  Let $\g _{\pm}(U)$ be the set of $x\in\sli (U)$
satisfying $\la xu,u'\ra =\pm\la u,xu'\ra$ and  let $\g _0(U)$ denote the
scalar
operators in $\gl (U)$.

\proclaim{\label Lemma} Suppose that $U$ is not an inner product space of
dimension two. Then $\gl _-(U)=\aut (U)$ and $\gl (U)=\g _-(U)\oplus \g
_0(U) \oplus \g _+(U)$  is an $\aut (U)$-invariant decomposition. The summands
are irreducible, except when $U$ is an inner product space of dimension $4$. If
$\dim U >2$, then  $[\g _+(U),\g _+(U)]=\g _ -(U)$ and for $\epsilon=\pm$,
there exists $Y\in\g _{\epsilon}(U)$ such that $Y^2\in\g _+(U)+\g _0(U)$ is not
a scalar and has nonzero trace.
\endproclaim
\demo{Proof}
The first statements are well-known. If $\dim U >2$, then
$[\g _+(U),\g _+(U)]$ is a nontrivial subspace of $\g _ -(U)$ and hence equal
to it (since $\g _ -(U)$ is irreducible). The last statement is an easy
exercise.
\enddemo

The following lemma describes the exceptional case of theorem \refer{7.3}. This
is (at least implicit) in the tables and in any case, the
proof is straightforward.

\proclaim{\label Lemma } Let $\s :=\gl ^{(1)}(V(6))+\gl ^{(5)}(V(6))$. This is
a
simple Lie subalgebra of $\gl (V(6))$ of type $G_2$ and any semisimple Lie
algebra of $\aut (V(6))$ that strictly contains $\slt$ contains $\s$.
The subspace $\g _+(V(6))$ is an irreducible representation of $\s$ (of
dimension $27$).
\endproclaim

We now treat the essential part of the case $k=2$ of the theorem. The proof is
however typical for the way we prove it in general.

\proclaim{\label Proposition}
Let $d,d'$ be positive integers and let $\g$ be a
semisimple Lie subalgebra of $\aut (V(d)\otimes V(d'))$ which contains $\slt
\times\slt$, but is not contained in $\aut (V(d))\times\aut (V(d'))$. Then $\g$
contains $\aut (V(d))\times\aut (V(d'))$.
\endproclaim
\demo{Proof} We abbreviate $V:=V(d)$  and $V':=V(d')$.

The irreducible $\slt\times\slt$-submodules of $\gl (V\otimes V')$ are
$\gl ^{(k)}(V)\otimes \gl ^{(l)}(V')$, where $0\le k\le d$ and
$0\le l\le d'$. The submodule  $\aut (V\otimes V')$ is the sum of the $\gl
^{(k)}(V)\otimes \gl ^{(l)}(V')$ with $k+l$ odd. We are given that $\g$
contains
the summands indexed by the pairs $(1,0)$, $(0,1)$
and some $(i,j)$ with $i$ and $j$ both positive and with odd sum.
Suppose $i$ is even and $j$ is odd.

We first prove that $\g\supset \aut (V)\otimes\bold{1}$.
According to \refer{7.2}, $u_j\in \gl ^{(j)}(V')$ is nonzero semisimple with
integral eigen values and so $\Tr (u^2)\not=0$.
Consider $h_i\otimes u _j^2=[e^i\otimes u_j,f^i\otimes u_j]\in\g$. This
element must have a nonzero component in some $\gl ^{(k)}(V)\otimes\bold{1}$
with $k$ odd and $\not= 1$.

If $d\not=6$, then theorem \refer{7.3} implies that $\aut
(V)\otimes\bold{1}\subset\g$.
If $d=6$, then it follows that $\g$ contains $\s\otimes\bold{1}$, with $\s$ a
Lie algebra of type $G_2$ as described in \refer{7.6}.
Since $i$ is even and positive, it follows from \refer{7.6} that $\g$ contains
$\g _+(V)\otimes\gl ^{(j)}(V')$. If  $X,Y\in \g _+(V)$, then
$[X,Y]\otimes h_j^2=[X\otimes h_j,Y\otimes h_j]\in\g$ and we find as
before that $[X,Y]\otimes\bold{1}\in\g$. The elements $[X,Y]$,
$X,Y\in\g _+(V)$, span $\g _-(V)$ by \refer{7.5} and so also in this case
$\g\supset \aut (V)\otimes\bold{1}$.

We next prove that $\bold{1}\otimes\aut (V')\subset\g$. If $d'=1$, there is
nothing to show, so suppose $d'\ge 2$.
In passing we have shown that $\g$ contains the summand $\gl
^{(k)}(V)\otimes u _j^2$ with $k$ odd. It is easily verified that for
$d'\ge 2$, $u_j^2$ is not a scalar and so $\g$ contains a summand $\gl
^{(k)}(V)\otimes \gl ^{(l)}(V')$ with $l>0$ (and necessarily even). So the same
argument as above (with $(i,j)$ replaced by $(k,l)$) shows that
$\bold{1}\otimes\aut (V')\subset\g$.
\enddemo

\proclaim{\label Proposition}
Let $d$ be a positive integer and $U$ a finite dimensional vector
space with a nondegenerate $\epsilon$-symmetric form that is not an inner
product space of dimension two. If $\g$ is a semisimple Lie subalgebra of $\aut
(U\otimes V(d))$ which contains $\aut (U)\times\slt$, but is not
contained in $\aut (U)\times\aut (V(d))$, then $\g$ contains
$\aut (U)\times\aut (V(d))$.
\endproclaim

The proof of this proposition is analogous to the proof of the
proposition preceding it (relying sometimes on \refer{7.5} instead of
\refer{7.2}) and so we omit it.

\demo{Proof of the theorem} As any $4$-dimensional inner product space is the
tensor product of two sumplectic planes, we may assume that no $U_i$ is of
that type. Then for $k=2$ the
assertion follows from the conjunction of \refer{7.7}, \refer{7.8} and theorem
\refer{7.4}. We proceed with induction on $k$ and assume $k\ge 3$.

In case $k=3$ and all factors
$U_1,U_2,U_3$ symplectic planes, then
$$
\aut (U_1\otimes U_2\otimes
U_3)=\sli (U_1)\times\sli (U_2)\times\sli (U_3)+\sli (U_1)\otimes\sli
(U_2)\otimes\sli (U_3).
$$
Since the last summand is irreducible as a
$\sli (U_1)\times\sli (U_2)\times\sli (U_3)$-module and contains a nonzero
element of $\g$, the theorem is then immediate.

Assume therefore that we are not in this situation.
First note that $\aut (U_1\otimes\cdots\otimes U_k)$ is the direct sum of the
summands $\g _{\epsilon _1}(U_1)\otimes\cdots\otimes\g _{\epsilon _k}(U_k)$
with $\epsilon _i\in \{ -,0,+\}$ with the value $-$ being taken an odd number
of
times (if $U_i$ is a symplectic plane, then $\epsilon _i$ cannot take the value
$+$). By assumption $\g$ contains a nonzero element in a summand
$\g _{\epsilon _1}(U_1)\otimes\cdots\otimes\g _{\epsilon _k}(U_k)$ with at
least
two nonzero $\epsilon _i$'s.

If $\epsilon _i=0$ for some $i$, then let $J$ denote
the set of indices $j$ with $\epsilon _j\not= 0$.  We wish to apply our
induction hypothesis to $U:=\otimes _{j\in J}U_j$. For this, we need to know
that  $\prod _{j\in J}\aut
(U_j)$ is contained in a simple  component of  $\g\cap\aut (U)$. The latter is
certainly reductive. The simple component of $\g\cap\aut (U)$  that contains
$\aut (U_j)$  ($j\in J)$ must also contain $\otimes _{l\in J}\g _{\epsilon
_l}(U_l)$ (since no nonzero element of $\g _{\epsilon _j}(U_j)$ commutes with
the elements of $\aut (U_j)$) and so the desired property holds. Our induction
hypothesis therfore applies and we conclude that $\g$ contains $\aut (U)$. Now
$U$ cannot be the tensor product of two symplectic planes $U_i\otimes U_j$, for
in that case $\aut (U)=\sli (U_i)\times \sli (U_j)$, so that $\epsilon  _i =0$
or $\epsilon _j=0$. Hence $U$ will not be a
$4$-dimensional inner product space and we may therefore apply the induction
hypothesis once more to the tensor product of $U$ and the $U_i$ with $\epsilon
_i=0$ and conclude that $\g =\aut (U_1\otimes\cdots\otimes U_k)$.

We next deal with the case when $\epsilon _i\not= 0$ for all $i$.
Suppose first that not all factors $U_i$ are symplectic planes. Since at least
one $\epsilon _i$ is $-$, we can by renumbering assume that
$\epsilon _1=-$ and $U_2$ is not a symplectic plane. We then show that
$\epsilon _3,\dots ,\epsilon _k$ can be made zero while keeping $\epsilon _1,
\epsilon _2$ nonzero, so that this takes us to the case considered above.
With the help of \refer{7.2} and \refer{7.5} we can find $X_1,X_2\in\g _-(U_1)$
with $[X_1,X_2]\not= 0$ and
$Y_j\in\g _{\epsilon}(U_i)$ ($j=2,\dots k$) such that $Y_2^2$ is not scalar,
$Y_j^2$ is not traceless for $j\ge 3$, and $X_i\otimes Y_2\otimes\cdots\otimes
Y_k\in\g$ ($i=1,2$).  Then $[X_1,X_2]\otimes Y_2^2\otimes\cdots\otimes Y_k^2$
is an element of $\g$ whose component in $\g _{-}(U_1)\otimes\g
_+(U_2)\otimes\bold{1}\otimes\cdots\otimes\bold{1}$ is nonzero.

It remains to do the case when all factors $U_i$ are symplectic planes
(and so all $\epsilon _i$'s are $-$). We then choose $X_i\in\sli (U_i)$ for
$i=1,2,3$ and $Y_j\in\sli (U_j)$ for $j=1,\dots ,k$ such that $\Tr
(Y_j^2)\not=0$ for $j\ge 4$ and $Z:=[X_1\otimes X_2\otimes X_3, Y_1\otimes
Y_2\otimes Y_3]\notin\sli (U_1)\times\sli (U_2)\times\sli (U_3)$. Then
$$
[X_1\otimes X_2\otimes X_3\otimes Y_4\otimes\cdots\otimes Y_k,
Y_1\otimes\cdots\otimes Y_k]=Z\otimes Y_4^2\cdots Y_k ^2
$$
is an element of $\g$ with a nonzero component in
$\sli (U_1)\otimes\sli (U_2)\otimes\sli
(U_3)\otimes\bold{1}\otimes\cdots\otimes\bold{1}$.
Therefore, $\g$ contains $\aut (U_1\otimes U_2\otimes U_3)$ and the induction
proceeds.
\enddemo

\enddocument

\Refs

\widestnumber\key{Verbitsky 1990}

\ref\key Beauville
\by A.~Beauville
\paper Vari\'et\'es K\"ahleriennes dont la premi\`ere classe de Chern est nulle
\jour J.\ Diff.\ Geom.
\vol 18
\pages 755--782
\yr 1983
\endref

\ref\key Bourbaki
\by N.~Bourbaki
\book Groupes et Alg\`ebres de Lie
\publ Hermann
\publaddr Paris
\yr 1975
\endref

\ref\key Chen-Ogiue
\by B.-Y.~Chen and K.~Ogiue
\paper Some characterizations of complex space forms in terms of Chern classes
\jour Q.~J.~of Math.
\vol 26\yr 1975
\pages 459--464
\endref

\ref\key Deligne 1968
\by P.~Deligne
\paper Th\'eor\`eme de Lefschetz et crit\`eres de d\'eg\'en\'erescence de
suites spectrales
\jour Inst\. Hautes \'Etudes Sci\. Publ\. Math\.
\vol 35
\yr 1968
\pages 259--126
\endref

\ref\key Deligne 1979
\by P.~Deligne
\paper Vari\'et\'es de Shimura
\inbook Automorphic Forms, Representations, and L-functions, Part 2
\bookinfo Proc.~Symp.~ in Pure Math.
\vol 33
\eds A.~Borel and W.~Casselman
\pages 247--290
\yr 1979
\publ Amer.~Math.~Soc.
\publaddr Providence, RI
\endref

\ref\key Dynkin 1952a
\by E.B.~Dynkin
\paper Semisimple subalgebras of semisimple Lie algebras
\inbook Amer.\ Math.\ Soc.\ Transl.
\bookinfo Series 2
\vol 6
\pages 111--244
\publ Amer.~Math.~Soc.
\publaddr Providence, RI
\endref

\ref\key Dynkin 1952b
\by E.B.~Dynkin
\paper The maximal subgroups of the classical groups
\inbook Amer.\ Math.\ Soc.\ Transl.
\bookinfo Series 2
\vol 6
\pages 245--378
\publ Amer.~Math.~Soc.
\publaddr Providence, RI
\endref

\ref\key Lange-Birk
\by H.~Lange and Ch.~Birkenhage
\book Complex Abelian Varieties
\publ Springer
\publaddr Berlin and New York
\yr 1992
\endref

\ref\key LiE
\by M.A.A.~van Leeuwen, A.M.~Cohen and B.~Lisser
\book Lie, a Package for Lie Group Computations
\publ Computer Algebra Nederland
\publaddr Amsterdam
\yr 1992
\endref

\ref\key Moonen-Zah
\by B.J.J.~Moonen and Yu.G.~Zahrin
\paper Hodge classes and tate classes on simple abelian fourfolds
\jour Duke Math J.
\vol 77
\yr 1995
\pages 553--581
\endref


\ref\key Okonek-vdV
\by Ch.~Okonek and A.~van de Ven
\paper Cubic forms and complex $3$-folds
\jour l'Ens.~Math.
\vol 41\yr 1995
\pages 297--333
\endref


\ref\key Saito
\by M.~Saito
\paper Mixed Hodge modules
\jour Publ\. Res\. Math\. Sci\.
\vol 26\yr 1990
\pages 221--333
\endref

\ref\key Salamon
\by S.M.~Salamon
\paper On the cohomology of K\"ahler and hyper-k\"ahler manifolds
\jour Topology
\vol 35
\yr 1996
\pages 137--155
\endref

\ref\key Shimura
\by G.~Shimura
\paper On analytic families of polarized abelian varieties and automorphic
functions
\jour Ann.~of Math.
\vol 78
\yr 1963
\pages 149--192
\endref

\ref\key Springer
\by T.A.~Springer
\book Jordan algebras and algebraic groups
\publ Springer
\publaddr Berlin and New York
\yr 1973
\endref

\ref\key Springer-St
\by T.A.~Springer and  R.~Steinberg
\paper Conjugacy classes
\inbook Seminar on Algebraic Groups and Related Finite Groups
\eds A.~ Borel et al.
\bookinfo Lecture Notes in Mathematics
\vol 131
\publ Springer
\publaddr Berlin and New York
\yr 1970
\endref


\ref\key Verbitsky 1990
\by M.~Verbitsky
\paper On the action of a Lie algebra $\so (5)$ on the cohomology of a
hyperk\"ahler manifold
\jour Func. An. and Appl.
\vol 24
\pages 70--71
\yr 1990
\endref

\ref\key Verbitsky 1995
\by M.~Verbitsky
\paper Cohomology of compact hyperk\"ahler manifolds
\paperinfo alg.~geom.~eprint 9501001
\endref

\ref\key Zucker
\by S.~Zucker
\paper Locally homogeneous variations of Hodge structure
\jour l'Ens. Math.
\vol 27
\yr 1981
\pages 243--276
\endref


\endRefs
\bye